\documentclass[preprint,11pt]{elsarticle}
\usepackage[bookmarks=true]{hyperref}  

\usepackage{url}	

\usepackage{ifpdf}

\usepackage{listings}
\usepackage{color}
\usepackage{amssymb,amsmath}
\usepackage{amsthm}

\usepackage{calc}
\usepackage[margin=1.in]{geometry}

\newcommand{\bse}{\begin{subequations}}
\newcommand{\ese}{\end{subequations}}

\def\ba#1\ea{\begin{align}#1\end{align}}

\def\bas#1\eas{\begin{align*}#1\end{align*}}

\def\bat#1\eat{\begin{alignat}{3}#1\end{alignat}}

\def\bats#1\eats{\begin{alignat*}{3}#1\end{alignat*}}

\newcommand{\Real}{\mathbb{R}}
\newcommand{\Integer}{\mathbb{J}}

\newcommand{\Kron}[1]{\delta[#1]}
\newcommand{\Dirac}[1]{\delta(#1)}

\newcommand{\imag}{\iota}

\newcommand{\be}{\begin{equation}}  
\newcommand{\ee}{\end{equation}}

\newcommand{\bes}{\begin{equation*}}  
\newcommand{\ees}{\end{equation*}}

\newcommand{\bi}{\begin{itemize}}
\newcommand{\ei}{\end{itemize}}

\newcommand{\bn}{\begin{enumerate}}
\newcommand{\en}{\end{enumerate}}

\newcommand{\R}{\bar{J}}
\newcommand{\RJ}{\bar{J}}
\newcommand{\RL}{\bar{L}}

\newcommand{\kx}{k^{(x)}}
\newcommand{\ky}{k^{(y)}}
\newcommand{\bvec}[1]{{\bf #1}}
\newcommand{\kv}{{\bvec{k}}}
\newcommand{\jv}{{\bvec{j}}}
\newcommand{\nv}{{\bvec{n}}}

\newcommand{\xv}{{\bvec{x}}}

\newcommand{\fk}{f^{(\kv)}}
\newcommand{\fka}{f^{(\kv_1)}}
\newcommand{\fkb}{f^{(\kv_2)}}

\newcommand{\ko}{0}

\newcommand{\dx}{\Delta x}

\newcommand{\dk}{\Delta k}
\newcommand{\dkx}{\Delta \kx}
\newcommand{\dky}{\Delta \ky}

\newcommand{\p}{\partial}

\newcommand\del[2]{\delta_{#1}^{#2}}
\newcommand\om[2]{\Omega_{#1}^{#2}}
\newcommand\al[2]{\alpha_{#1}^{#2}}

\newcommand{\ak}{a_{\kv}}
\newcommand{\nk}{n_{\kv}}
\newcommand{\nkj}{n_{{\jv}}}
\newcommand{\wk}{w_{\kv}}
\newcommand{\akb}{\bar{a}_{\kv}}

\newcommand*\no[1]{\tikz[baseline=(char.base)]{
        \node[fill=pink,shape=circle,draw, inner sep=.1pt] (char){\tiny$\displaystyle{#1}$};}}

\usepackage{tikz}

\begin{document}


\begin{frontmatter}
\title{A New Approach to Direct Discretization of Wave Kinetic Equations with Application to a Nonlinear Schr{\"o}dinger System in 2D}

\author[rpi]{J. W. Banks\corref{cor1}\fnref{simonsThanks}}
\ead{banksj3@rpi.edu}

\author[nyu]{Jalal Shatah\fnref{simonsThanks}}
\ead{js14@nyu.edu}

\address[rpi]{Department of Mathematical Sciences, Rensselaer Polytechnic Institute, 110 8th Street, Troy, New York 12180, USA}

\address[nyu]{Courant Institute of Mathematical Sciences, New York University, 251 Mercer Street, New York, New York 10012, USA}

\cortext[cor1]{Corresponding author.}

\fntext[simonsThanks]{The authors are grateful for support from the Simons Collaboration on Wave Turbulence (Awards No. 617006 and No. 651459).}

\begin{abstract}
Wave Kinetic Equations (WKEs) are often used to describe the evolution of ensemble averaged wave amplitudes for nonlinear wave systems. In the present manuscript we describe a new approach to direct numerical simulation of solutions to WKEs. This new method relies on a piecewise polynomial approximation of the resonant manifold, followed by numerical quadrature of the collision integral. The approach is general in nature, and is discussed in detail here for a particular nonlinear Schr{\"o}dinger model in 2 spatial dimensions. Detailed convergence studies demonstrate 2nd-order accuracy for model collision integrals, and self-convergence studies for the WKE show near 2nd-order rates. Furthermore, comparison of the WKE approximation to ensemble averages of the NLS illustrate the efficacy of the method and the validity of the WKE, for both isotropic and an-isotropic solutions. 
\end{abstract}


\end{frontmatter}


\section{Introduction}
\label{sec:intro}

Nonlinear dispersive equations are ubiquitous in the mathematical modeling of physical systems where both wave dispersion and nonlinearity play an essential role. A general framework is provided by \emph{weakly nonlinear dispersive equations} of the form
\begin{equation}
\partial_t\mathbf{u} = \mathcal{L}(\mathbf{u}) + q(\mathbf{u}),
\label{eq:governingGeneral}
\end{equation}
where $\mathbf{u}(\mathbf{x},t)$ is a state vector depending on the spatial coordinate $\mathbf{x}\in\mathbb{R}^d$, with $d$ the spatial dimension, $t$ the time, $\mathcal{L}(\cdot)$ a linear differential operator describing dispersion, and $q(\cdot)$ a nonlinear operator of small amplitude. Such models arise naturally in fluids~\cite{Ablowitz_2011,whitham74,Lannes2013}, plasmas~\cite{Galtier2022}, and nonlinear optical media~\cite{Agrawal2019,Drazin1989}, and are used in real-life situations such as describing ocean waves and interactions with wind~\cite{komen1994,janssen2004,janssen2007}.

A key question concerns the \emph{long-time behavior} of solutions to~\eqref{eq:governingGeneral}. In realistic situations, the initial data is often only partially known or even completely unknown, owing to experimental limitations, intrinsic fluctuations, or stochastic forcing. In such cases, a purely deterministic description is of limited predictive value, and it becomes natural to adopt a \emph{statistical} viewpoint, focusing on ensemble averages and spectral quantities as the primary descriptors of the dynamics.

From this perspective, one seeks to derive \emph{effective statistical equations} that capture the evolution of the probability distribution or spectrum of solutions. This statistical reduction from microscopic dynamics to macroscopic evolution mirrors the celebrated derivation of the Boltzmann equation in kinetic gas theory~\cite{bodineau2024,boltzman1872,deng2025}, where molecular chaos assumptions yield a closed kinetic description for particle distributions. For wave systems, the analogous outcome is the \emph{Wave Kinetic Equation} (WKE), which governs the slow evolution of the wave spectrum and emerges naturally in the framework of \emph{Wave Turbulence Theory} (WTT)~\cite{Nazarenko2011,Zakharov1992}.

In WTT, one recasts~\eqref{eq:governingGeneral} in Fourier space and, for Hamiltonian systems, in canonical variables. This representation makes the structure of resonant wave interactions explicit. These resonant wave interactions are the dominant contributors to long-time spectral transfer. The WKE plays for waves, the same role the Boltzmann equation plays for particles: it bridges the gap between the underlying deterministic (Hamiltonian) dynamics and the emergent statistical behavior.

A particularly striking prediction of WTT is the existence of exact stationary power-law solutions of the WKE, known as the \emph{Kolmogorov-Zakharov} (KZ) spectra~\cite{Zakharov1992,Nazarenko2011}. These spectra describe scale-invariant cascades of conserved quantities such as energy or wave action across Fourier space, analogous to Kolmogorov's $k^{-5/3}$ law in hydrodynamic turbulence. The direction of these cascades, direct or inverse, depends on the conserved quantity and the sign of the dispersion. In applications, KZ spectra have been observed or predicted in systems ranging from the nonlinear Schr{\"o}dinger equation~\cite{buckmaster2021,banks22_WKE,deng2023} to surface gravity and capillary waves~\cite{pan2017}.

Although the derivation of the WKE is now classical, solving it, especially in its full time-dependent form, remains computationally challenging. It is a nonlinear integro-differential equation with potentially singular collision kernels supported on resonant manifolds determined by conservation of energy and momentum. These features make numerical approximation difficult, particularly in capturing the geometry of resonances and accurately evaluating the multidimensional integrals. Several numerical strategies have been proposed to address this. In the pioneering work of Webb and Resio~\cite{webb1978,Resio1991}, assumed symmetries are leveraged to generate special polar grids where resonance conditions are satisfied at grid points. In the work of Krstulovic et al.~\cite{krstulovic2025}, the authors develop a method based on an analytic parametrization of the resonant manifold, onto which the wave spectrum is interpolated before applying quadrature techniques to evaluate the collision integral. Walton and Tran~\cite{walton2023,walton2025} consider isotropic formulations and recast the WKE as a conservation law, which exposes forward and backward collision operators and leads to a stiff system treated using implicit time-stepping. The approach due to Engquist, Tornberg, and Tsai~\cite{engquist2005}, employs a level-set representation of the resonant manifold via signed distance functions to regularize and localize the collision integral.

In this manuscript, we introduce a new numerical method for computing the time evolution of the WKE. Our approach avoids the need for analytic parametrization and is designed to be both flexible and efficient across a wide range of physical systems. It enables accurate and systematic benchmarking against theoretical predictions, including stationary solutions, as well as direct numerical simulations of the underlying wave equations. This work represents an advancement in the computational study of wave turbulence and provides a versatile tool for probing the statistical dynamics of nonlinear dispersive systems far from equilibrium.

The remainder of this manuscript proceeds as follows. Section~\ref{sec:governing} presents a model governing dynamical equation in the form of a particular nonlinear Schr{\"o}dinger equation, as well as the corresponding wave kinetic equation which is the primary subject of the work. The core numerical approach which is used to discretize the right-hand-side of the WKE is then described for simple models in Section~\ref{sec:method}, including for 1D in Section~\ref{sec:1Dmodel} and 2D in Section~\ref{sec:2D}. In both 1D an 2D the numerical approximations are compared against known exact solutions, and the approach is shown to be second-order accurate. Section~\ref{sec:verification} then applies the new technique to the time-dependent WKE. First the full method-of-lines discretization is presented in Section~\ref{sec:WKEManipulations}, and baseline test cases as well as approximate solutions for both isotropic and anisotropic cases are given in Section~\ref{sec:ICs}. Self-convergence studies for these cases are presented in Section~\ref{sec:selfConvergence}, followed by a comparison to ensembles of simulation results for the original NLS equation in Section~\ref{sec:NLSComparison}. Concluding remarks are given in Section~\ref{sec:conclusions}. Finally for completeness, the derivation of the WKE for this particular NLS system is presented in~\ref{sec:WKEDerive}.

\section{Governing Equations}
\label{sec:governing}
The focus of the present manuscript is on the development of new numerical methods for a broad class of wave kinetic equations. Nevertheless, it is useful to present the approach in a concrete setting, and so let us consider the equations of motion
\ba
  \imag\p_t u = (I-\Delta) u + u^2+2u\bar{u}
  \label{eq:NLS2D}
\ea
for $\xv\in\Real^2$, where $\imag=\sqrt{-1}$ is the imaginary unit. Note that the appearance of $\imag$ on the left-hand-side of the equation is simply for historical reasons. Note also that in~\eqref{eq:NLS2D}, the strength of the nonlinearity is defined by the size of the solution rather than a preceding multiplicative constant. This is simply a choice of convenience for the derivation of the corresponding WKE 
\ba
  \p_{\tau}\nk = &\frac{1}{\pi}\int\int
    \delta(\kv_1+\kv_2-\kv)\left(n_{\kv_1}n_{\kv_2}-\nk n_{\kv_1}-\nk n_{\kv_2}\right)
      \delta(\omega_{\kv_1}+\omega_{\kv_2}-\omega_{\kv})\nonumber\\
    & \qquad +\delta(\kv_1-\kv_2-\kv)\left(2n_{\kv_1}n_{\kv_2}-2\nk n_{\kv_2}+2\nk n_{\kv_1}\right)
      \delta(\omega_{\kv_1}-\omega_{\kv_2}-\omega_{\kv})\,d\kv_2\, d\kv_1,
  \label{eq:NLSWKE}
\ea
where $\nk(\tau)=n(\kv,\tau)$ is an average wave amplitude for wave vector $\kv\in\Real^2$, $\tau$ is the so-called kinetic time\footnote{See~\ref{sec:WKEDerive}) for the exact definition of the wave amplitude and $\tau$ with respect to $t$, both of which involve a box size, $L$, which is tending to infinity.}, the linear dispersion relation defines $\omega_\kv=1+\|\kv^2\|$, the Dirac delta function is $\delta(\cdot)$, and $\kv_1$ and $\kv_2$ are dummy vector variables of integration.  A derivation of~\eqref{eq:NLSWKE} from~\eqref{eq:NLS2D} is presented in~\ref{sec:WKEDerive}. Equation~\eqref{eq:NLSWKE} is the targeted governing equation for the newly developed discretizations that are the subject of this manuscript.

The focus of the present manuscript is new numerical methods to approximate the solution to WKEs. Therefore, the model given by~\eqref{eq:NLS2D} and WKE~\eqref{eq:NLSWKE} is chosen not for its connection to any particular physical system, but rather due to its simplicity and the resulting clarity in describing the numerical technique. This has the advantage that here we present a self-contained manuscript encompassing governing equation~\eqref{eq:NLS2D}, derivation of the corresponding WKE in~\ref{sec:WKEDerive}, description of new numerical method in Section~\ref{sec:method}, and finally results in Section~\ref{sec:verification}. These results include a comparison of the WKE solution to ensemble averages of solutions to~\eqref{eq:NLS2D}. In addition, future work can use the model~\eqref{eq:NLS2D} and~\eqref{eq:NLSWKE}, in combination with our results as a benchmark for comparison of new methods and techniques.

\section{Development of the Numerical Approach}
\label{sec:method}
Before confronting the WKE directly, we describe the generic approach in an even more constrained environment. Using simple models, we are able to more clearly discuss the approach, and derive exact solutions for which we perform grid convergence studies. These models are constructed to incorporate the various levels of complexity existent in general WKEs, and they can serve as good tests for comparison of future numerical methods. 

\subsection{1D Model of RHS}
\label{sec:1Dmodel}
To describe the basic methodology, consider the following 1D model for the right-hand side of a WKE
\begin{align}
  C(k) = \int_{-\infty}^{\infty} f(\xi;k)\delta(g(\xi;k))\, d\xi,
  \label{eq:model1D}
\end{align}
where $k$ is essentially a parameter in the RHS, and to simplify the presentation the dummy variable of integration has been changed to $\xi$ rather than $k_1$ which was used in~\eqref{eq:NLSWKE}. This choice of notation will be extended throughout Section~\ref{sec:method}, where the wave-number space is indicated with Roman $k$, while the integration space with a Greek $\xi$. The basic task is to approximately evaluate $C(k)$ with known asymptotic accuracy knowing only how to evaluate $f(\xi;k)$, $g(\xi;k)$, and $g_{\xi}(\xi;k)$. In practice, the infinite domain will be truncated to a finite one, and so it is useful to take $f$ to be compactly supported (or nearly so) on a finite interval $k\in[k_a,k_b]$ so that boundary effects are negligible. Therefore take $k_j=k_a+j\dk$ for $j=0,1,\ldots,N$, $\dk=(k_b-k_a)/N$ to be a uniform grid and consider the approximation of $C(k_j)$. 

To proceed, $g(\xi;k_j)$ is approximated by a piecewise polynomial, $\tilde{g}(\xi,k_j) \approx g(\xi;k_j)$, whose zero level set is taken as an approximate resonant manifold. Using this manifold, integrations are then performed either exactly or approximately, which ultimately defines the scheme. To begin, consider $\tilde{g}$ to be a piecewise linear interpolant, i.e. the piecewise linear function fitting the discretely evaluated $g_{J,j}=g(\xi_J,k_j)$ for $J=0,1,...,N$. Note the notation used here will be applied more generally throughout where lower case refers to the grid index for wave number space $k$, while upper case refers to the grid index for the integration variable $\xi$. Now assume that $\bar{\xi}$ is a root of $\tilde{g}$, i.e. the zero level set of $\tilde{g}$, with $\bar{\xi}\in[\xi_{\R},\xi_{\R+1}]$. Note that such intervals or cells can be identified by the feature that $g(\xi_{\R},k_j)$ and $g(\xi_{\R+1},k_j)$ have different signs. In this setting, the interpolant of $g$ for $\xi\in[\xi_{\R},\xi_{\R+1}]$ is
\begin{align}
  L_g(\xi;k_j) \equiv g_{\R,j}+\frac{g_{\R+1,j}-g_{\R,j}}{\dk}\left(\xi-\xi_{\R}\right)\label{eq:gApprox1D}.
\end{align}
%
Using this approximation in~\eqref{eq:model1D} then obtain
\bse
\label{eq:1DBasic}
\begin{align}
  \int_{\xi_{\R}}^{\xi_{\R+1}} f(\xi;k_j)\delta(g(\xi;k_j))\, d\xi 
    &\approx \int_{\xi_{\R}}^{\xi_{\R+1}} f(\xi;k_j)\delta(\tilde{g}(\xi;k_j))\, d\xi.\\
  \tilde{g}(\xi;k_j) & = L_g(\xi;k_j).
\end{align}
\ese
One can evaluate the approximation suggested in~\eqref{eq:1DBasic} to give
\bse
\label{eq:scheme1D_p1}
\begin{align}
  \int_{\xi_{\R}}^{\xi_{\R+1}} f(\xi;k_j)\delta(\tilde{g}(\xi;k_j))\, d\xi 
    & = \frac{f(\bar{\xi};k_j)}{\tilde{g}_{\xi}(\bar{\xi};k_j)}\\
 \tilde{g}_{\xi}(\bar{\xi};k_j) 
   & = \frac{{g}_{\R+1,j}-{g}_{\R,j}}{\dk}.
\end{align}
\ese
This is, in principle, a complete approach for the model as presented in~\eqref{eq:model1D}, but it is deficient for a number of reasons. Firstly in the larger context, $f$ will involve discrete grid quantities which eliminates the ability to evaluate $f$ anywhere at will. As a result, it becomes important to consider interpolation from grid quantities. Secondly, the scheme~\eqref{eq:scheme1D_p1} has low accuracy because $\tilde{g}_{\xi}(\bar{\xi};k_j) -g_{\xi}(\bar{\xi};k_j)=O(\dk)$, i.e. is only 1st order accurate, unless $\bar{\xi}$ is a midpoint, which is generally not true. To simultaneously address both of these concerns, piecewise linear interpolation is applied to both $f$ and $g_{\xi}$ to give
\bse
\begin{align}
  f(\xi;k_j) & \approx L_f(\xi;k_j) \equiv f_{\R,j}+\frac{f_{\R+1,j}-f_{\R,j}}{\dk}\left(\xi-\xi_{\R}\right)\label{eq:fApprox1D}\\
  g_{\xi}(\xi;k_j) & \approx L_{g_{\xi}}(\xi;k_j) \equiv {g_{\xi}}_{\R,j}+\frac{{g_{\xi}}_{\R+1,j}-{g_{\xi}}_{\R,j}}{\xi_{\R+1}-\xi_{\R}}\left(\xi-\xi_{\R}\right). \label{eq:gDerivApprox1D}
 \end{align}
\ese
Here the basic assumption is that $g$ is evaluated on the grid to define a piecewise linear interpolant $L_{g}$ which then defines $\bar{\xi}$. Similarly $f$ and $g_{\xi}$ are evaluated on the grid to define piecewise linear interpolants $L_{f}$ and $L_{g_{\xi}}$. Using these three then gives the approximation
%
\begin{align}
  \int_{\xi_{\R}}^{\xi_{\R+1}} f(\xi;k_j)\delta(g(\xi;k_j))\, d\xi 
    &\approx \frac{L_f(\bar{\xi};k_j)}{L_{g_{\xi}}(\bar{\xi};k_j)}\label{eq:scheme1D_p2}
\end{align}
%
which is second-order accurate, as will be demonstrated below.

\subsubsection{Numerical Results for Simple 1D Model}
\label{sec:results2DSimple}
\begin{figure}[h]
  \begin{center}
    \begin{tikzpicture}[scale=1]
        \draw (0cm,0cm) node[anchor=south]{\includegraphics[width=6cm]{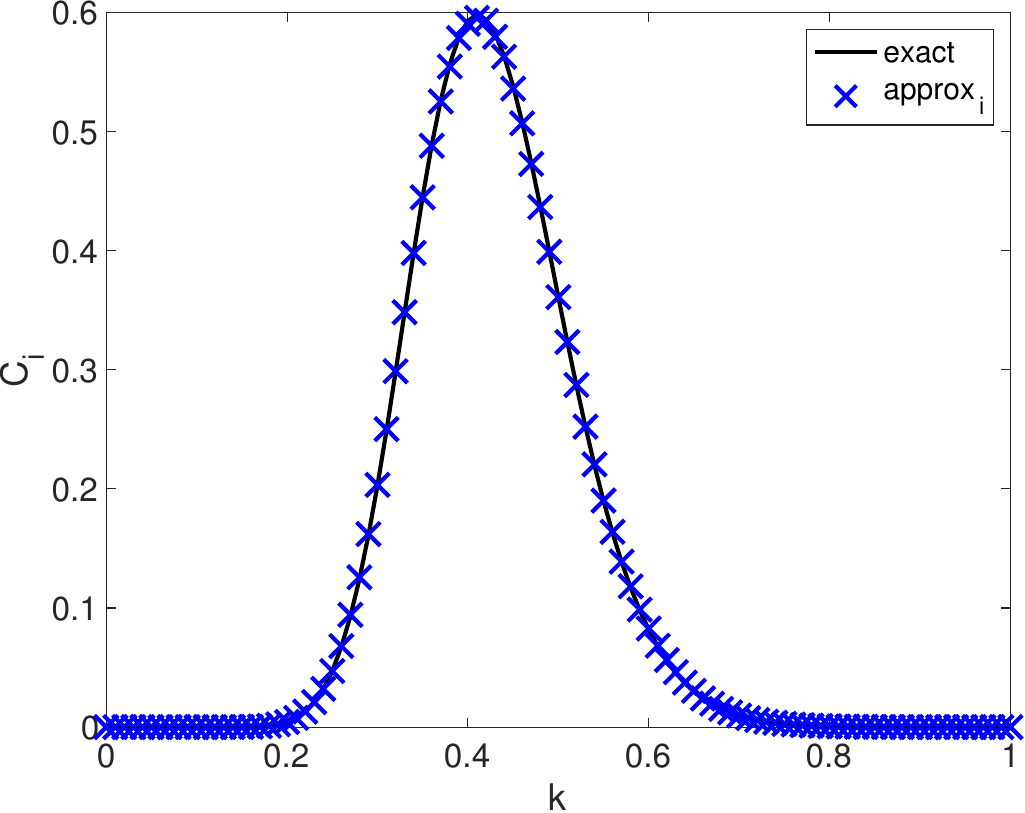}};
        \draw (7cm,0cm) node[anchor=south]{\includegraphics[width=6cm]{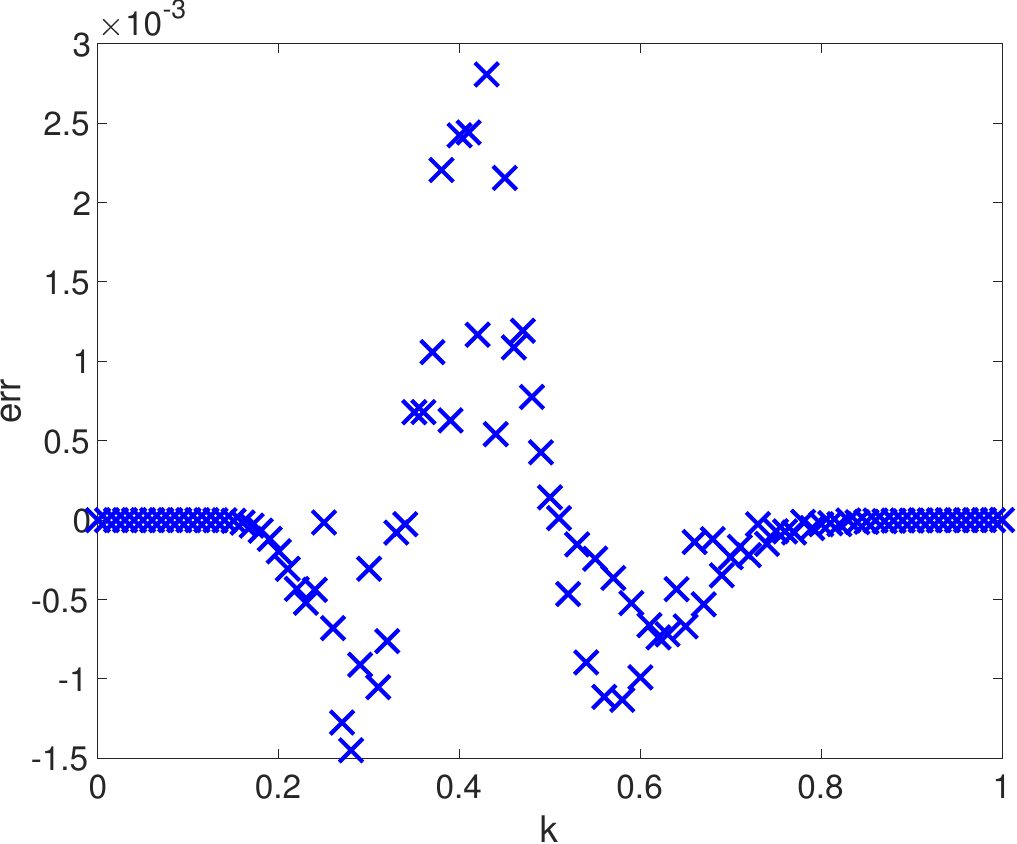}};
      \end{tikzpicture}
  \caption{Plot of exact soln and approximation (left), and the error in the approximation (right) for Test i.} 
  \label{fig:T1_soln_meth2}
  \end{center}
  \end{figure}

\begin{figure}[h]
  \begin{center}
    \begin{tikzpicture}[scale=1]
        \draw (0cm,0cm) node[anchor=south]{\includegraphics[width=6cm]{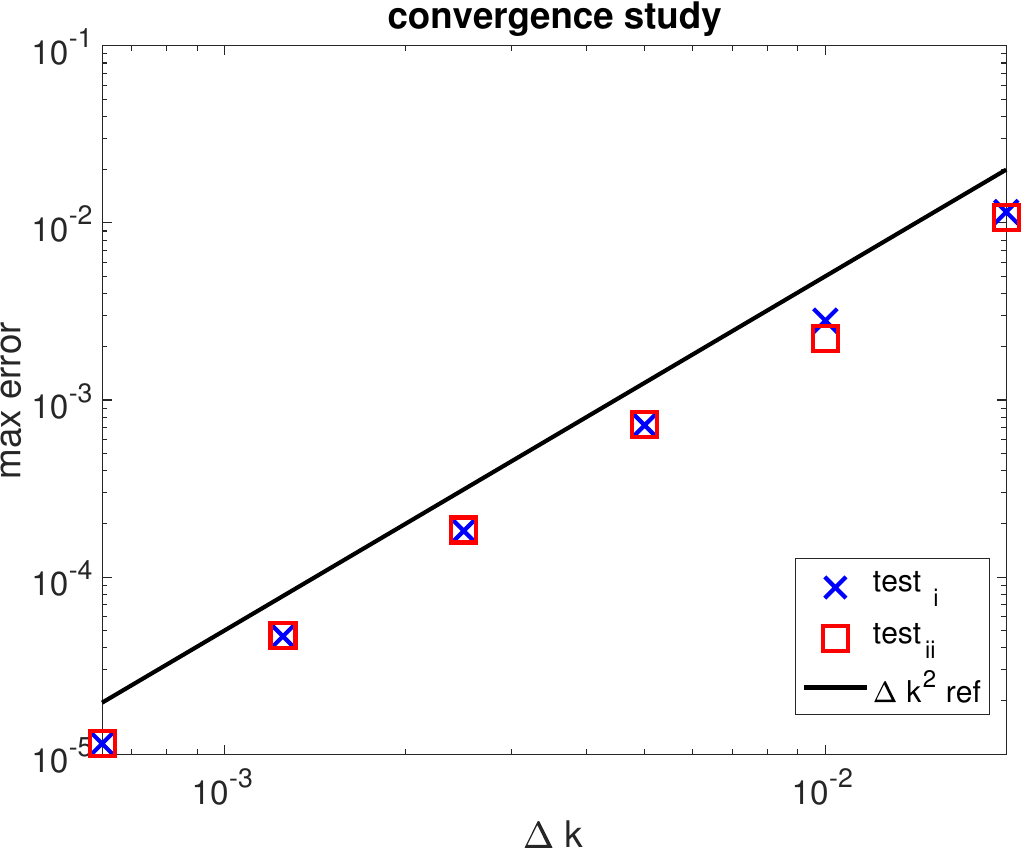}};
      \end{tikzpicture}
  \caption{Max-norm convergence results showing error vs. grid size for Tests i and ii. Clear 2nd order convergence is demonstrated for both, with no appreciable differences observed.} 
  \label{fig:1D_conv}
  \end{center}  
  \end{figure}

To illustrate the behavior of the scheme~\eqref{eq:scheme1D_p2}, consider the concrete test defined using
\bse
\begin{align}
  f(\xi,k) & = e^{-200(\xi-\frac{3}{4})^2}\\
  g(\xi,k) & = \xi^3-k\\
  g_{\xi}(\xi,k) & = 3\xi^2
\end{align}
for which the model problem can be written
\begin{align}
  C_i(k) = \int_{-\infty}^{\infty} e^{-200(\xi-\frac{3}{4})^2}\delta\left(\xi^3-k\right)\, d\xi,
  \label{eq:modelI}
\end{align}
where the ``i'' subscript indicates model problem i. Making the change of variables $w=\xi^3-k$ gives the exact solution
\begin{align}
  C_i(k) = \frac{1}{3\sqrt[3]{k^2}}e^{-200(\sqrt[3]{k}-\frac{3}{4})^2}.
\end{align}
\ese
Clearly this solution is not defined for $k=0$, and so to avoid this singularity take the domain to be $k_a=10^{-16}$, $k_b=1$. Note that in this case the resonant set is always contained in the computational domain. Figure~\ref{fig:T1_soln_meth2}  shows computed results with $N=100$ discretization points. The figure includes the exact solution, computed approximation, and errors for test a using the 2nd order method presented in~\eqref{eq:scheme1D_p2}. Furthermore, Figure~\ref{fig:1D_conv} presents results for a convergence study using the max-norm for problem i (and problem ii to be introduced in Section~\ref{sec:extended1D}) which demonstrates the expected 2nd order accurate convergence.
  
\subsubsection{An Extended 1D Model}
\label{sec:extended1D}
The structure of WKEs presents certain challenges and opportunities that are not obvious from the form of the model in Section~\ref{sec:1Dmodel}. Specifically there are often portions of $f$ whose functional dependence is not directly $\xi$, but instead a shifted version $\xi-\hat{\xi}$. This shift comes about, for example, upon performing the $\kv_1$ integration for the first term from whence dependence on $\kv_1$ will become dependence on ${\kv-\kv_2}$. This of course would not be an issue if $f$ were considered a closed-form expression which could be evaluated anywhere, but here there are portions taken as given at grid points. Note that for more general WKEs there are often portions of $f$ that can be evaluated directly (i.e. there is a closed form expression defining them). In the present model~\eqref{eq:NLSWKE}, this portion is simply constant, and so it will be disregarded here, although its inclusion is a triviality. To model this situation in a close analog to the WKE~\eqref{eq:NLSWKE}, one might consider $f$ as the following product
\begin{align}
  f(\xi;k) = \fk(k)\fkb(\xi)\fka(\xi-\hat{\xi}).
  \label{eq:model1DExtended}
\end{align}
Here $\fk$ is an analog to $n_{\kv}$, $\fkb$ is an analog to $n_{\kv_2}$, and $\fka$ is an analog to $n_{\kv_1}$. Note that in general the $\hat{\xi}$ shift could easily take the interpolation stencil outside the computational domain where data is assumed known. In this case we will rely on the idea that $f$ is essentially compactly supported in the computational domain and so the interpolation point is taken to be the closest point on the boundary of the computational domain. From a computational perspective, the $\fk(k)$ term is not algorithmically interesting since it simply moves outside the integral and scales the result. Thus for the present test case take $\fk=1$. A model problem nearly identical to model i from Section~\ref{sec:results2DSimple} can be constructed by taking 
\bse
\begin{align}
  \fk(k) & = 1\\
  \fkb(\xi) & = e^{-100(\xi-\frac{3}{4})^2}\\
  \fka(\eta) & = e^{-100(\tilde{\xi}-\frac{3}{8})^2}\\
  \tilde{\xi} & = \xi-\hat{\xi}\\
  g(\xi,x) & = \xi^3-k\\
  g_{\xi}(\xi,x) & = 3\xi^2.
\end{align}
\ese
By setting $\hat{\xi}=\frac{3}{8}$ one obtains the model 
\bse
\begin{align}
  C_{ii}(k) 
    & = \int_{-\infty}^{\infty} e^{-100(\xi-\frac{3}{4})^2}e^{-100(\tilde{\xi}-\frac{3}{8})^2}\delta\left(\xi^3-k\right)\, d\xi,\\
    & = \int_{-\infty}^{\infty} e^{-200(\xi-\frac{3}{4})^2}\delta\left(\xi^3-k\right)\, d\xi,
\end{align}
\ese
which is identical to~\eqref{eq:modelI} whose exact solution is $C_{i}(k)$.

\begin{figure}[h]
  \begin{center}
    \begin{tikzpicture}[scale=1]
        \draw (0cm,0cm) node[anchor=south]{\includegraphics[width=6cm]{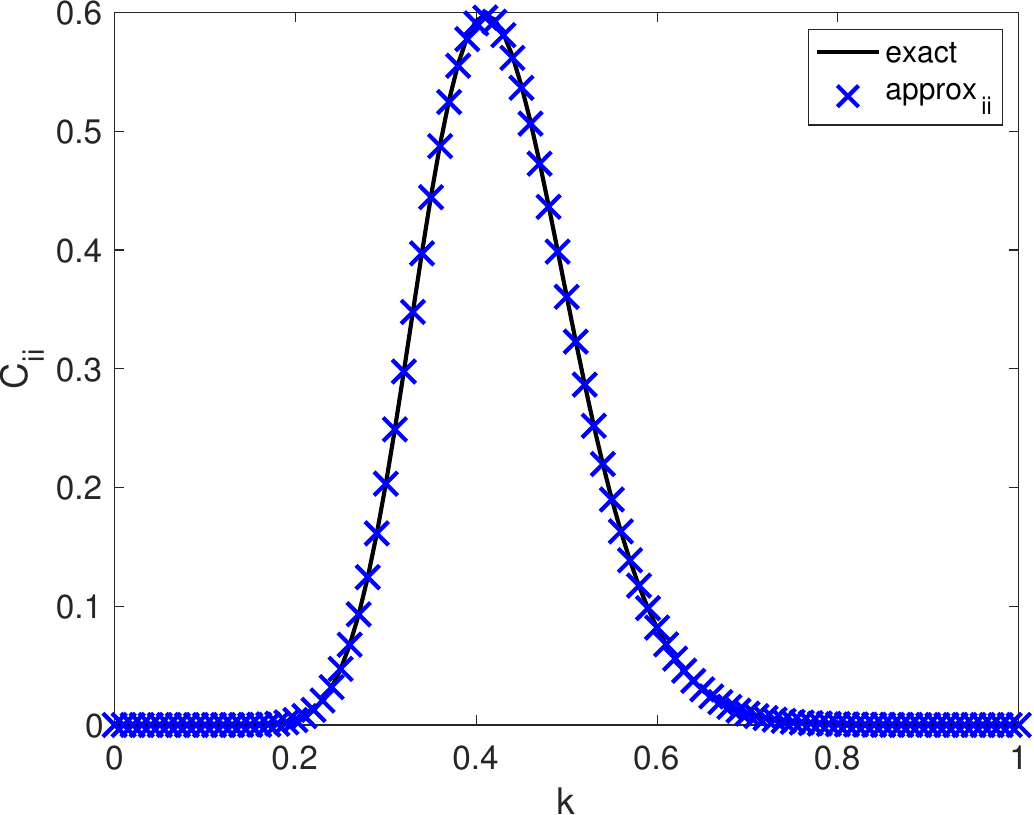}};
        \draw (7cm,0cm) node[anchor=south]{\includegraphics[width=6cm]{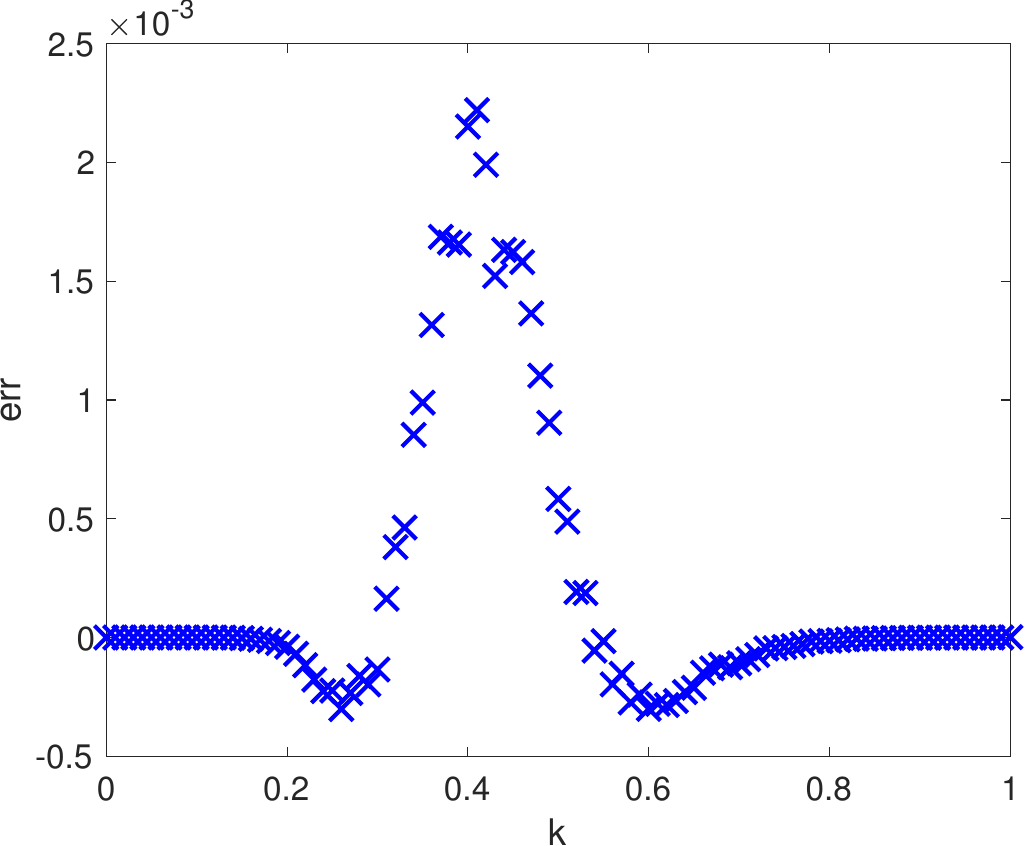}};
      \end{tikzpicture}
  \caption{Plot of exact soln and approximation (left) and the error in the approximation (right) for Test ii. The results are satisfyingly similar, although clearly different, to those in Figure~\ref{fig:T1_soln_meth2}. } 
  \label{fig:T2_soln_meth2}
  \end{center}
  \end{figure}

As before the exact solution is not defined for $k=0$, and so to avoid the singularity we again take $k_a=10^{-16}$, $k_b=1$. Although the exact solution for this case is the same as in Section~\ref{sec:results2DSimple}, the numerical approximation is clearly different.  Figure~\ref{fig:T2_soln_meth2}  shows computed results with $N=100$ discretization points including the exact solution, computed approximation, and errors for the 2nd order method for this test. As mentioned previously, Figure~\ref{fig:1D_conv} also presents results for a convergence study for this case which demonstrates 2nd order accuracy and shows nearly no difference to the results for model i.
  
\subsection{Extension to Two Space Dimensions}
\label{sec:2D}
The fundamental approach is the same in 2D, although there are a number of subtleties requiring discussion and clarification. To facilitate, consider the following continuous model
\begin{align}
  C(\kx,\ky) = \int_{-\infty}^{\infty}\int_{-\infty}^{\infty} f(\xi,\eta;\kx,\ky)\delta(g(\xi,\eta;\kx,\ky))\, d\xi d\eta,
  \label{eq:model2D}
\end{align}
and introduce a uniform grid $\kx_j=\kx_a+j\dkx$ for $j=0,1,\ldots,N_{\kx}$, $\dkx=(\kx_b-\kx_a)/N_{\kx}$, $\ky_l=\ky_a+l\dky$ for $l=0,1,\ldots,N_{\ky}$, $\dky=(\ky_b-\ky_a)/N_{\ky}$. Note that similar to 1D,  the dummy variables of integration have been changed to $(\xi,\eta)$ rather than $\kv_1=(\kx_1,\ky_1)$ in order to simplify the presentation.

The algorithm is again built around the idea of approximating the resonant manifold, i.e. the zero level contour of $g$,  as a piecewise polynomial. However, extension of the prior ideas from 1D leaves certain choices that need to be clarified to fully describe the numerical procedure. In particular, if $g$ is evaluated on a 2D rectangular mesh, there is no unique linear approximation to $g$ over a given cell with 4 corners. The approach taken here is to proceed in a tensor product manner using bilinear interpolation. Specifically, suppose that $\RJ$ and $\RL$ indicate the lower-left corner of a cell for which $g(\xi_{\RJ},\eta_{\RL};\kx_j,\ky_l)$, $g(\xi_{\RJ+1},\eta_{\RL};\kx_j,\ky_l)$, $g(\xi_{\RJ},\eta_{\RL+1};\kx_j,\ky_l)$ and $g(\xi_{\RJ+1},\eta_{\RL+1};\kx_j,\ky_l)$ do not all have the same sign. That is to say that the cell $[\xi_{\RJ},\xi_{\RJ+1}]\times[\xi_{\RL},\xi_{\RL+1}]$ contains part of the resonant manifold. The bi-linear interpolant of $g$ over this cell can be compactly expressed (where we have temporarily suppressed dependence on $\kx_j,\ky_l$ to simplify the notation) as
\begin{align}
  L_g(\xi,\eta) \equiv \frac{1}{\dkx\dky}
    \begin{bmatrix}\xi_{\RJ+1}-\xi\quad \xi-\xi_{\RJ}\end{bmatrix}
    \begin{bmatrix}
      g(\xi_{\RJ},\eta_{\RL}) \quad g(\xi_{\RJ},\eta_{\RL+1})\\
      g(\xi_{\RJ+1},\eta_{\RL}) \quad g(\xi_{\RJ+1},\eta_{\RL+1})
    \end{bmatrix}
    \begin{bmatrix}\eta_{\RL+1}-\eta\\ \eta-\eta_{\RL}\end{bmatrix}.
    \label{eq:gApprox2D}
\end{align}
One can notice that $L_g(\xi,\eta)=0$ is not linear, and that depending on the arrangement of the signs of $g$ on the corners, any given cell can contain 0, 1, or 2 distinct (nonintersecting) zero level curves. For simplicity of presentation the following assumes a single segment, but the extension to 2 segments is straightforward and the implementation accounts for this\footnote{In fact the implementation groups all connected segments of $L_g(\xi,\eta)=0$ into a single piecewise curve, and then loops over those independent continuous pieces.}.  Suppose that the curve defined by $L_g(\xi,\eta)=0$ has endpoints $(\bar{\xi}_a,\bar{\eta}_a)$ and $(\bar{\xi}_b,\bar{\eta}_b)$ on the perimiter of the cell. The approximate resonant manifold is then viewed as a piecewise linear curve connecting those points, e.g. see Figure~\ref{fig:2dContours}. The requisite integration is then approximated over the segment using midpoint quadrature. Now suppose the length, normal, and midpoint to this segment, $\Delta K$, $\nv$, and $(\bar{\xi},\bar{\eta})$ respectively, are defined
\bse
\begin{align}
  \Delta K = \sqrt{(\bar{\xi}_b-\bar{\xi}_a)^2+(\bar{\eta}_b-\bar{\eta}_a)^2},\\
  \nv = \frac{1}{\Delta K}(\bar{\eta}_b-\bar{\eta}_a, \bar{\xi}_b-\bar{\xi}_a),\\
  (\bar{\xi},\bar{\eta}) = \frac{1}{2}(\bar{\xi}_a+\bar{\xi}_b, \bar{\eta}_a+\bar{\eta}_b).
\end{align}
\ese
With this notation in place, the second-order accurate midpoint approximation that defines the scheme is simply (again with dependence on $\kx_j,\ky_l$ suppressed to simplify the notation)
\begin{align}
  \int_{\eta_{\RL}}^{\eta_{\RL+1}}\int_{\xi_{\RJ}}^{\xi_{\RJ+1}} f(\xi,\eta)\delta(g(\xi,\eta))\, d\xi 
    &\approx \frac{L_f(\bar{\xi},\bar{\eta})\Delta K}{|L_{\nabla_{\nv} g}(\bar{\xi},\bar{\eta})|},
    \label{eq:scheme2D_p2}
\end{align}
where $\nabla_{\nv}g=\nv\cdot\nabla g$ is simply the $\nv$ component of the gradient of $g$, and $L_{\nabla_{\nv}g}$ is its bilinear interpolant.

 \subsubsection{Simple 2D Examples}
\begin{figure}[h]
  \begin{center}
    \begin{tikzpicture}[scale=1]
        \draw (0cm,0cm) node[anchor=south]{\includegraphics[width=5cm]{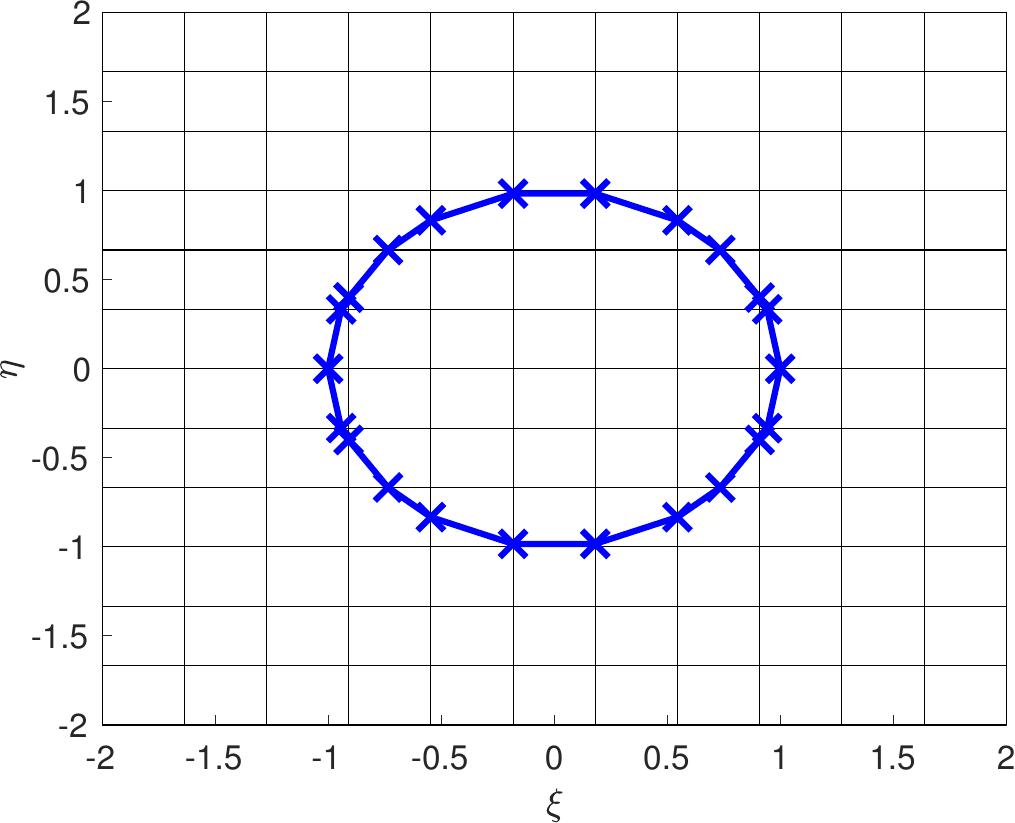}};
        \draw (6cm,0cm) node[anchor=south]{\includegraphics[width=5cm]{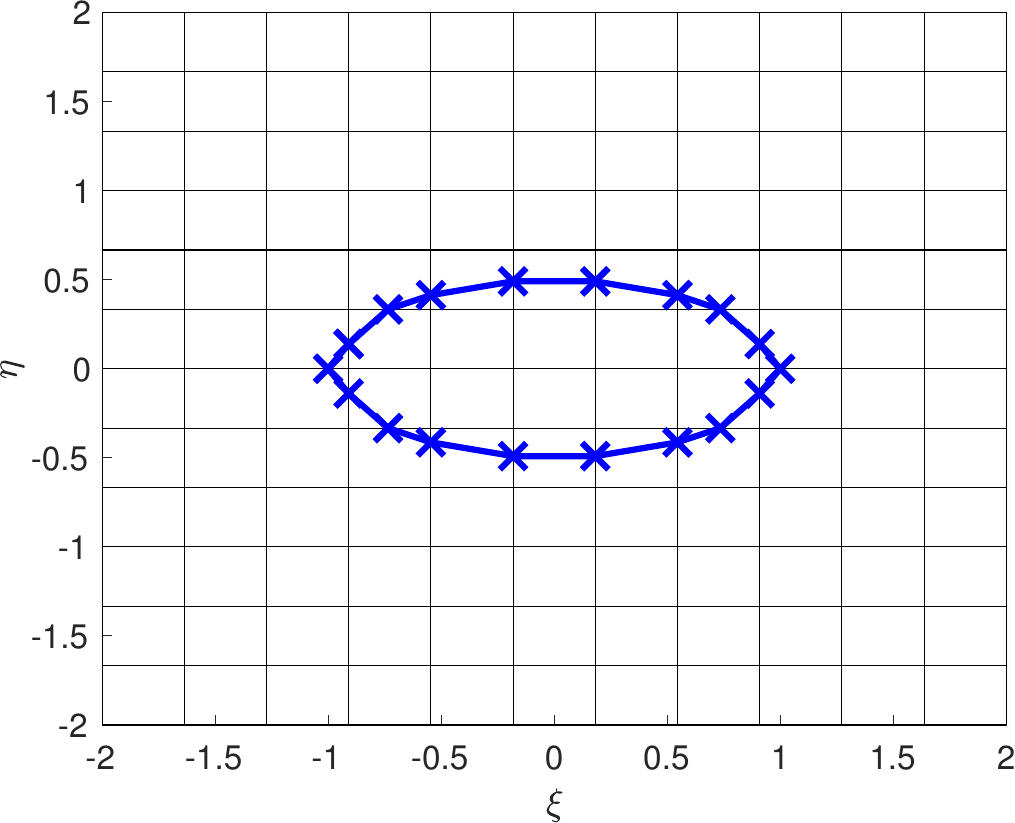}};
        \draw (12cm,0cm) node[anchor=south]{\includegraphics[width=5cm]{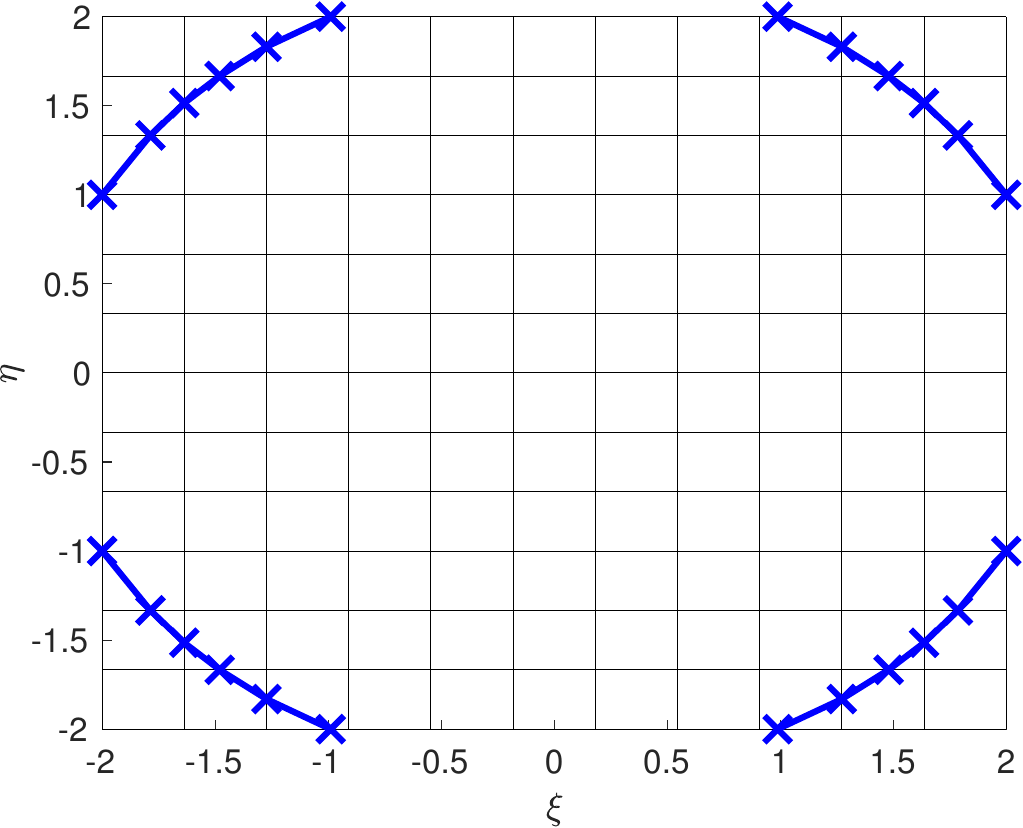}};
      \end{tikzpicture}
  \caption{Plot of approximate resonant manifolds for 2D Test I, V, and VI. Note that for case VI, the resonant manifold extends beyond the computational grid, and so the approximate manifold ``seen'' by the method consists of 4 independent continuous pieces.} 
  \label{fig:2dContours}
  \end{center}
  \end{figure}

\begin{figure}[h]
  \begin{center}
    \begin{tikzpicture}[scale=1]
        \draw (0cm,0cm) node[anchor=south]{\includegraphics[width=8cm]{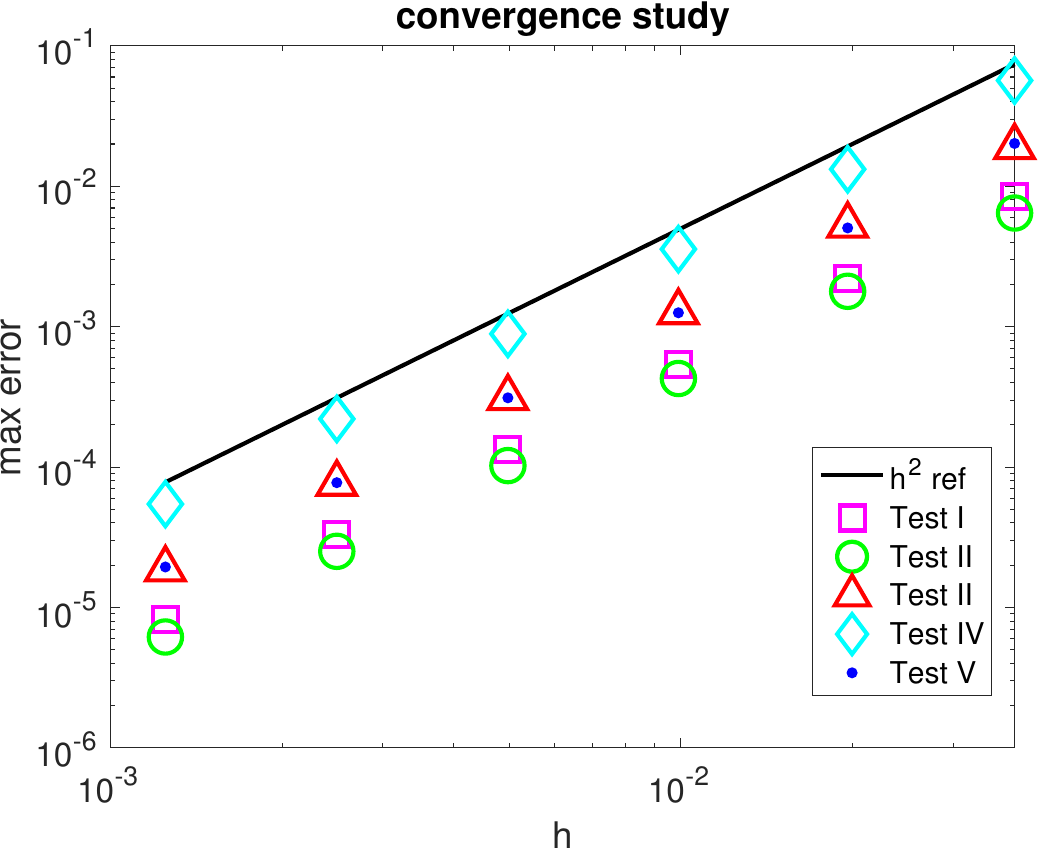}};
      \end{tikzpicture}
  \caption{Convergence vs. grid size for 2D tests independent of $(\kx,\ky)$. Clear second-order convergence is observed in all cases.} 
  \label{fig:2dTestsConv}
  \end{center}  
  \end{figure}

 As a way to help clearly illustrate the mechanics and behavior of the method, let us turn now to a few tests which do not exhibit any dependence on $(\kx,\ky)$. Table~\ref{table:2DExamples} gives definitions for $f$ and $g$ for 6 tests, along with the corresponding exact solutions where relevant. To illustrate the approximate resonant sets, Figure~\ref{fig:2dContours} plots the piecewise linear representation of the resonant manifold for Tests I, V, and IV on a coarse grid. These three examples show tests with a closed circular contour (Test I), a closed elliptical contour (Test V), and a circular contour which extends beyond the assumed grid (Test VI). The last test is included to illustrate the mechanics of the method in the case where the resonant manifold extends beyond the computational domain, and so the notion of an exact solution and convergence studies are somewhat irrelevant. However this case is important in practice since contours formed as a series of disconnected pieces are inevitable in the eventual WKE computations. Figure~\ref{fig:2dTestsConv} shows max-norm convergence studies for Tests I through V, and illustrates the expected 2nd order accuracy. The convergence studies are performed on a domain with $\kx\in[-2,2]$ and $\ky\in[-2,2]$ using grids with $N_{\kx}=25 \cdot 2^m+1$ and $N_{\ky}=26\cdot 2^m+1$ where $m=0,1,\ldots,5$ indicates the refinement of the grid.
\begin{table}
\begin{center}
\begin{tabular}{|c||c|c|c|}
    \hline
    \multicolumn{4}{|c|}{Examples which are independent of $(\kx,\ky)$} \\ \hline
    & $f$ & $g$ & exact result\\ \hline\hline
    Test I & $1$ & $\sqrt{\xi^2+\eta^2}-1$ & $2\pi$ \\ \hline
    Test II & $1$ & $\xi^2+\eta^2-1$ & $\pi$ \\ \hline
    Test III & $\xi^2\eta^4+1$ & $\xi^2+\eta^2-1$ & $\frac{17\pi}{16}$ \\ \hline
    Test IV & $\xi^2\eta^4+1$ & $\xi^2+\eta^2-2$ & $\frac{3\pi}{2}$ \\ \hline
    Test V & $1$ & $\sqrt{\xi^2+4\eta^2}-1$ & $\pi$ \\ \hline
    Test VI & $\xi^2\eta^4+1$ & $\xi^2+\eta^2-5$ & - \\ \hline
  \end{tabular}
\end{center}
\caption{Definition of tests cases for 2D which are independent of $(\kx,\ky)$.}
\label{table:2DExamples}
\end{table}

  \subsubsection{A General 2D RHS Model}
  \label{sec:2DModel}
\begin{figure}[h]
  \begin{center}
    \begin{tikzpicture}[scale=1]
        \draw (0cm,0cm) node[anchor=south]{\includegraphics[width=6cm]{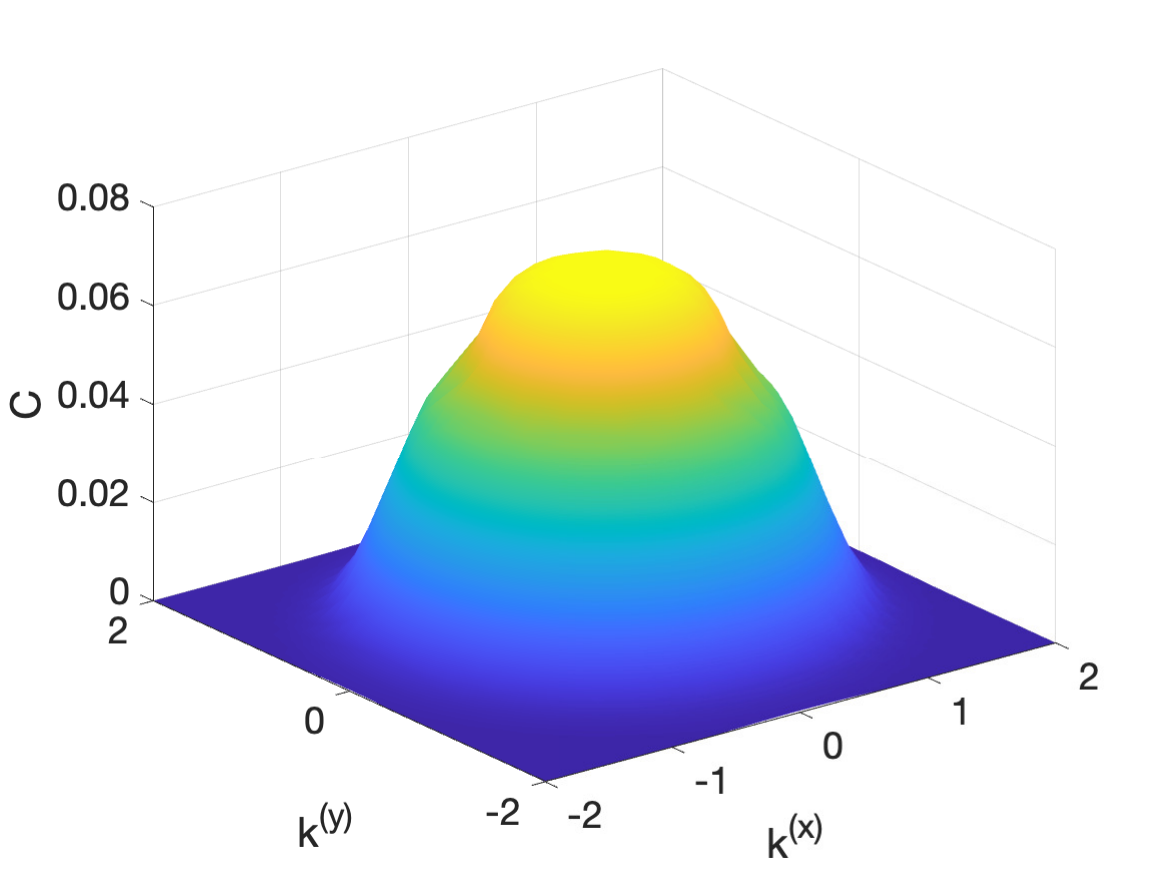}};
        \draw (7cm,0cm) node[anchor=south]{\includegraphics[width=6cm]{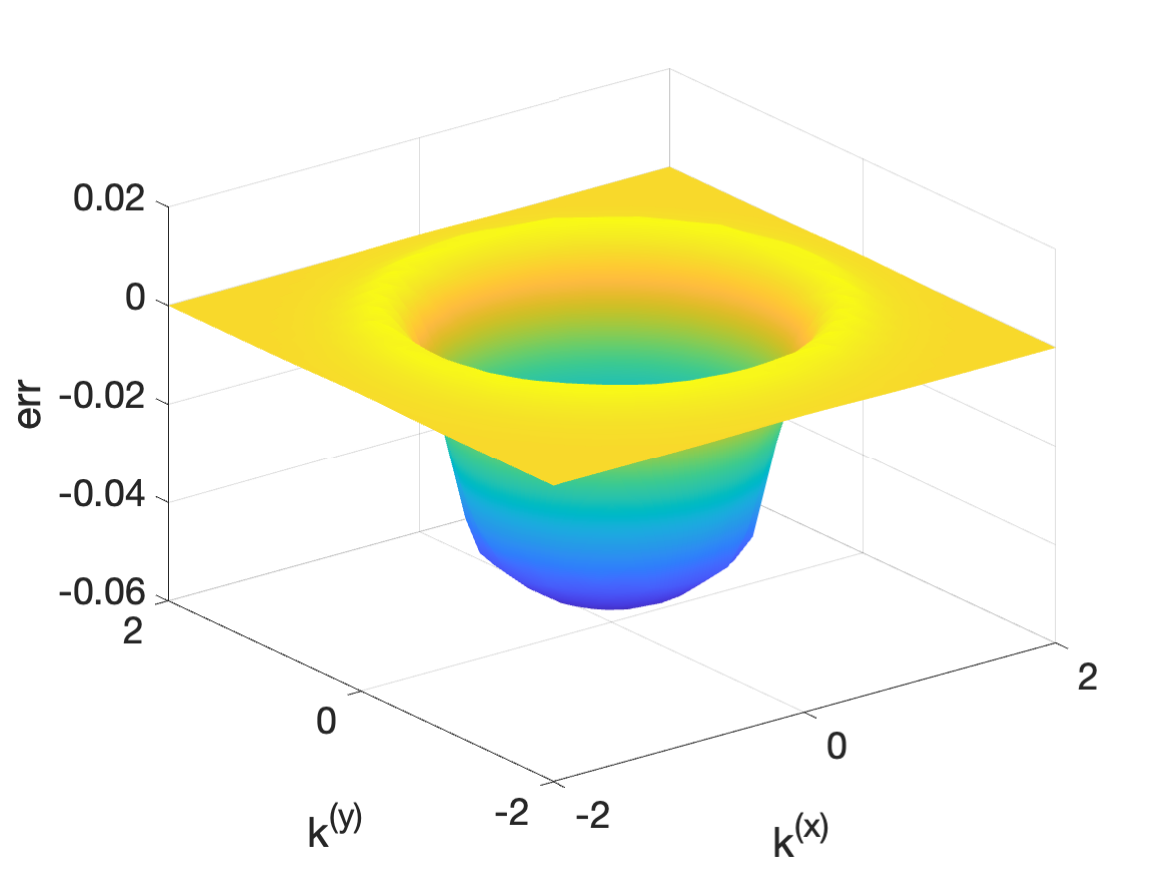}};
      \end{tikzpicture}
  \caption{Plot of computed solution (left), and error (right) for 2D model of Section~\ref{sec:2DModel} with $(\kx,\ky)$-dependence, potentially open integration contours, decaying solutions which are vanishingly small near the artificially truncated boundary, and parts of the kernel that involve a shift, as defined in~\eqref{eq:test2_2DXYExtended} on a grid  with $N_{\kx}=26, N_{\ky}=27$.} 
  \label{fig:test2_2DXYExtended}
  \end{center}
  \end{figure}

\begin{figure}[h]
  \begin{center}
    \begin{tikzpicture}[scale=1]
        \draw (0cm,0cm) node[anchor=south]{\includegraphics[width=8cm]{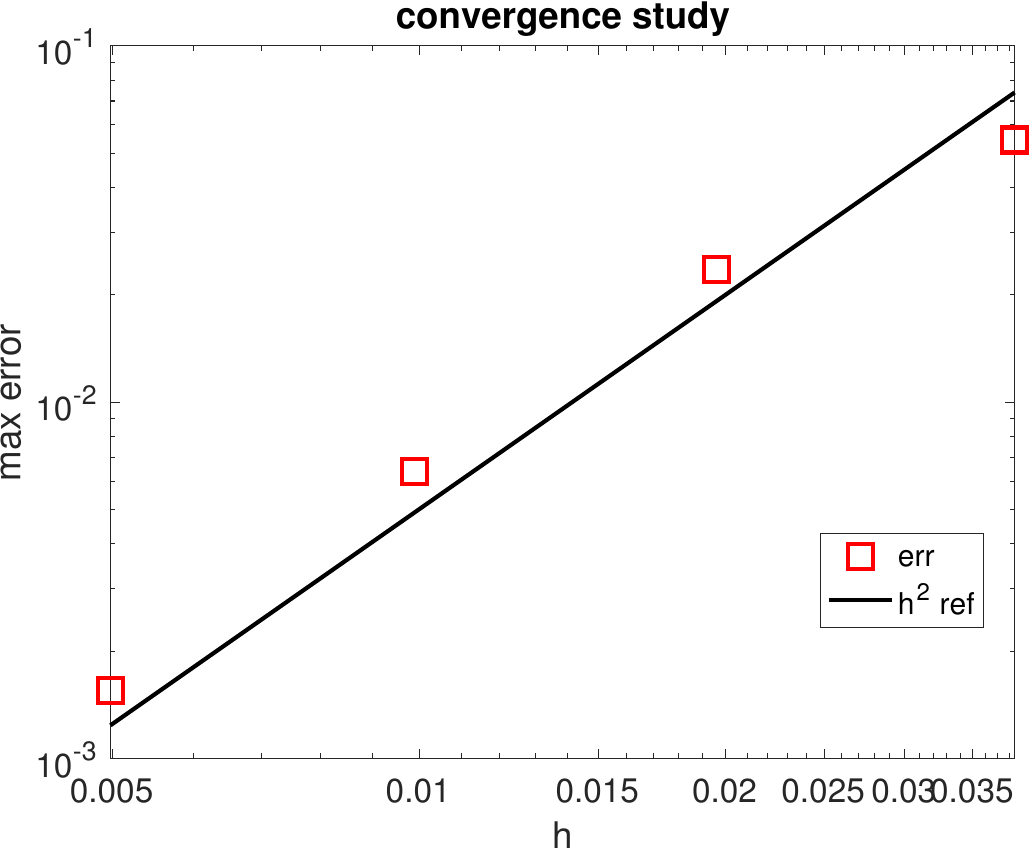}};
      \end{tikzpicture}
  \caption{Convergence vs. grid size for 2D model of Section~\ref{sec:2DModel} with $(\kx,\ky)$-dependence, potentially open integration contours, decaying solutions which are vanishingly small near the artificially truncated boundary, and parts of the kernel that involve a shift, as defined in~\eqref{eq:test2_2DXYExtended}. Here $h$ is a measure of the grid spacing, taken to be $h=4/(N_{\kx}-1)$. Clear second-order convergence is observed.} 
  \label{fig:2dTest2ExtendedConv}
  \end{center}  
  \end{figure}

  Similar to the 1D model described in Section~\ref{sec:extended1D}, here we apply the method to a 2D model with all relevant complications that will be experienced in WKE's. In particular this includes $(\kx,\ky)$-dependence, potentially open integration contours, decaying solutions which are vanishingly small near the artificially truncated boundary, and parts of the kernel that involve a shift (potentially outside the domain). To do so, consider
  \begin{align}
    f(\xi;\kx,\ky) = \fk(\kx,\ky)\fkb(\xi,\eta)\fka(\xi-\hat{\xi},\eta-\hat{\eta}),
    \label{eq:model2DExtended}
  \end{align}
where
\bse
\label{eq:test2_2DXYExtended}
\begin{align}
  \fk(\kx,\ky) & = 1\\
  \fkb(\xi,\eta) & = e^{-50(\left((\xi-\xi_0)^2+(\eta-\eta_0)^2\right)}\\
  \fka(\tilde{\xi},\tilde{\eta}) & = e^{-50(\left(\tilde{\xi}^2+\tilde{\eta}^2\right)}\\
  \tilde{\xi} & = \xi-\hat{\xi}\\
  \tilde{\eta} & = \eta-\hat{\eta}.
\end{align}
Taking $\hat{\xi}=\xi_0$ and $\hat{\eta}=\eta_0$ and assuming an infinite integration domain gives the exact solution
\begin{align}
  C(\kx,\ky) & = \frac{\sqrt{2\pi}}{20}e^{-50\alpha\left((\kx-\xi_0)^2+(\ky-\eta_0)^2\right)}.
\end{align}
\ese
Figure~\ref{fig:test2_2DXYExtended} shows the computed solution and error on a truncated domain with $\kx\in[-2,2]$ and $\ky\in[-2,2]$ using a grid with $N_{\kx}=26$ and $N_{\ky}=27$ points respectively. Here the method is executed at each point in the $(\kx,\ky)$ space, including the evaluation of integration contours for the resonant manifolds and approximate quadrature. Note that the error shows no spurious behavior near the artificially truncated boundaries, indicating that the present truncation method is sufficient in these cases. Figure~\ref{fig:2dTest2ExtendedConv} shows the related convergence study which clearly illustrates 2nd order convergence. The study uses $\kx\in[-2,2]$ and $\ky\in[-2,2]$ with $N_{\kx}=25 \cdot 2^m+1$ and $N_{\ky}=26\cdot 2^m+1$ where $m=0,1,\ldots,3$.

\section{Time Dependent Solutions to the WKE}
\label{sec:verification}
Having developed the core of a numerical approach, and demonstrated its behavior on idealized collision integral models in 1D and 2D, we proceed to investigate the behavior of the method for the WKE equation~\eqref{eq:NLSWKE}. This discussion begins in Section~\ref{sec:WKEManipulations} with a brief description of the manipulations taken so that the method described in Section~\ref{sec:2D} is applicable. The approach uses a method-of-lines type of discrete formulation, and the time-stepping method that is used is also described in Section~\ref{sec:WKEManipulations}. A number of test cases are defined, and sample calculations are shown in Section~\ref{sec:ICs}. Section~\ref{sec:selfConvergence} then presents results of various self-convergence studies for isotropic and anisotropic solutions. Finally, Section~\ref{sec:NLSComparison} relates the WKE solution to ensemble averages of the NLS equations of motion~\eqref{eq:NLS2D} for the isotropic and anisotropic cases.

\subsection{Discretization of the WKE}
\label{sec:WKEManipulations}
The general approach to discretizing the governing WKE~\eqref{eq:NLSWKE} will rely on a method-of-lines  formulation and explicit time-stepping. However, the techniques described in Section~\ref{sec:2D} are not immediately applicable to the right-hand-side of~\eqref{eq:NLSWKE}. The approach taken here is to eliminate the linear resonant condition, and thereby one of the integrations, by solving for $\kv_1$ as a function of $\kv$ and $\kv_2$, e.g. $\kv_1+\kv_2-\kv=0$ implies $\kv_1=\kv-\kv_2$. As a result,~\eqref{eq:NLSWKE} becomes
\ba
  \p_{\tau}\nk = &\frac{1}{\pi}\int
    \left(n_{\kv-\kv_2}n_{\kv_2}-\nk n_{\kv-\kv_2}-\nk n_{\kv_2}\right)
      \delta(\omega_{\kv-\kv_2}+\omega_{\kv_2}-\omega_{\kv})\nonumber\\
    & \qquad +\left(2n_{\kv+\kv_2}n_{\kv_2}-2\nk n_{\kv_2}+2\nk n_{\kv+\kv_2}\right)
      \delta(\omega_{\kv+\kv_2}-\omega_{\kv_2}-\omega_{\kv})\,d\kv_2.
  \label{eq:NLSWKE_1Resonance}
\ea
The right-hand side of~\eqref{eq:NLSWKE_1Resonance} consists of a summation of terms which are each of the form considered in Section~\ref{sec:2D}, and in particular in Equation~\eqref{eq:test2_2DXYExtended}. As a result, the methods described in Section~\ref{sec:2D} are directly applicable. To be specific, again introduce the uniform grid $\kx_j=\kx_a+j\dkx$ for $j=0,1,\ldots,N_{\kx}$, $\dkx=(\kx_b-\kx_a)/N_{\kx}$, $\ky_l=\ky_a+l\dky$ for $l=0,1,\ldots,N_{\ky}$, $\dky=(\ky_b-\ky_a)/N_{\ky}$. Further, let $\jv=(j,l)$ be a multi index so that $\kv_{\jv}=(\kx_j,\ky_l)$ and $\omega_{\jv}=\omega_{\kv_{\jv}}$, and the numerical approximation of the WKE solution at $\kv_{\jv}$ is then indicated as $\nkj\approx n_{\kv_{\jv}}$. 

As an illustrative example, consider the first term on the right-hand-side of~\eqref{eq:NLSWKE_1Resonance}. After truncating to the finite computational domain, the integral is broken into a sum of integrals over cells as 
\ba
    \int n_{\kv-\kv_2}n_{\kv_2}&\delta(\omega_{\kv-\kv_2}+\omega_{\kv_2}-\omega_{\kv})\,d\kv_2
     \approx
    \sum_{\jv_2}\int_{\ky_{l_2-1}}^{\ky_{l_2}}\int_{\kx_{j_2-1}}^{\kx_{j_2}}n_{\kv-\kv_2}n_{\kv_2}\delta(\omega_{\kv-\kv_2}+\omega_{\kv_2}-\omega_{\kv})\,d\kx_2 d\ky_2
\ea
where $\sum_{\jv_2}=\sum_{j_2=1}^{N_{\kx}}\sum_{l_2=1}^{N_{\ky}}$ is used as notational shorthand. Identifying $f_1=n_{\kv-\kv_2}n_{\kv_2}$ and $g_1=\omega_{\kv-\kv_2}+\omega_{\kv_2}-\omega_{\kv}$, the numerical approximation in~\eqref{eq:scheme2D_p2} is then applied to give
\ba
    \int n_{\kv-\kv_2}n_{\kv_2}&\delta(\omega_{\kv-\kv_2}+\omega_{\kv_2}-\omega_{\kv})\,d\kv_2 \approx
    \sum_{\jv_2}\frac{L_{f_1}(\overline{\kv}_{\jv_2})\Delta K_{\jv_2}}{|L_{\nabla_{\nv} {g_1}}(\overline{\kv}_{\jv_2})|}.
\ea
Similar manipulations are performed for the remaining 5 terms in~\eqref{eq:NLSWKE_1Resonance}. The spatially discrete, $\tau$-continuous set of ODEs for the WKE, i.e. the method-of-lines formulation, is then given as
\bse
\label{eq:MOL}
\ba
  \p_{\tau}\nkj &=
    \sum_{\jv_2}\frac{L_{f_a}(\overline{\kv}_{\jv_2})\Delta K_{\jv_2}}{|L_{\nabla_{\nv} {g_a}}(\overline{\kv}_{\jv_2})|}+
    \sum_{\jv_2}\frac{L_{f_b}(\overline{\kv}_{\jv_2})\Delta K_{\jv_2}}{|L_{\nabla_{\nv} {g_b}}(\overline{\kv}_{\jv_2})|},\\
    f_a &= n_{\jv-\jv_2}n_{\jv_2}-\nkj n_{\jv-\jv_2}-\nkj n_{\jv_2}, \quad\qquad   g_a=\omega_{\jv-\jv_2}+\omega_{\jv_2}-\omega_{\jv},\\
  f_b &= 2n_{\jv+\jv_2}n_{\jv_2}-2\nkj n_{\jv_2}+2\nkj n_{\jv+\jv_2}, \qquad   g_b=\omega_{\jv+\jv_2}-\omega_{\jv_2}-\omega_{\jv}.
\ea
\ese

For time advancement of the equations~\eqref{eq:MOL}, classical 4th-order accurate explicit Runge-Kutta is then applied~\cite{ascher98}. The only remaining item to complete the description of the discretization approach is the selection of the time step. Note that WKE~\eqref{eq:NLSWKE} does not contain any $\kv$ derivatives, and so formally from a stability perspective the time step is independent of $\dkx$ and $\dky$. Nevertheless, in order that the overall scheme converges as the $\kv$-space mesh is refined, the time step is selected according to
\bse
\ba
  \Delta \tau = \Lambda_{S} \frac{h}{\lambda_h},
\ea
where $\Lambda_{S}$ is a ``safety'' factor (here set as $\Lambda_{S}=2.25$), $h=\min(\dkx,\dky)$, and $\lambda_h$ is chosen as an approximation of the linear time stepping eigenvalue. Within reasonable limits, the precise choice of $\lambda_h$ has little effect on the overall performance, and the present work uses
\ba
  \lambda_h & = \frac{1}{\pi}\max_{\jv}\left|
    \sum_{\jv_2} \frac{L_{f_{\lambda_1}}(\overline{\kv}_{\jv_2})\Delta K_{\jv_2}}{|L_{\nabla_{\nv} {g_{a}}}(\overline{\kv}_{\jv_2})|}+
    \sum_{\jv_2} \frac{L_{f_{\lambda_2}}(\overline{\kv}_{\jv_2})\Delta K_{\jv_2}}{|L_{\nabla_{\nv} {g_{b}}}(\overline{\kv}_{\jv_2})|}
  \right|,\\
  f_{\lambda_1} & = -n_{\kv-\kv_2}-n_{\kv_2},\\
  f_{\lambda_2} & = -2 n_{\kv_2}+2n_{\kv+\kv_2}.
\ea
\ese
This choice simply uses an approximation $\lambda_h\approx\lambda$ where 
\[
  \lambda = \frac{1}{\pi}\max_{\kv}\left|\left[ \int
    \left(-n_{\kv-\kv_2}-n_{\kv_2}\right)
      \delta(\omega_{\kv-\kv_2}+\omega_{\kv_2}-\omega_{\kv})
      +\left(-2 n_{\kv_2}+2n_{\kv+\kv_2}\right)
      \delta(\omega_{\kv+\kv_2}-\omega_{\kv_2}-\omega_{\kv})\,d\kv_2\right]\right|,
\]
and is computed along side the discretization itself.

\subsection{Test Cases and Approximate Solutions}
\label{sec:ICs}
\begin{figure}[h]
  \begin{center}
    \begin{tikzpicture}[scale=1]
        \draw (0cm,0cm) node[anchor=south]{\includegraphics[trim=2.4cm 7.6cm 2cm 6.75cm,clip,width=7cm]{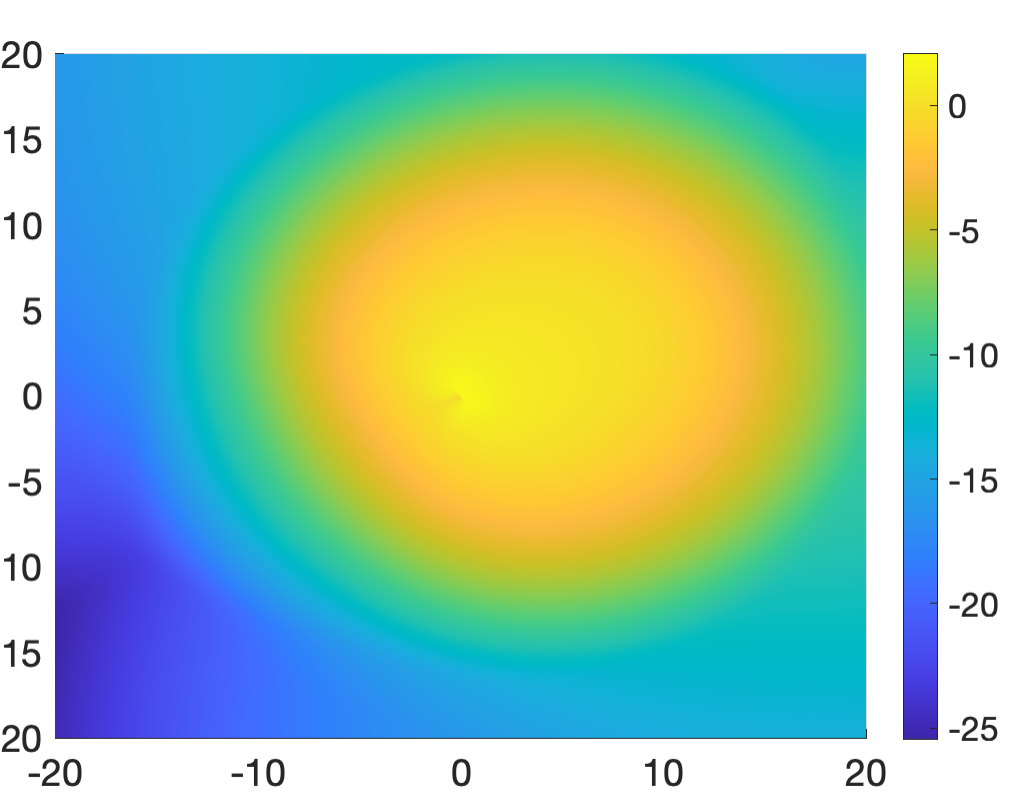}};
        \draw (9cm,0cm) node[anchor=south]{\includegraphics[trim=2.4cm 7.6cm 2cm 6.75cm,clip,width=7cm]{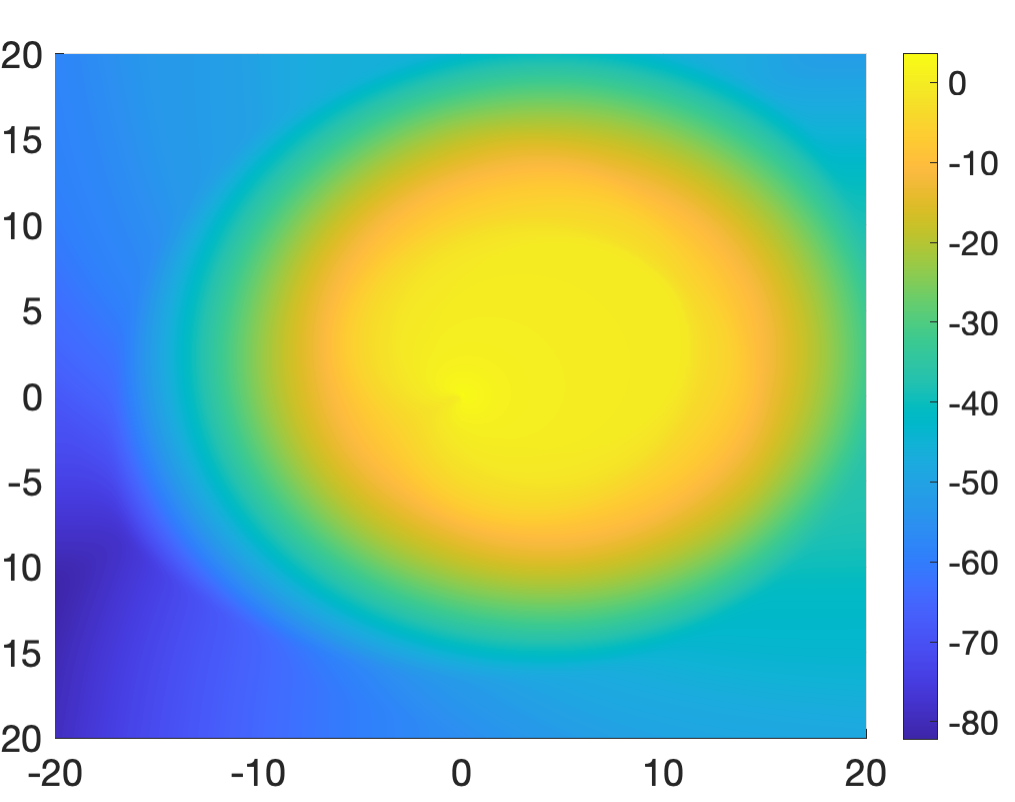}};
        \draw (-.3,0.1)node[anchor=north,fill=white]{\tiny $k_x$};
        \draw (-4,3.)node[anchor=north,fill=white,rotate=90]{\tiny $k_y$};
        \draw (8.7,0.1)node[anchor=north,fill=white]{\tiny $k_x$};
        \draw (5,3.)node[anchor=north,fill=white,rotate=90]{\tiny $k_y$};
        \draw (-.2,5.6)node[anchor=north,draw,fill=white,rounded corners]{\tiny Offset Gaussian ($\tau=1$)};
        \draw (8.8,5.6)node[anchor=north,draw,fill=white,rounded corners]{\tiny Offset top-hat ($\tau=1$)};
      \end{tikzpicture}
  \caption{Plot of $\wk=\log_{10}(\nk)$ at $\tau=1$ for initial Gaussian initial condition with $A=5$ (left), and top-hat initial condition with $r_0=2\pi$ (right). Both cases are anisotropic and use $\kx_0=5$ and $\ky_0=3$. } 
  \label{fig:timeEvolve}
  \end{center}
  \end{figure}
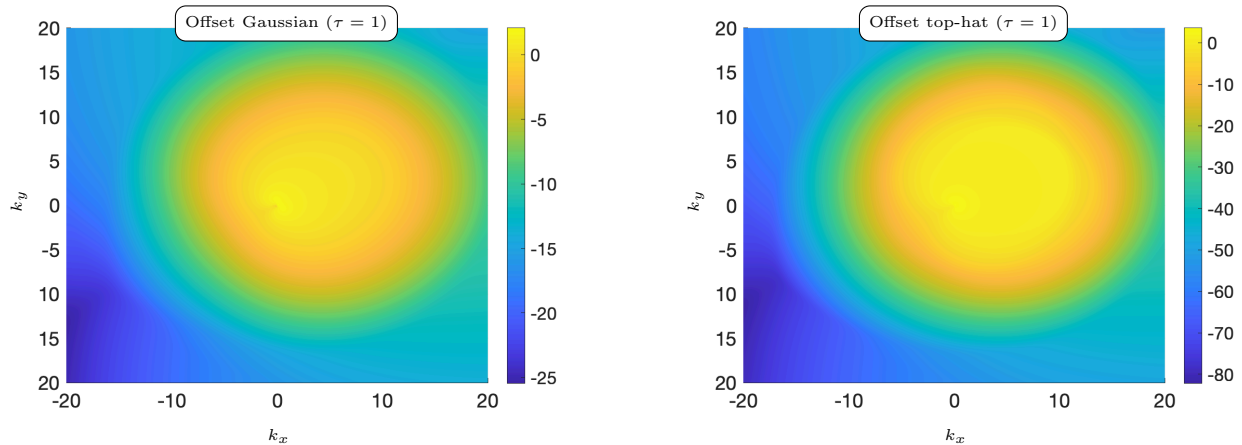

The full discretization of the WKE~\eqref{eq:NLSWKE} will be applied to a variety of case, both isotropic and anisotropic. The initial conditions for these runs will be defined by either $\nk(0)=\mathcal{N}_{G}$ or $\nk(0)=\mathcal{N}_{TH}$ where
\bse
\label{eq:ICs}
\ba
  \mathcal{N}_{G}(\kv) & = Ae^{-\frac{1}{10}\|\kv-\kv_0\|^2}, \\
  \mathcal{N}_{TH}(\kv) & = 
  \begin{cases}
    1 & \qquad \hbox{if $\|\kv-\kv_0\|\le r_0$}\\
    0 & \qquad \hbox{else}.
  \end{cases}
\ea
\ese
Here $\mathcal{N}_{G}$ is a Gaussian with amplitude $A$, while $\mathcal{N}_{TH}$ is a {\em top-hat} function of radius $r_0$. In both cases  and $\kv_0=(\kx_0,\ky_0)$ represents a centering so that $\kv_0=0$ will correspond to an isotropic case, while $\kv_0\ne0$ will give an anisotropic case. These two functions are chosen as prototypical smooth and non-smooth profiles respectively, with the top-hat case used as a more strenuous numerical test. Here, and in the remainder of the manuscript, $\log_{10}(\nk)$ is considered rather than $\nk$ for two reasons. First, the log is more often used in practice. Second, the log more clearly illustrates various features present in the solutions. Also related to the second point, the log-space solution has magnitude which is more evenly distributed so that error norms tend to more accurately reflect the general trend across the domain rather than the behavior at a single point (this is a similar effect to considering point-wise relative norms which is not done here). To make this more clear in the presentation, define
\ba
  \wk\equiv\log_{10}(\nk).
\ea
To give a sense of the time-evolved solution for these two initial conditions, Figure~\ref{fig:timeEvolve} shows $\wk(\tau=1)$ computed with $(\kx,\ky)\in[-20,20]^2$ on a grid with $N_{\kx}=200, N_{\ky}=208$ for two cases. Note that the grid size in $\kx$ is chosen different than $\ky$ for essentially no other reason than to more fully exercise the code\footnote{ Note that for this WKE (and the discretization of it), the solution at $\kv=0$ is constant in time which leads to a singular solution, as a result the grids are deliberately selected so that there is no grid point at $\kv=0$.}. Figure~\ref{fig:timeEvolve} shows only the anisotropic cases with $\kx_0=5$ and $\ky_0=3$ since they are more interesting than the related isotropic solutions. The Gaussian initial condition uses $A=5$, while the top-hat uses $r_0=2\pi$. These parameters are chosen without any real intent, other than to show nontrivial dynamics at $\tau=1$.

\subsection{Self-Convergence Studies}
\label{sec:selfConvergence}
Given initial conditions of the form~\eqref{eq:ICs}, there is no known closed-form exact solution for the WKE. Therefore to quantify the convergence behavior of the discretization, Richardson extrapolation self-convergence studies are performed using the procedure described in~\cite{pog2008a,banks09a,banks13b_RELDD}. Recall that this procedure assumes $\wk(\tau)=w_{E}(\kv,\tau)+c(\kv,\tau)h^p$ where $w_E$ is the exact solution, $c$ is a constant function of $\kv$ and $\tau$ but independent of grid spacing, and $h$ is a measure of the grid spacing. For a given discrete norm, computations at 3 levels of refinement (here $50\times52, 100\times104$, and $200\times208$) are performed, from whence estimates for $w_E$, $c$, and $p$ are derived. Note that the discrete max norm tends to be the most revealing, but rates for weaker norms (e.g. $l_1$-norm or $l_2$-norm) are also computed for completeness. 
 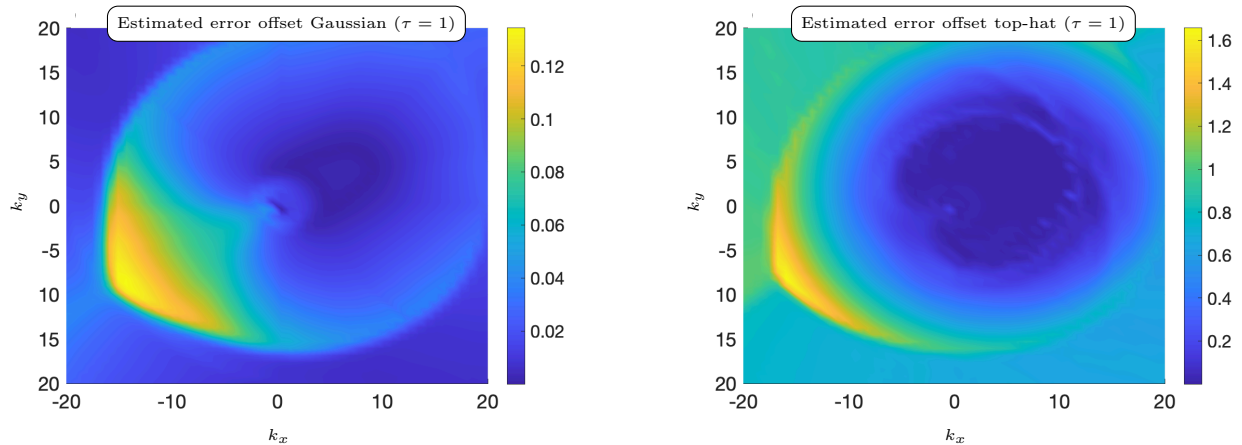
\begin{figure}[h]
  \begin{center}
    \begin{tikzpicture}[scale=1]
        \draw (0cm,0cm) node[anchor=south]{\includegraphics[trim=2.4cm 7.6cm 2cm 6.75cm,clip,width=7cm]{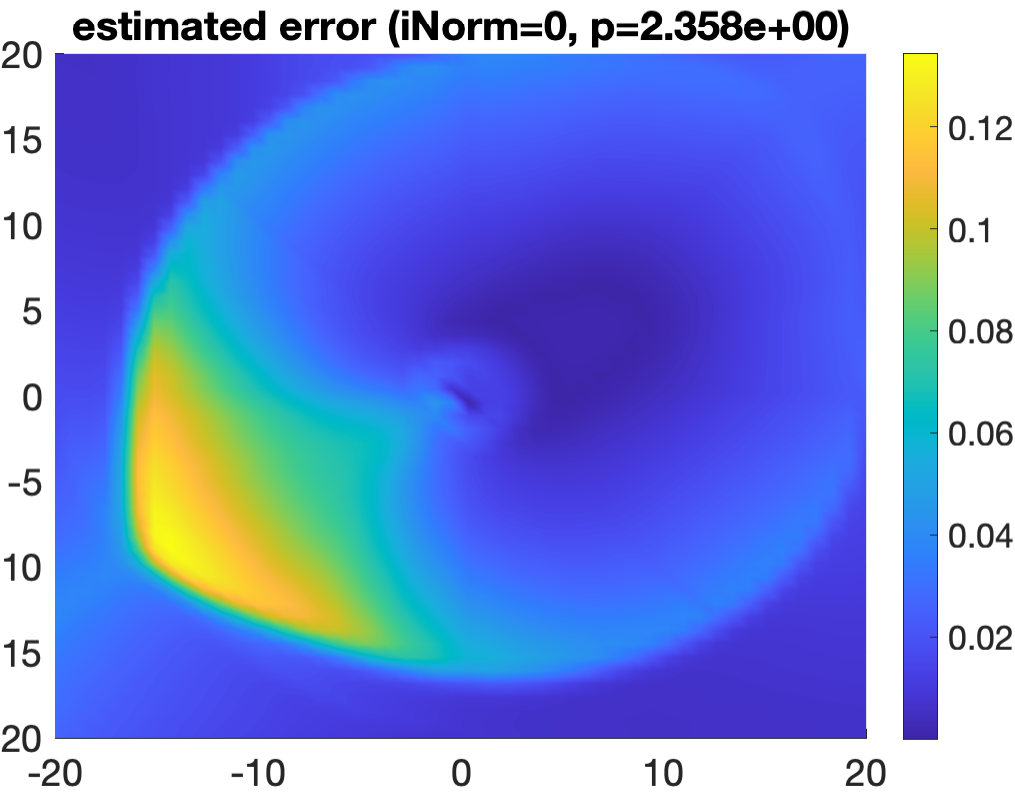}};
        \draw (9cm,0cm) node[anchor=south]{\includegraphics[trim=2.4cm 7.6cm 2cm 6.75cm,clip,width=7cm]{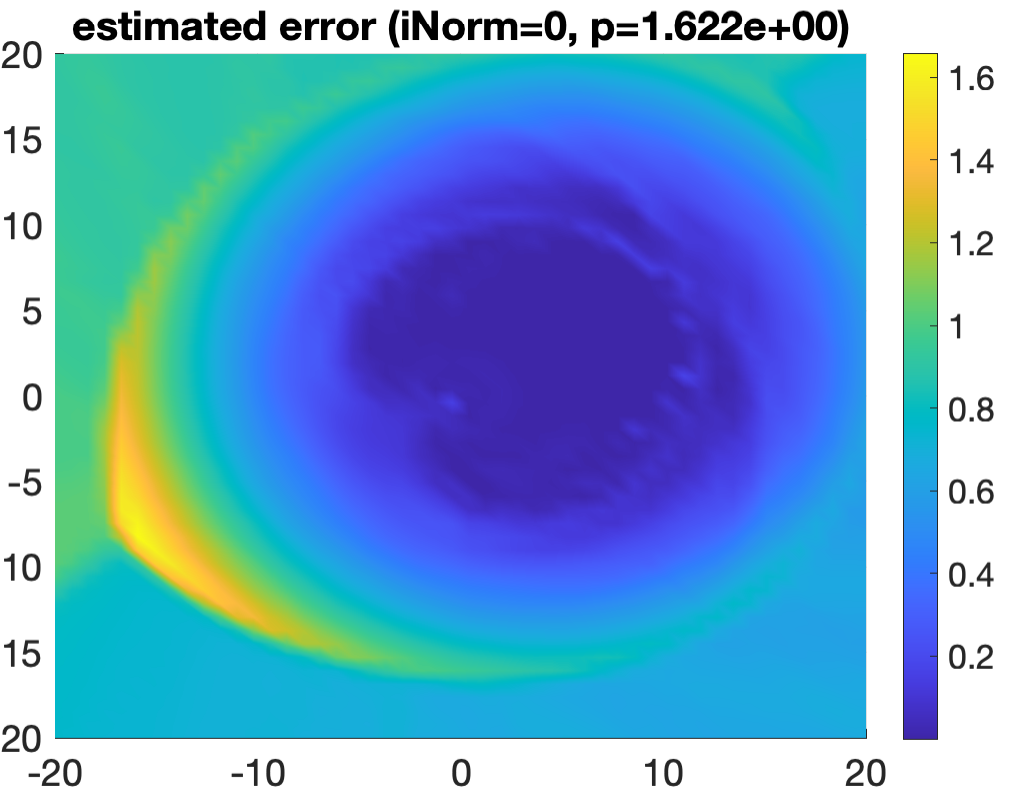}};
        \draw (-.3,0.1)node[anchor=north,fill=white]{\tiny $k_x$};
        \draw (-4,3.)node[anchor=north,fill=white,rotate=90]{\tiny $k_y$};
        \draw (8.7,0.1)node[anchor=north,fill=white]{\tiny $k_x$};
        \draw (5,3.)node[anchor=north,fill=white,rotate=90]{\tiny $k_y$};
        \draw [fill=white, draw=none](-3,5.32) rectangle (2.5,5.6);
        \draw [fill=white, draw=none](6,5.32) rectangle (11.5,5.6);
        \draw (-.2,5.6)node[anchor=north,draw,fill=white,rounded corners]{\tiny Estimated error offset Gaussian ($\tau=1$)};
        \draw (8.8,5.6)node[anchor=north,draw,fill=white,rounded corners]{\tiny Estimated error offset top-hat ($\tau=1$)};
      \end{tikzpicture}
  \caption{Plot of estimated max-norm error in $\wk$ at $\tau=1$ for the anisotropic Gaussian (left), and anisotropic top-hat (right). In both cases, the maximum error occurs in the domain interior, and the error appears to be well-behaved throughout the domain. For the top-hat, one can observe effects from the discontinuous initial conditions, which speaks to the difficulty of this as a quantitative test.} 
  \label{fig:RE}
  \end{center}
  \end{figure}
Figure~\ref{fig:RE} shows the estimated max-norm errors in $\wk$ on the finest grid. For both cases, the error appears to be reasonably well behaved, which indicates that the numerical approximation is of good quality. However, the top-hat case does retain some features associated with the use of a discontinuous initial condition. Table~\ref{table:REEstimates} presents estimated convergence rates for all four cases (both isotropic and anisotropic, using both Gaussian and top-hat initial conditions) in the discrete $l_1$-norm, discrete $l_2$-norm, and max-norm. In all cases the rates are observed to be reasonably close to the expected 2nd order, with the smooth cases somewhat better than the non-smooth cases, as is expected. 
\begin{table}
\begin{center}
\begin{tabular}{|c||c|c|c|c|}
    \hline
    				& \multicolumn{2}{|c|}{Gaussian} & \multicolumn{2}{|c|}{top-hat} \\ 
     \strut 			& isotropic & anisotropic &  isotropic & anisotropic\\ \hline\hline
    $l_1$-norm 		& $2.47$ 	& $2.36$	& $1.69$	& $1.62$\\ \hline
    $l_2$-norm 		& $2.34$ 	& $2.45$	& $1.69$	& $1.63$\\ \hline
    $l_{\infty}$-norm 	& $2.35$ 	& $2.42$	& $1.69$	& $1.62$\\ \hline
  \end{tabular}\\\bigskip
\end{center}
\caption{Estimated convergence rates for $w=\log_{10}(\nk)$ using Richardson extrapolation as described in~\cite{pog2008a,banks09a,banks13b_RELDD}.}
\label{table:REEstimates}
\end{table}

\subsection{Comparison to Ensemble Averages of Dynamical Equations}
\label{sec:NLSComparison}
As a final study of the newly described discretization technique for the WKE~\eqref{eq:NLSWKE}, direct comparison to ensemble averaged solutions of the original nonlinear Schr{\"o}dinger equation~\eqref{eq:NLS2D} are now performed using the top-hat initial condition. This study is similar to the one presented in~\cite{banks22_WKE}, which probed the validity of the wave kinetic theory for the quintic NLS in 1D. In the present work, the intent of the comparison is primarily benchmarking of the new discretization technique, although the study does serve the dual role of probing the relationship of wave kinetic theory to direct ensemble averages. To that end, the derivation of the WKE~\eqref{eq:NLSWKE} from the governing NLS system~\eqref{eq:NLS2D}, which is presented in~\ref{sec:WKEDerive}, provides the necessary details.

The WKE is derived in the limit of increasingly large spatially periodic systems. However, for numerical purposes, simulations on finite sized boxes must be used, and so, as in~\cite{banks22_WKE}, the procedure will be to consider simulations of the NLS~\eqref{eq:NLS2D} on increasingly large boxes. To ensure validity of the wave kinetic theory, the large box limit must be taken in conjunction with the real time satisfying $t\to\infty$ as in~\eqref{eq:tToInfinity}, the magnitude of the nonlinearity tending to zero with $\epsilon\to 0$ (where here $\epsilon$ is a pointwise measure of the magnitude of the solution) see~\eqref{eq:akBasic}, and for order 1 kinetic times satisfying $\tau=t\epsilon^2=O(1)$ as in~\eqref{eq:tauToZero}. To obtain this distinguished limit, and following the description in~\cite{banks22_WKE}, the setup for a given ensemble solution of~\eqref{eq:NLS2D} on a periodic box with $\xv\in[-L/2,L/2]^2$ proceeds as follows. Let a parameter $P\in[0,1]$ be given, and set 
\bse
\ba
  \epsilon & =\frac{1}{L^P},\\
  t_F & = \frac{1}{\epsilon^2} = L^{2P},
\ea
\ese
where $t_F$ is the real final time corresponding to the kinetic time $\tau=1$. For all cases in the present manuscript, $P=\frac{3}{4}$ is chosen. The initial condition for a single member of the ensemble is then taken as 
\bse
\ba
  u(\xv,0) & = \frac{\epsilon}{L}\sum_{\jv} a_{\kv_{\jv}}(0) e^{i(\kv_{\jv}\cdot\xv+\phi_{\kv})},\\
  a_{\kv_{\jv}}(0) & = \mathcal{N}_{TH}\left(\kv_{\jv}\right),
\ea
\ese
where $j\in \mathbb Z^2$, $\kv_{\jv}=\frac{2\pi\jv}{L}$, and $\phi_{\kv}$ is a random phase taken from a uniform distribution in $[0,2\pi]$. To make the comparison with the WKE, the expectation (denoted using angle brackets $<.>$) is taken over an ensemble with $N_e$ members, and so 
\ba
  \nk(\tau) = <\ak(\tau),\akb(\tau)>.
\ea

Although it is not the focus of the present manuscript, it is essential to at least briefly describe the numerical setup and techniques used to discretize the NLS~\eqref{eq:NLS2D}. For the domain $\xv\in[-L/2,L/2]^2$, a uniform grid with $M_x \times M_y$ points is used, and periodic boundary conditions applied. The spatial grid uses $8th$ order accurate central finite differences~\cite{banks25FD,fornberg96,fornberg98} with no artificial dissipation. For time-stepping, traditional explicit $4th$-order accurate Runge-Kutta is used~\cite{ascher98}. This choice results in a maximum stable time step $\Delta t\sim\dx^2$ where $\dx$ is the spatial grid size. As a result, the order-of-accuracy for the method is $O(\Delta x^8,\Delta t^4)=O(\Delta x^8)$. To keep the cost manageable, relatively small ensembles are run with only $180$ valid members. As in~\cite{banks22_WKE}, for certain instances of initialization, particularly when $L$ is small, it is possible for a given simulation to fail to run to completion. The typical case for failure results when, in a particular realization of random data, an accidental alignment of phases leads to a locally large solution. In such a case, either the time step decreases to an unmanageably small number, or numerical overflow occurs due to instability from insufficiently small time steps. In either case the run is considered to have failed and is removed from the ensemble. Also note that as in~\cite{banks22_WKE}, as $L$ increases this situation becomes less common and for sufficiently large $L$ all instances run to completion. 

Before proceeding to comparisons between ensemble averages of the NLS~\eqref{eq:NLS2D}, and approximate solutions of~\eqref{eq:NLSWKE}, let us briefly discuss some of the practical challenges associated with the ensemble averages. The existence of various limit processes, both analytical and numerical, leads to rather significant numerical expense in order to sufficiently probe the limit. First, the analytical limit of $L\to\infty$ obviously requires sufficiently large physical domains. Then for any discretization on a finite domain, one would in principle require $M_x\to \infty$ and $M_y\to\infty$. Put together, one would therefore expect that as $L$ increases, both $M_x$ and $M_y$ must increase faster. Adding to the challenge is the fact that in frequency space, the solutions appear to possess quickly decaying tails, and due to finite precision effects one might expect that any range of wave amplitudes approaching machine precision, i.e. 15 orders of magnitude, should invariably be viewed with some skepticism (recall that in real-space simulations of the NLS, small amplitude waves will inevitably be added to large ones). Furthermore, wrapped around this is the requirement of an infinite, or at least very large, ensemble. These observations provide a kind of baseline that should guide the subsequent comparisons, and we will attempt to probe some of these issues.

\subsubsection{Isotropic Case}
\begin{figure}
  \begin{center}
    \begin{tikzpicture}[scale=1]
      \useasboundingbox (-2.75,-.3) rectangle (13,4.5);
        \draw (0cm,0cm) node[anchor=south]{\includegraphics[trim=2.5cm 7.75cm 4cm 6.75cm,clip,width=4.5cm]{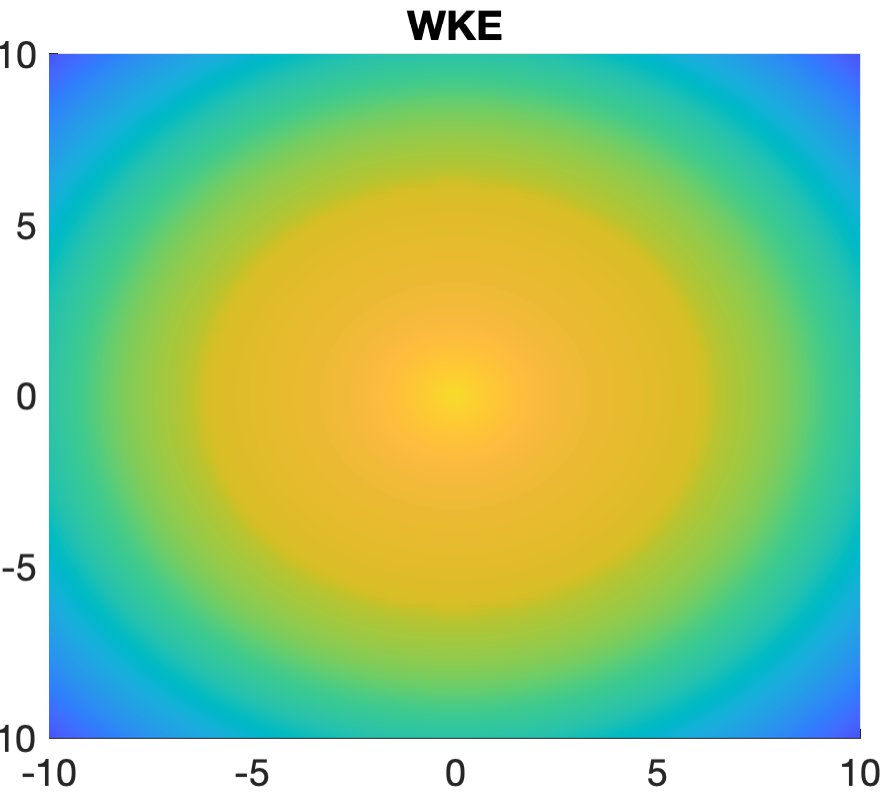}};
        \draw (5cm,0cm) node[anchor=south]{\includegraphics[trim=2.5cm 7.75cm 4cm 6.75cm,clip,width=4.5cm]{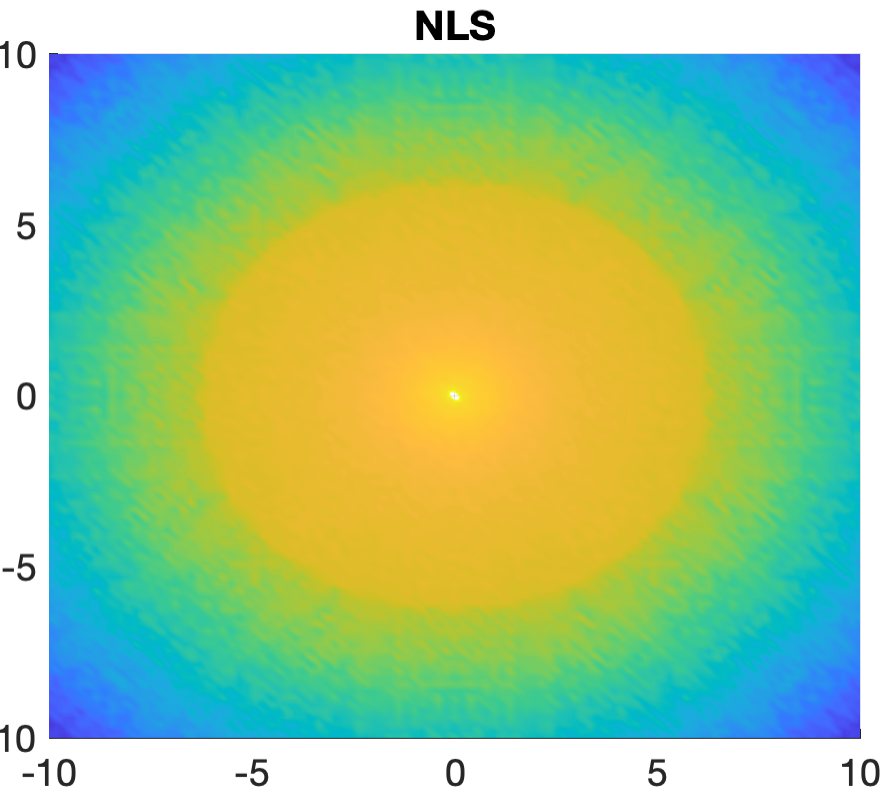}};
        \draw (10cm,0cm) node[anchor=south]{\includegraphics[trim=2.5cm 7.75cm 4cm 6.75cm,clip,width=4.5cm]{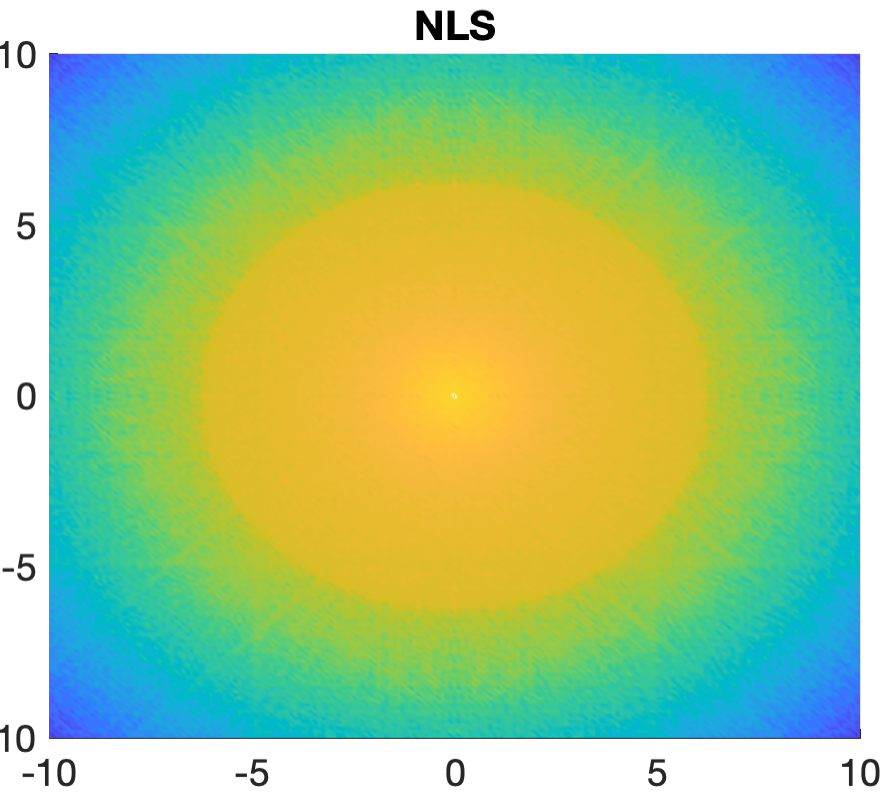}};
        \draw (12.5cm,.25cm) node[anchor=south]{\includegraphics[trim=17.5cm 7.75cm 2.5cm 6.75cm,clip,height=3.5cm]{images/Avg_L64_N1000}};
        \draw (0.05,4.25)node[anchor=north,draw,fill=white,rounded corners]{\tiny WKE ($\tau=1$)};
        \draw (5.05,4.25)node[anchor=north,draw,fill=white,rounded corners]{\tiny NLS ($L=32$, $M_x=M_y=301$)};
        \draw (10.05,4.25)node[anchor=north,draw,fill=white,rounded corners]{\tiny NLS ($L=64$, $M_x=M_y=1001$)};
        \draw (0.05,0.1)node[anchor=north,fill=white]{\tiny $k_x$};
        \draw (5.05,0.1)node[anchor=north,fill=white]{\tiny $k_x$};
        \draw (10.05,0.1)node[anchor=north,fill=white]{\tiny $k_x$};
        \draw (-2.75,2.25)node[anchor=north,fill=white,rotate=90]{\tiny $k_y$};
        %
      \end{tikzpicture}
  \caption{At left is the solution of the WKE at $\tau=1$. At center is a 180 member ensemble for the NLS with $L=32$ and $M_x=M_y=301$. At right is a 180 member ensemble for the NLS with $L=64$ and $M_x=M_y=1001$. All results are plotted using $\log_{10}$, which is also indicated by the color table. The similarity of the solutions is clear, as is the roughness inherent to the relatively small ensembles.} 
  \label{fig:Comp2D}
  \end{center}
  \end{figure}
Consider the isotropic top-hat case given by initial condition with $\nk(0)=\mathcal{N}_{TH}(\kv)$ where $r_0=2\pi$ and $\kv_0=0$. Figure~\ref{fig:Comp2D} shows computed results for the WKE~\eqref{eq:NLSWKE} with $\kx\in[-20,20]$, $\ky\in[-20,20]$, $N_{\kx}=200$, $N_{\ky}=202$, and computed to a final time $\tau=1$. For plotting the domain has been restricted to $[-10,10]^2$, where the majority of the dynamics are occurring. Figure~\ref{fig:Comp2D} also shows the result of ensemble averages with 180 members for the NLS~\eqref{eq:NLS2D} for 2 cases, the first with $L=32$ and $M_x=M_y=301$, and the second with $L=64$ and $M_x=M_y=1001$.

\begin{figure}
  \begin{center}
    \begin{tikzpicture}[scale=.7]
      \useasboundingbox (-4.5,2.25) rectangle (18,7);
        \draw (0cm,0cm) node[anchor=south]{\includegraphics[width=5cm]{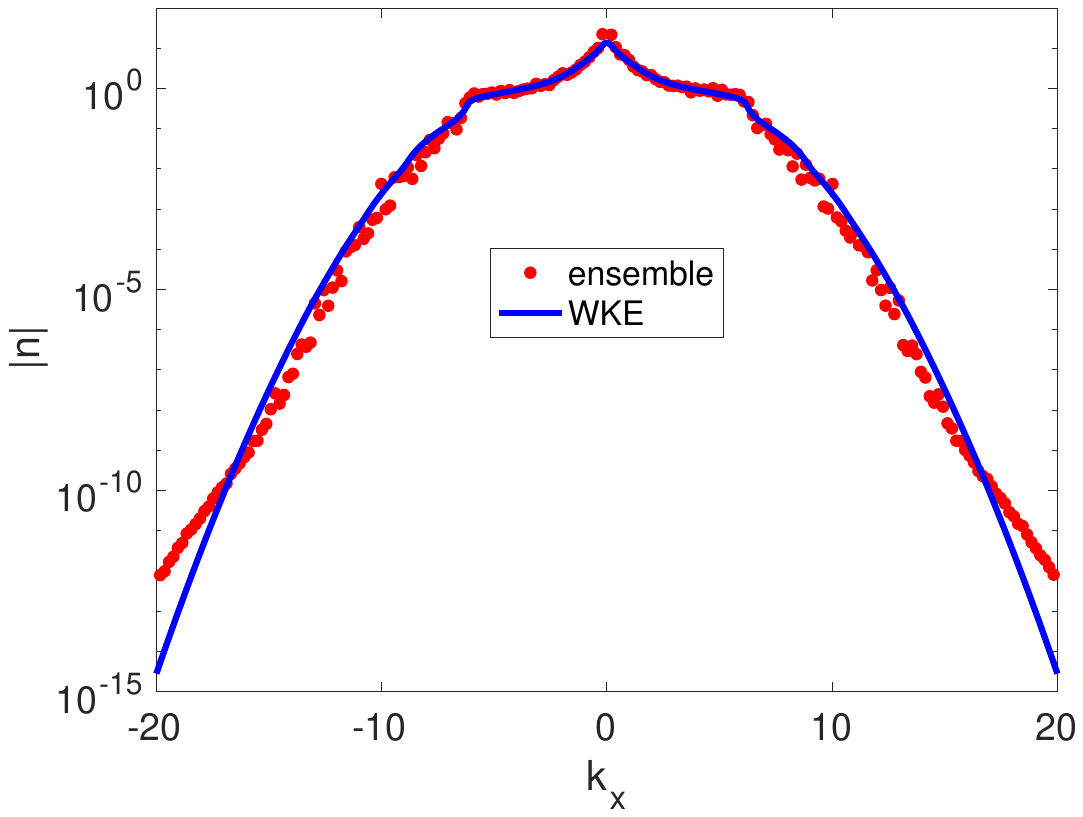}};
        \draw (7cm,0cm) node[anchor=south]{\includegraphics[width=5cm]{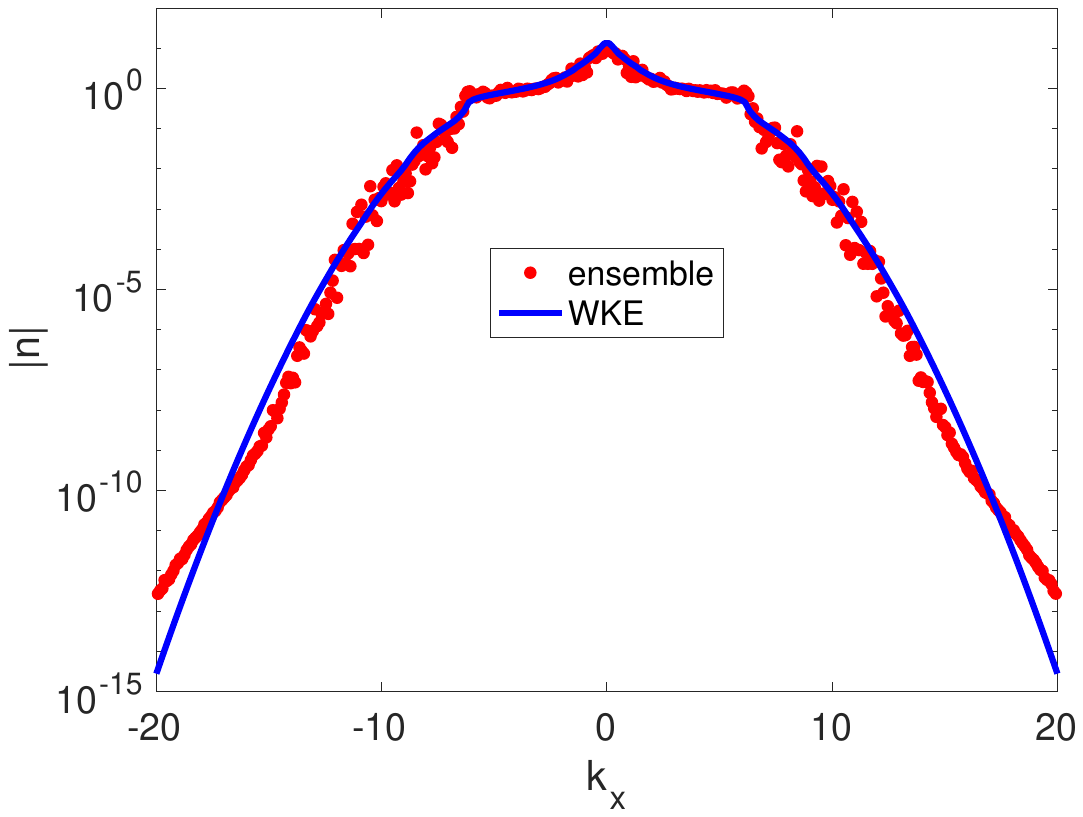}};
        \draw (14cm,0cm) node[anchor=south]{\includegraphics[width=5cm]{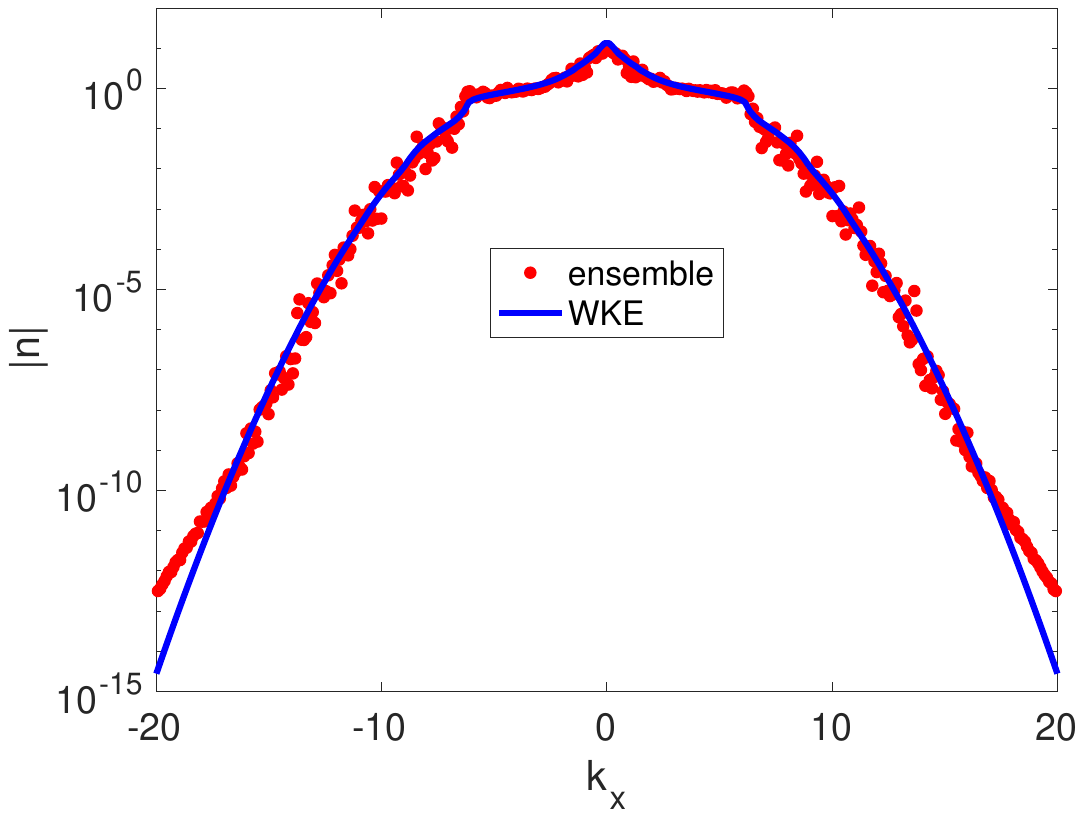}};
        \draw (0.1,7.55)node[anchor=north,draw,fill=white,rounded corners]{\tiny $L=32$, $M_x=M_y=301$};
        \draw (7.1,7.55)node[anchor=north,draw,fill=white,rounded corners]{\tiny $L=64$, $M_x=M_y=601$};
        \draw (14.1,7.55)node[anchor=north,draw,fill=white,rounded corners]{\tiny $L=64$, $M_x=M_y=1001$};
      \end{tikzpicture}
  \caption{Line slices along $\ky=0$ comparing the WKE solution to NLS ensemble averages. At left the NLS runs use $L=32$ and $M_x=M_y=301$. At center the NLS runs use $L=64$ and $M_x=M_y=601$. At right the NLS runs use  $L=64$ and $M_x=M_y=1001$.} 
  \label{fig:Comp1D_y0}
  \end{center}
  \end{figure}

\begin{figure}
  \begin{center}
    \begin{tikzpicture}[scale=.7]
      \useasboundingbox (-4.5,2.25) rectangle (18,7);
        \draw (0cm,0cm) node[anchor=south]{\includegraphics[width=5cm]{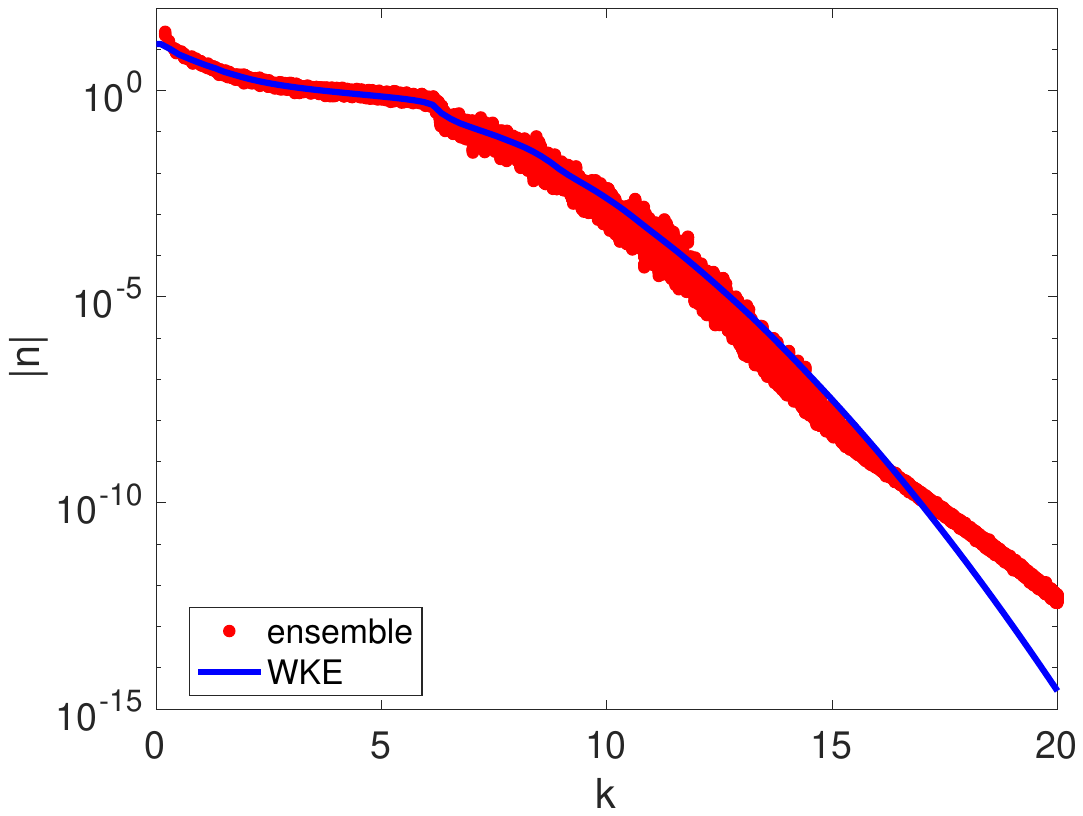}};
        \draw (7cm,0cm) node[anchor=south]{\includegraphics[width=5cm]{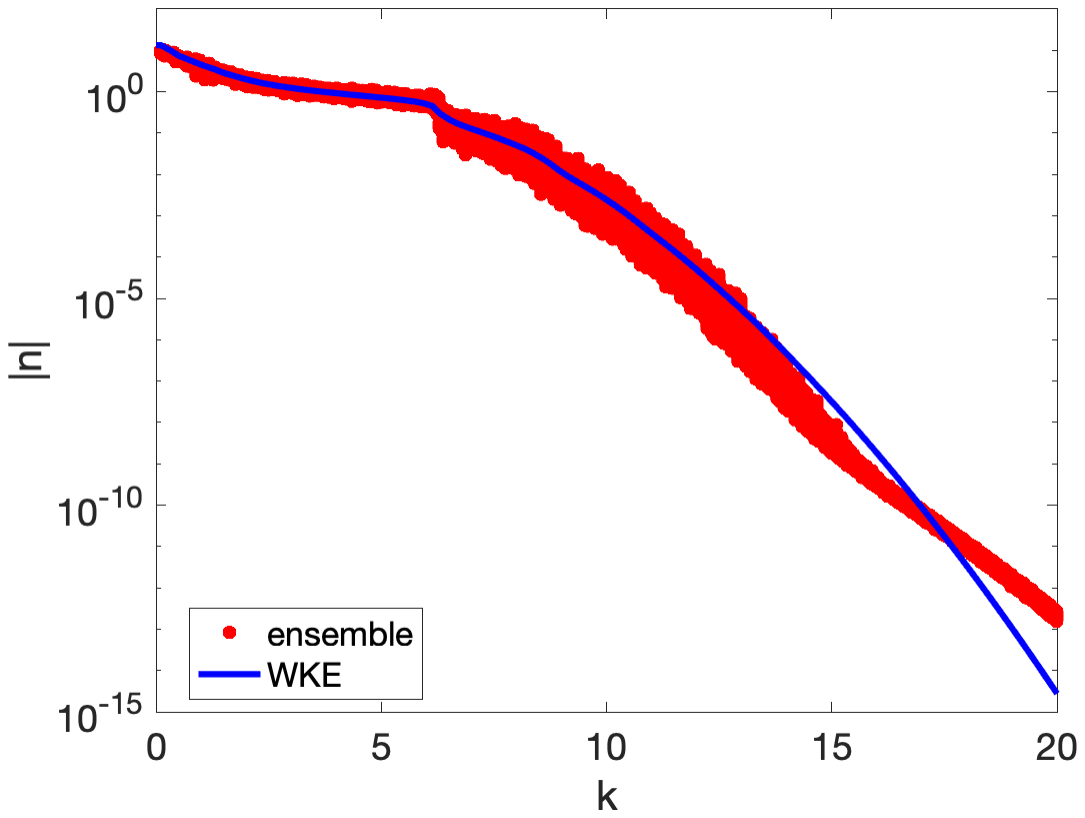}};
        \draw (14cm,0cm) node[anchor=south]{\includegraphics[width=5cm]{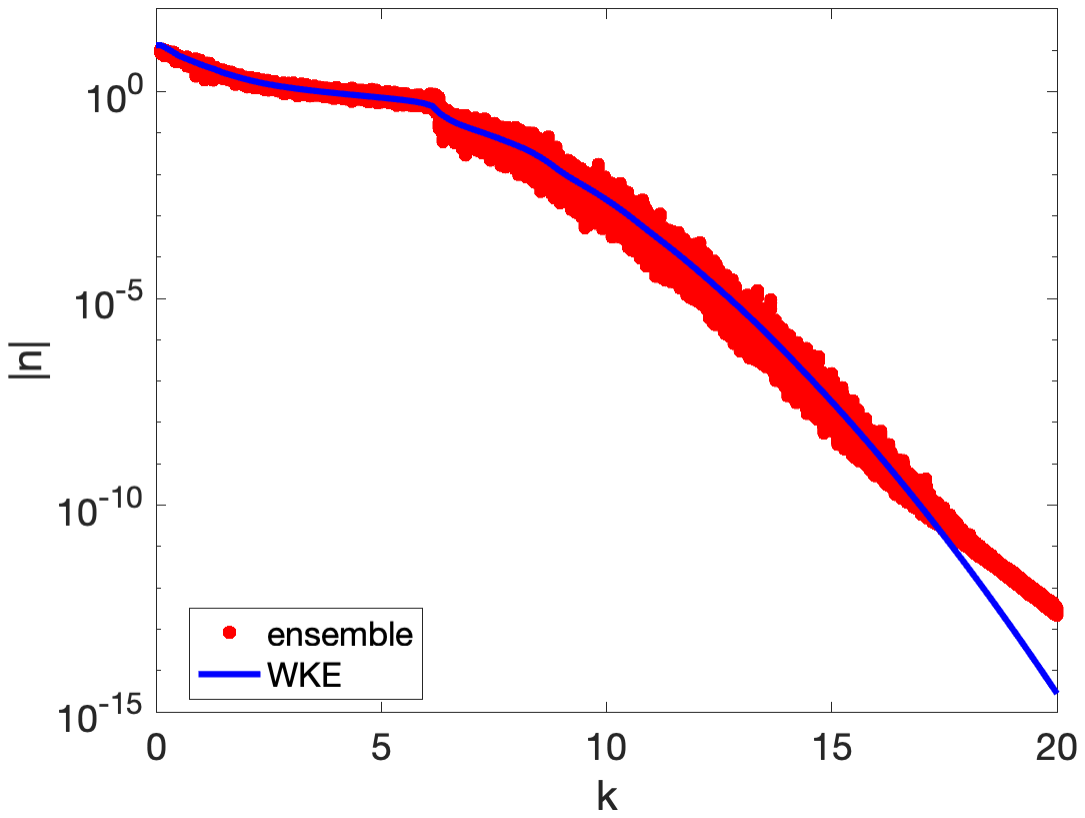}};
        \draw (0.1,7.55)node[anchor=north,draw,fill=white,rounded corners]{\tiny $L=32$, $M_x=M_y=301$};
        \draw (7.1,7.55)node[anchor=north,draw,fill=white,rounded corners]{\tiny $L=64$, $M_x=M_y=601$};
        \draw (14.1,7.55)node[anchor=north,draw,fill=white,rounded corners]{\tiny $L=64$, $M_x=M_y=1001$};
      \end{tikzpicture}
  \caption{Comparison of the WKE solution for $\ky=0$ to the NLS solutions plotted vs. radius $k$ for every k-space point. At left the NLS runs use $L=32$ and $M_x=M_y=301$. At center the NLS runs use $L=64$ and $M_x=M_y=601$. At right the NLS runs use  $L=64$ and $M_x=M_y=1001$.} 
  \label{fig:scatComp1D_y2}
  \end{center}
  \end{figure}

To probe the results discussed in Figure~\ref{fig:Comp2D}, Figure~\ref{fig:Comp1D_y0} presents slices of the various solutions along the line $\ky=0$ for $\kx\in[-20,20]$. Recall that the WKE is computed for $\kx\in[-20,20]$, but the $\kx$ bounds for the ensemble NLS solutions may be much larger and is determined by $L$ and $M_x$. The figure presents results for 3 cases; the first with $L=32$ and $M_x=M_y=301$, the second with $L=64$ and $M_x=M_y=601$, and the third with $L=64$ and $M_x=M_y=1001$. For the first to the second, the length of the box is increased while keeping the $\xv$-grid spacing the same. Moving to the third case then the larger box is retained, and the spatial resolution is increase. Increasing the box size has the effect of increasing the k-space resolution for small $k$ whereas increasing the spatial resolution has the effect to move the high-k boundary further out. Indeed the plots make clear that both larger boxes and finer $\xv$-grids improve the results, as expected. However it is perhaps surprising that $L=64$ and $M_x=M_y=1001$, which corresponds to roughly $\kx\in[-50,50]$, is required for reasonable agreement of the ensemble and WKE solutions, even for $\kx\in[-7,7]$. Also note that the NLS solutions begin to deviate significantly from the WKE solution once the variation of wave amplitudes becomes too large, here roughly 12 orders of magnitude. This is a consistent phenomenon across all our simulations. 
As a final comparison for the isotropic case, isotropy is leveraged and the ensemble NLS solutions are plotted vs. radius $k$ for every k-space point. Comparison is made to the WKE solution for $\ky=0$ (which is essentially equivalent to any other radial line). Figure~\ref{fig:scatComp1D_y2} presents these results for the same cases as presented in Figure~\ref{fig:Comp1D_y0}, and the prior conclusions regarding $L$ and $\xv$-grid resolution are even more convincingly demonstrated.

\subsubsection{Anisotropic Case}
\begin{figure}[h]
  \begin{center}
    \begin{tikzpicture}[scale=1]
        \draw (0cm,0cm) node[anchor=south]{\includegraphics[trim=2.5cm 7.5cm 4cm 6.75cm,clip,width=4.5cm]{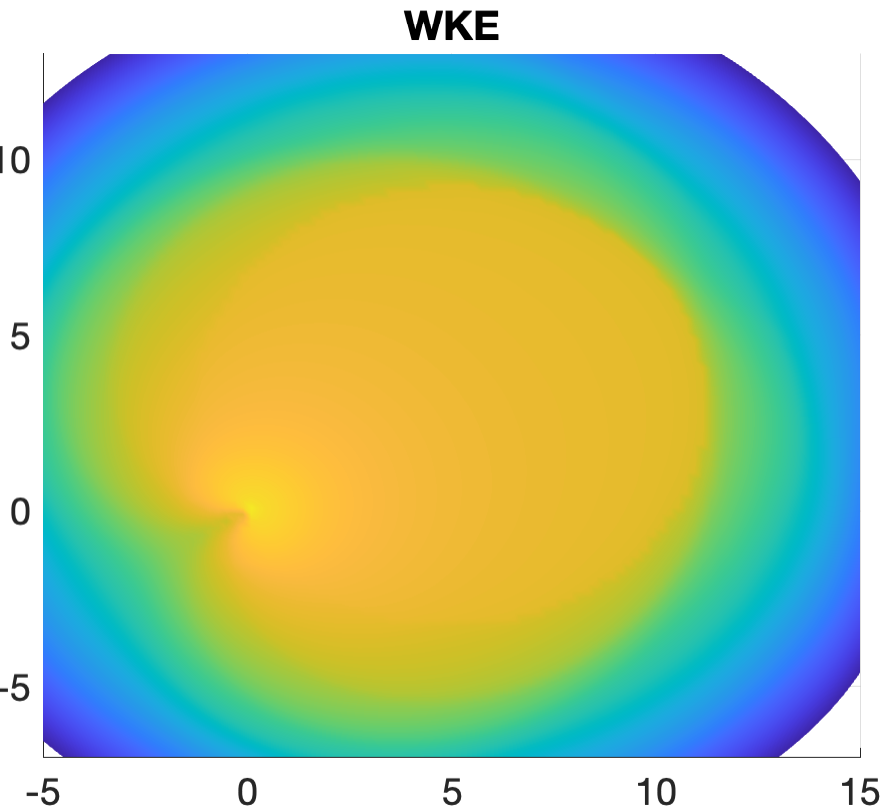}};
        \draw (5cm,0cm) node[anchor=south]{\includegraphics[trim=2.5cm 7.5cm 4cm 6.75cm,clip,width=4.5cm]{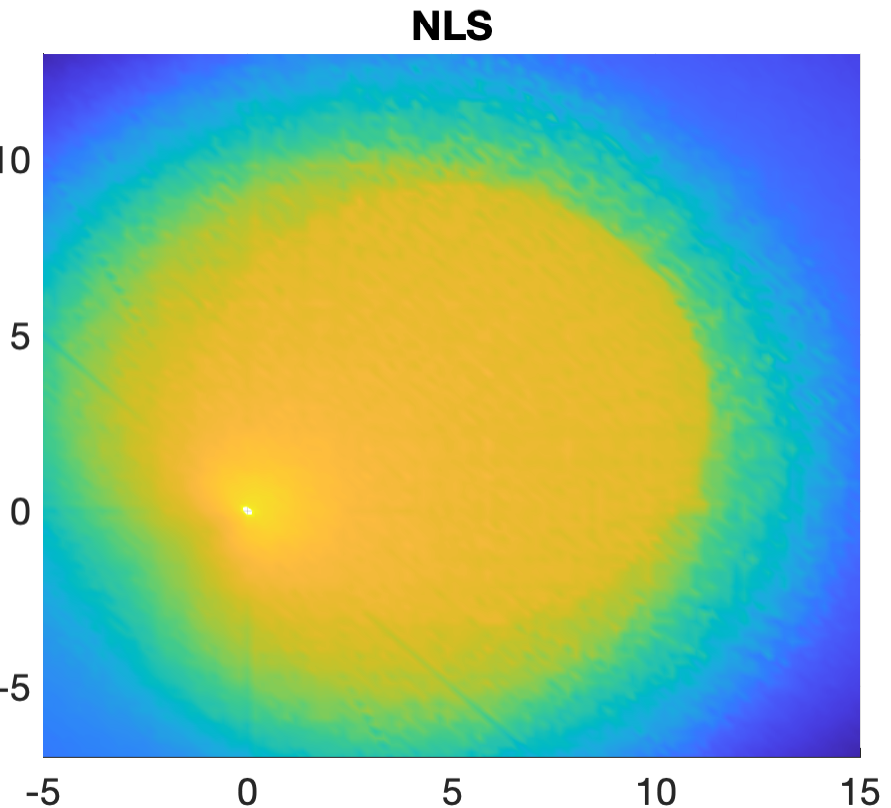}};
        \draw (10cm,0cm) node[anchor=south]{\includegraphics[trim=2.5cm 7.5cm 4cm 6.75cm,clip,width=4.5cm]{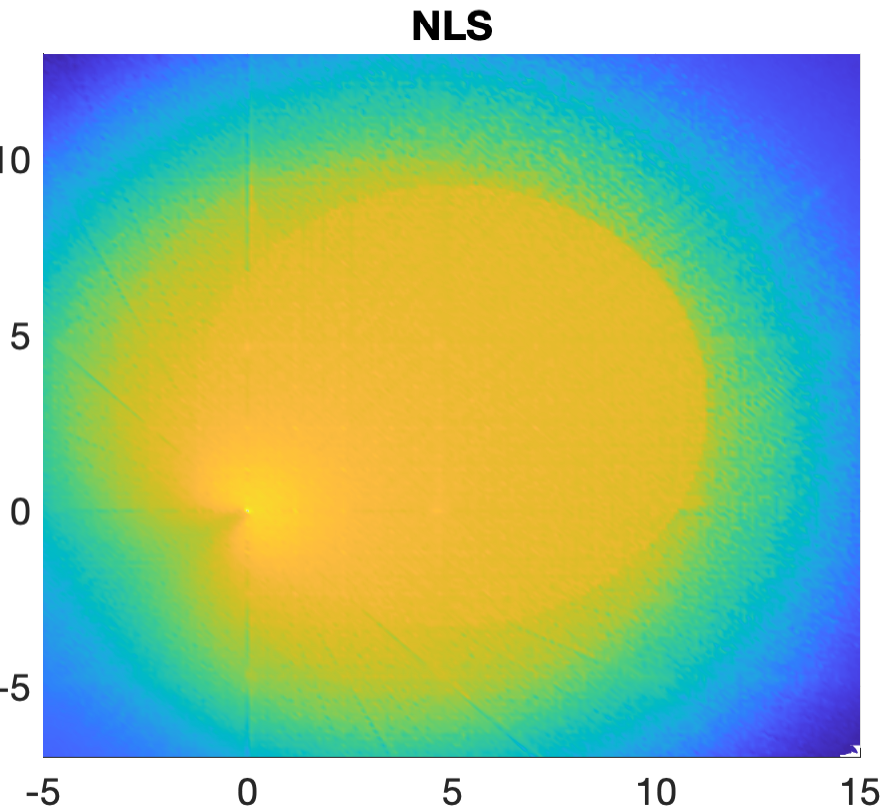}};
        \draw (12.5cm,.25cm) node[anchor=south]{\includegraphics[trim=17.5cm 7.75cm 2.5cm 6.75cm,clip,height=3.5cm]{images/AvgOffset_L64_N1000}};
        \draw (0.05,4.25)node[anchor=north,draw,fill=white,rounded corners]{\tiny WKE ($\tau=1$)};
        \draw (5.05,4.25)node[anchor=north,draw,fill=white,rounded corners]{\tiny NLS ($L=32$, $M_x=M_y=301$)};
        \draw (10.05,4.25)node[anchor=north,draw,fill=white,rounded corners]{\tiny NLS ($L=64$, $M_x=M_y=1001$)};
        \draw (0.05,0.1)node[anchor=north,fill=white]{\tiny $k_x$};
        \draw (5.05,0.1)node[anchor=north,fill=white]{\tiny $k_x$};
        \draw (10.05,0.1)node[anchor=north,fill=white]{\tiny $k_x$};
        \draw (-2.75,2.25)node[anchor=north,fill=white,rotate=90]{\tiny $k_y$};
      \end{tikzpicture}
  \caption{At left is the solution of the WKE at $\tau=1$ for the anisotropic top-hat. At center is a 180 member ensemble for the NLS with $L=32$ and $M_x=M_y=301$. At right is a 180 member ensemble for the NLS with $L=64$ and $M_x=M_y=1001$. All results are plotted using $\log_{10}$, which is also indicated by the color table. The similarity of the various solutions is clear, with a definite improvement by increasing $L$ and $\xv$-grid resolution.} 
  \label{fig:CompOffset2D}
  \end{center}
  \end{figure}

\begin{figure}[h]
  \begin{center}
    \begin{tikzpicture}[scale=.7]
      \useasboundingbox (-4.5,2.25) rectangle (18,7);
        \draw (0cm,0cm) node[anchor=south]{\includegraphics[width=5cm]{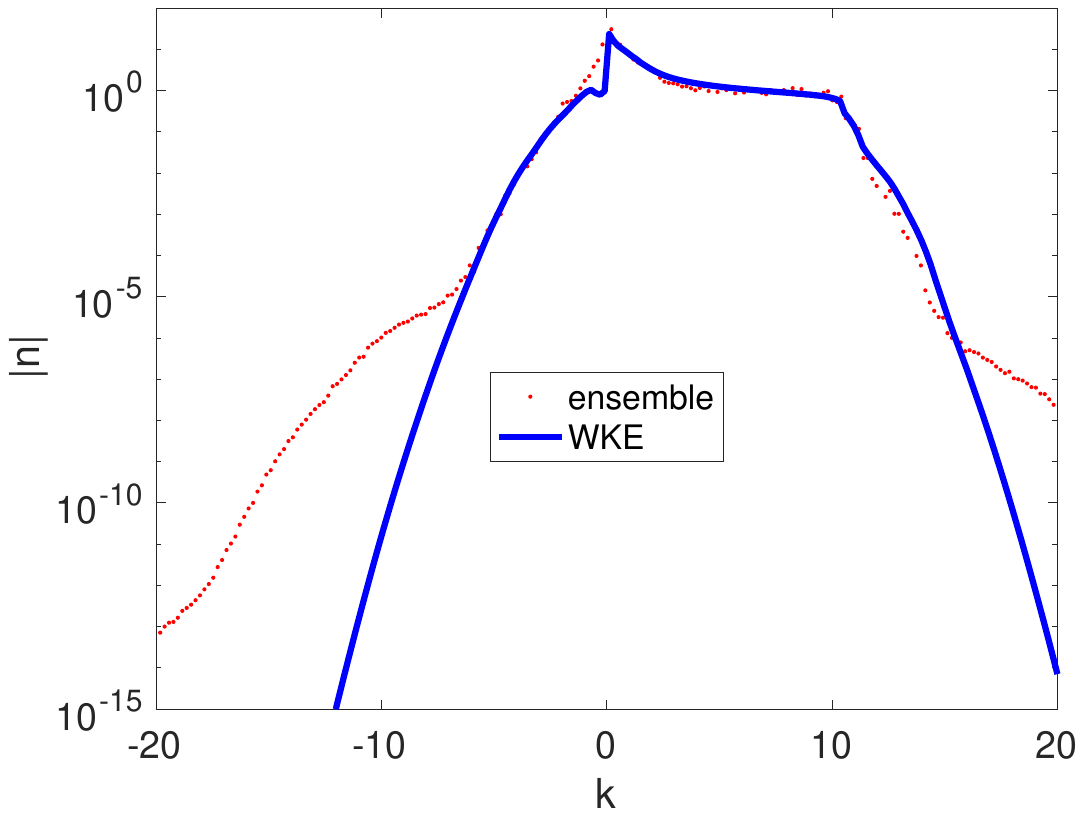}};
        \draw (7cm,0cm) node[anchor=south]{\includegraphics[width=5cm]{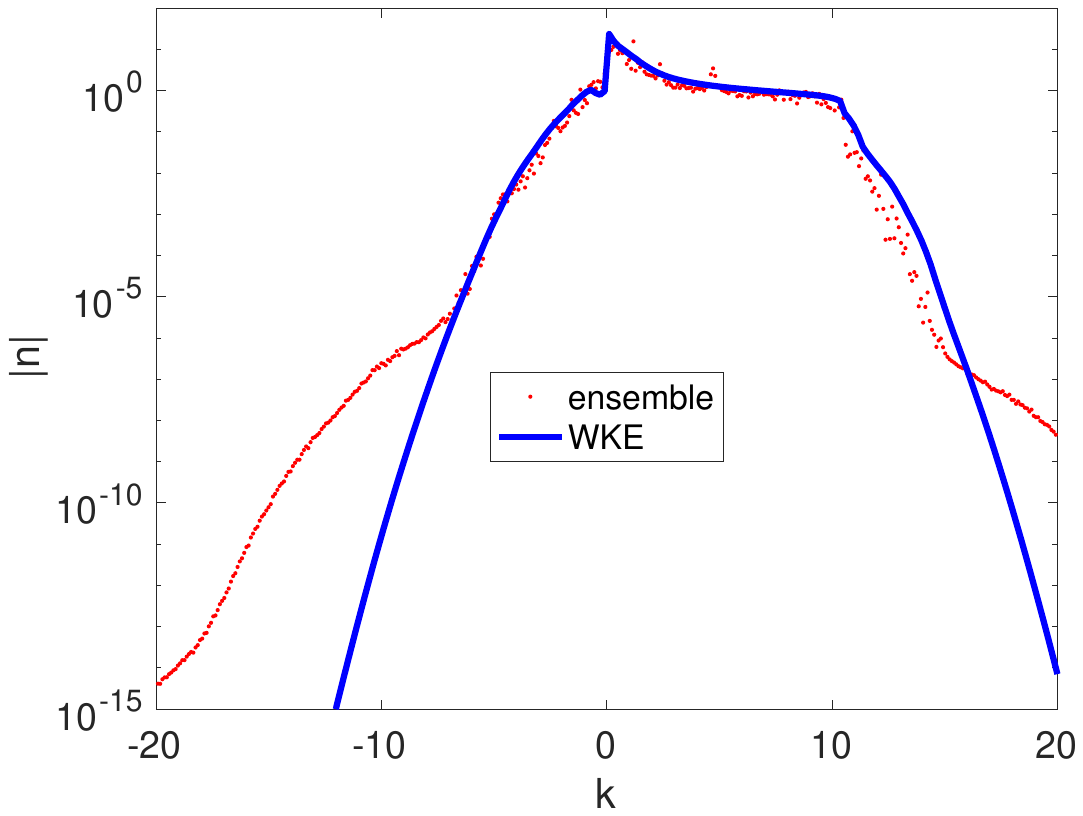}};
        \draw (14cm,0cm) node[anchor=south]{\includegraphics[width=5cm]{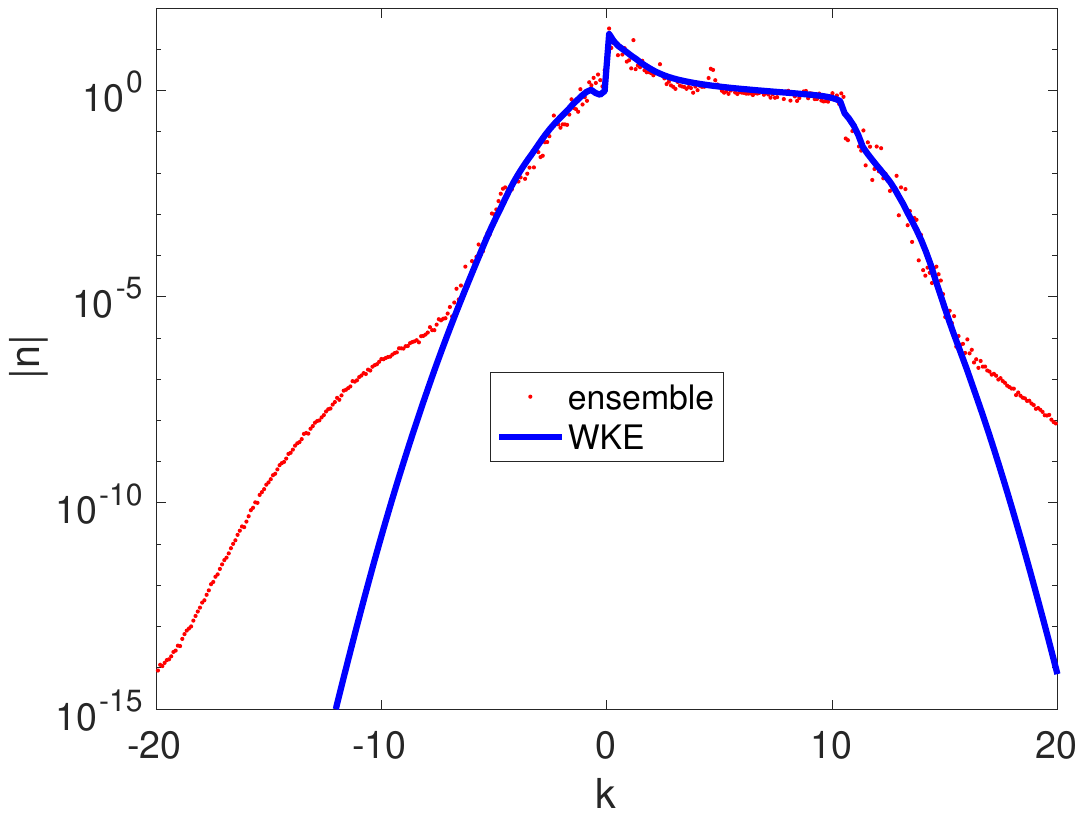}};
        \draw (0.1,7.55)node[anchor=north,draw,fill=white,rounded corners]{\tiny $L=32$, $M_x=M_y=301$};
        \draw (7.1,7.55)node[anchor=north,draw,fill=white,rounded corners]{\tiny $L=64$, $M_x=M_y=601$};
        \draw (14.1,7.55)node[anchor=north,draw,fill=white,rounded corners]{\tiny $L=64$, $M_x=M_y=1001$};
      \end{tikzpicture}
  \caption{Line slices along $\ky=0$ comparing the WKE solution to NLS ensemble averages. The general agreement is quite good for wave numbers in roughly $\kx\in[-8,12]$, with improvement for both increasing $L$ and $\xv$-grid resolution.} 
  \label{fig:CompOffset1D_y0}
  \end{center}
  \end{figure}

\begin{figure}[h]
  \begin{center}
    \begin{tikzpicture}[scale=.7]
      \useasboundingbox (-4.5,2.25) rectangle (18,7);
        \draw (0cm,0cm) node[anchor=south]{\includegraphics[width=5cm]{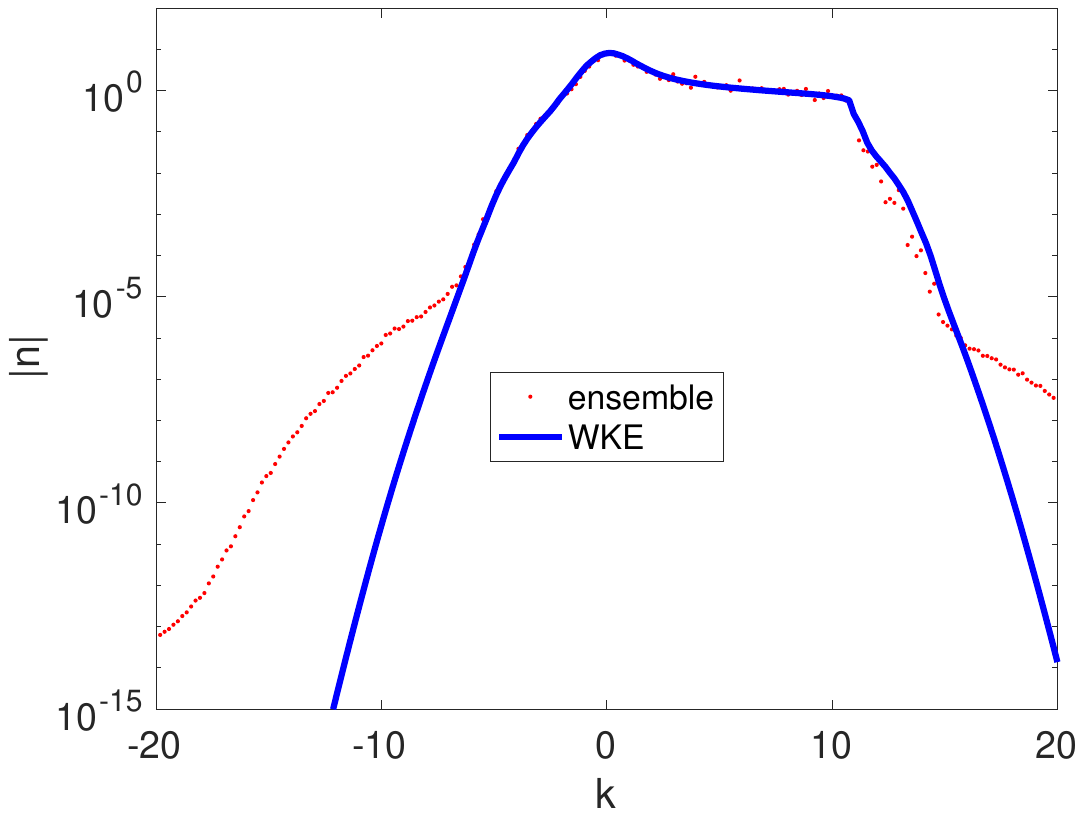}};
        \draw (7cm,0cm) node[anchor=south]{\includegraphics[width=5cm]{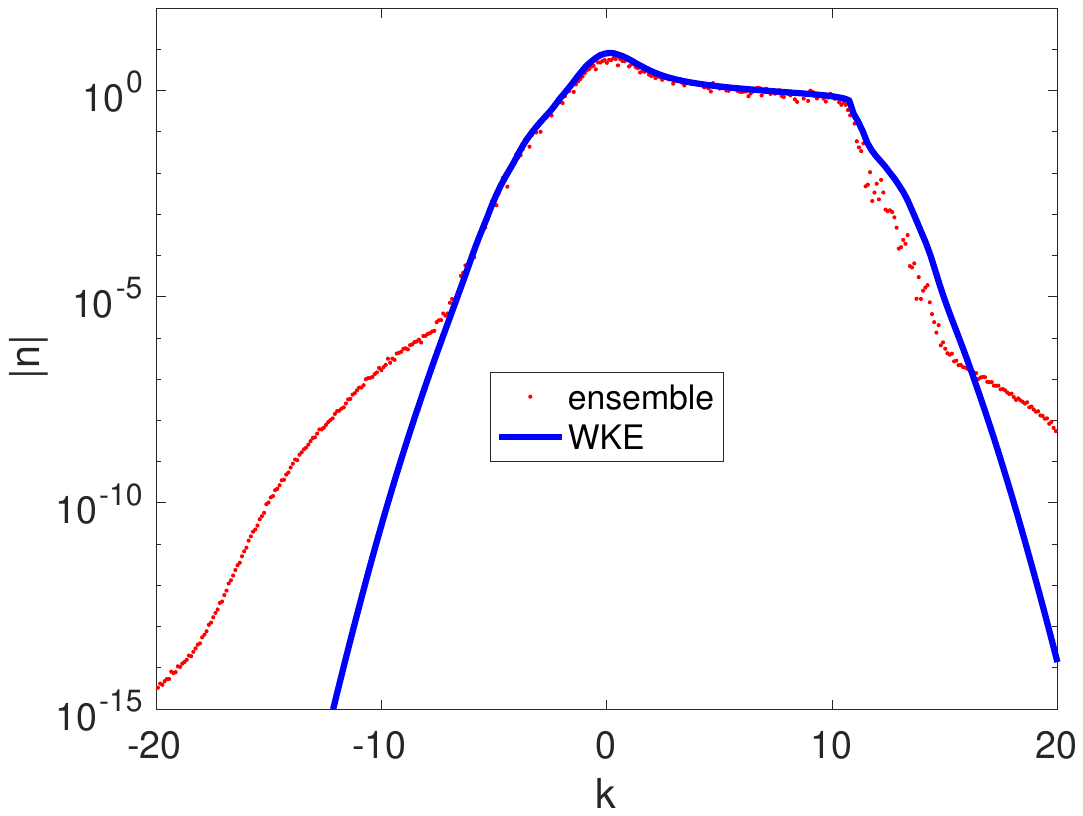}};
        \draw (14cm,0cm) node[anchor=south]{\includegraphics[width=5cm]{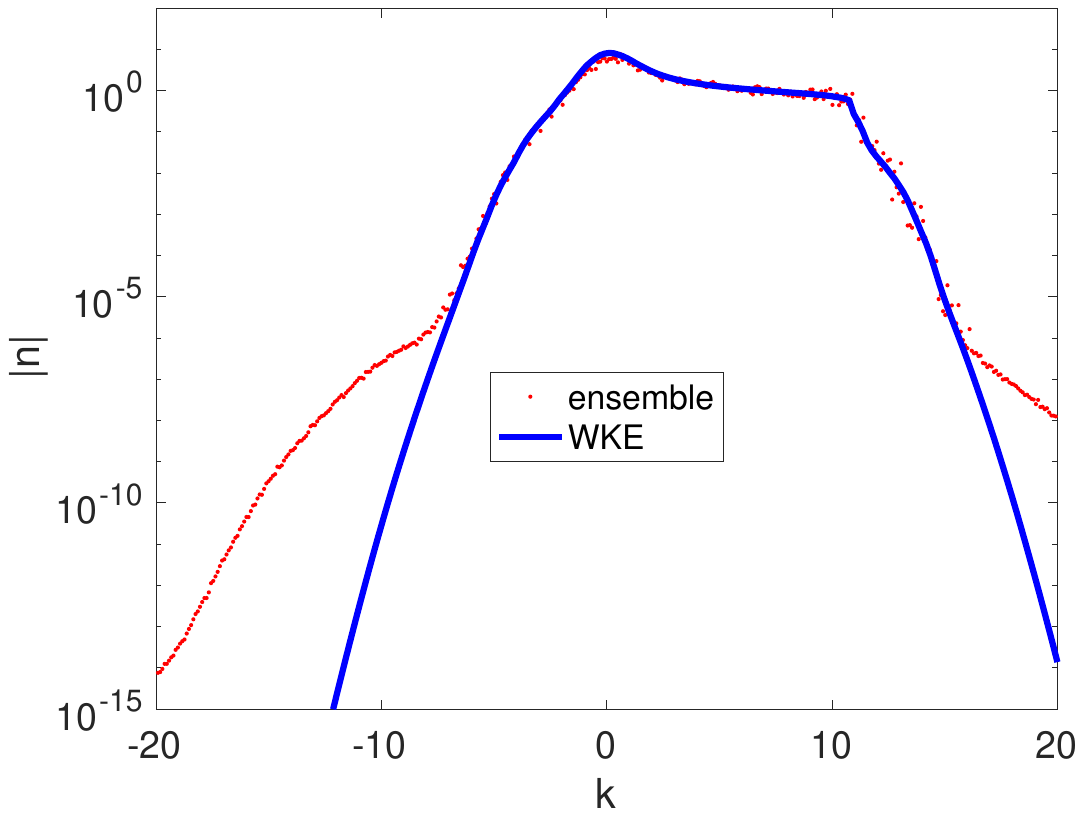}};
        \draw (0.1,7.55)node[anchor=north,draw,fill=white,rounded corners]{\tiny $L=32$, $M_x=M_y=301$};
        \draw (7.1,7.55)node[anchor=north,draw,fill=white,rounded corners]{\tiny $L=64$, $M_x=M_y=601$};
        \draw (14.1,7.55)node[anchor=north,draw,fill=white,rounded corners]{\tiny $L=64$, $M_x=M_y=1001$};
      \end{tikzpicture}
  \caption{Line slices along $\ky=2$ comparing the WKE solution to NLS ensemble averages. The general agreement remains quite good for wave numbers in roughly $\kx\in[-8,12]$, with improvement for both increasing $L$ and $\xv$-grid resolution.} 
  \label{fig:CompOffset1D_y1}
  \end{center}
  \end{figure}

Consider now the isotropic top-hat case which is given by initial condition with $\nk(0)=\mathcal{N}_{TH}(\kv)$ where $r_0=2\pi$ and $\kx_0=5$ and $\ky_0=3$. Figure~\ref{fig:CompOffset2D} shows computed results for the WKE~\eqref{eq:NLSWKE} with $\kx\in[-20,20]$, $\ky\in[-20,20]$, $N_{\kx}=200$, $N_{\ky}=202$, and computed to a final time $\tau=1$. For plotting, the domain has been restricted to $[-5,15]\times[-7,13]$, which is a box centered on  the initial top-hat and where the majority of the dynamics occur. Figure~\ref{fig:CompOffset2D} also shows the result of ensemble averages with 180 members for the NLS~\eqref{eq:NLS2D} for 2 cases, the first with $L=32$ and $M_x=M_y=301$, and the second with $L=64$ and $M_x=M_y=1001$. There is clear agreement between the WKE and the NLS ensembles, with a definite improvement by increasing $L$ and $\xv$-grid resolution.
Figures~\ref{fig:CompOffset1D_y0} and~\ref{fig:CompOffset1D_y1} present slices of the various solutions along the lines $\ky=0$ and $\ky=1$ respectively for $\kx\in[-20,20]$. As in the isotropic case, the WKE is computed for $\kx\in[-20,20]$, but the $\kx$ bounds for the ensemble NLS solutions are determined by $L$ and $M_x$, and are well outside $\pm20$. The figures presents ensemble NLS results for 3 cases; the first with $L=32$ and $M_x=M_y=301$, the second with $L=64$ and $M_x=M_y=601$, and the third with $L=64$ and $M_x=M_y=1001$. The general agreement between WKE and NLS approximations is seen to be quite good. One observation in comparing the limits of larger $L$ vs. increased $\xv$-grid resolution (which moves the high-$k$ boundaries further out) is that higher grid resolution has arguably more effect for these runs than further increases in $L$. Furthermore, the previously observed divergence for high $k$, which we attribute simultaneously to insufficient $\xv$-grid resolution and rounding effects, appears here as well.

\section{Conclusions}
\label{sec:conclusions}

In the present manuscript, we have discussed a new class of discretization techniques for time-dependent wave kinetic equations. Although the approach is rather general and could be applied to a wide class of WKEs, the discussion is grounded in a particular 2D nonlinear Schr{\"o}dinger equation and its associated WKE. The resonant manifold is approximated using a piecewise linear polynomial, and midpoint quadrature is applied to evaluate the collision integral. A method-of-lines approach is taken, and standard explicit Runge Kutta is used for time evolution. The efficacy of the scheme is demonstrated for model collision integrals in 1D and 2D, where second-order accuracy is demonstrated. For the full WKE, self-convergence studies for smooth and non-smooth initial conditions, where both isotropic and anisotropic cases show near second order convergence. Finally, we compare the approximate WKE solutions to ensemble averages of direct numerical approximations of the original NLS equations of motion. The convergence of the NLS ensembles and the WKE approximation is observed in the limit of large domains and well refined NLS discretizations. This latter comparison, i.e. solutions to the WKE and ensemble averages of NLS, is carefully crafted to respect the various limits and assumptions of the WKE derivation, which for completeness we present in an appendix. There are many potential avenues for future work following this line of research. At present, the methods naturally apply to equations with a quadratic nonlinear term. Relaxing this assumption is a natural next step that would open the door many more physical systems. At the same time, even in this restricted setting, the technique could fruitfully be applied to many models, e.g. capillary waves and gravity waves. Finally, there are numerical improvements that should be investigated, including higher-order accurate methods, implicit time stepping, and potentially to alternate domain truncation techniques. 

\bibliographystyle{elsart-num}
\bibliography{journal-ISI,WKE_Disc.bib}

\appendix

\section{Derivation of WKE for NLS model}
\label{sec:WKEDerive}
Here we present a derivation of the WKE~\eqref{eq:NLSWKE} from the governing equations of motion~\eqref{eq:NLS2D}, which are repeated here as 
\ba
  \imag\p_t u = (I-\Delta) u +u^2+2u\bar{u},
\ea
posed on a 2D periodic box with $\xv\in[0,L]^2$, where $\imag=\sqrt{-1}$ is the imaginary unit. Taking the DFT of both sides using $u=\frac{1}{L}\sum_{\jv}\hat{u}_{\kv}e^{\imag\kv\cdot x}$, $\bar{u}=\frac{1}{L}\sum_{\jv}\bar{\hat{u}}_{\kv}e^{-\imag\kv\cdot x}$, $\jv\in\Integer^2$, $\kv=2\pi\jv/L$, and $k=\|\kv\|$ yields
\bse
\ba
  \imag\p_t\hat{u}_{\kv} 
    & = (1+k^2)\hat{u}_{\kv}+\frac{1}{L}\sum_{\jv_1}\hat{u}_{\kv_1}\hat{u}_{\kv-\kv_1}+\frac{2}{L}\sum_{\jv_1}\hat{u}_{\kv_1}\bar{\hat{u}}_{\kv_1-\kv}.
\ea
Performing a change of variables with $\kv_2=\kv-\kv_1$ on the first convolution, and $\kv_2=\kv_1-\kv$ on the second, gives
\ba
  \imag\p_t\hat{u}_{\kv} 
    & = (1+k^2)\hat{u}_{\kv}+\frac{1}{L}\sum_{\jv_1}\sum_{\jv_2}\left[\hat{u}_{\kv_1}\hat{u}_{\kv_2}\Kron{\kv_1+\kv_2-\kv}+2\hat{u}_{\kv_1}\bar{\hat{u}}_{\kv_2}\Kron{\kv_1-\kv_2-\kv}\right],
\ea
where $\Kron{\cdot}$ is defined
\ba
  \Kron{x} = 
  \begin{cases}
    1 \qquad \hbox{if $x=0$}\\
    0 \qquad \hbox{else}.
  \end{cases}
\ea
Despite the slight abuse in nomenclature, $\Kron{\cdot}$ is sometimes referred to as a Kronecker delta function, and we will adopt this here. Note that square brackets are used to distinguish the Kronecker delta used here from the Dirac delta which will be used later and notated $\Dirac{\cdot}$. Using $\hat{u}_{\kv} = \epsilon a_\kv e^{-\imag\omega_\kv t}$, $\bar{\hat{u}}_{\kv} = \epsilon\bar{a}_\kv e^{\imag\omega_\kv t}$, and setting $\omega_\kv=1+k^2$ from the liner dispersion relation then yields
\ba
  &\imag\p_t\ak e^{-\imag\omega_\kv t}+\omega_\kv a_\kv e^{-\imag\omega_\kv t}\nonumber\\
   &\quad  = (1+k^2)a_\kv e^{-\imag\omega_\kv t}+\frac{\epsilon}{L}\sum_{\jv_1}\sum_{\jv_2}a_{\kv_1}a_{\kv_2}\Kron{\kv_1+\kv_2-\kv}e^{-\imag(\omega_{\kv_1}+\omega_{\kv_2})t}\nonumber\\
   &\hspace{1.25in}+\frac{\epsilon}{L}\sum_{\jv_1}\sum_{\jv_2}2a_{\kv_1}\bar{a}_{\kv_2}\Kron{\kv_1-\kv_2-\kv}e^{-\imag(\omega_{\kv_1}-\omega_{\kv_2})t},\\
  &\imag\p_t\ak 
    = \frac{\epsilon}{L}\sum_{\jv_1}\sum_{\jv_2}a_{\kv_1}a_{\kv_2}\Kron{\kv_1+\kv_2-\kv}e^{-\imag(\omega_{\kv_1}+\omega_{\kv_2}-\omega_\kv)t}\nonumber\\
    &\hspace{1.25in}+\frac{\epsilon}{L}\sum_{\jv_1}\sum_{\jv_2}2a_{\kv_1}\bar{a}_{\kv_2}\Kron{\kv_1-\kv_2-\kv}^{-\imag(\omega_{\kv_1}-\omega_{\kv_2}-\omega_\kv)t}. \label{eq:akBasic}
\ea
\ese

To simplify the presentation, we adopt some simplifying notation as below, where  $e$, $f$, and $g$ are integers
\bse
\ba
  \mu & = \frac{\epsilon}{L},\\
  \del{ef}{g} & = \Kron{\kv_e+\kv_f-\kv_g},\\
  \del{e}{fg} & = \Kron{\kv_e-\kv_f-\kv_g},\\
  \om{ef}{g} & = \omega_{\kv_e}+\omega_{\kv_f}-\omega_{\kv_g},\\
  \om{e}{fg} & = \omega_{\kv_e}-\omega_{\kv_f}-\omega_{\kv_g},\\
  \al{ef}{g} & = a_{\kv_e}a_{\kv_f}\bar{a}_{\kv_g}\del{ef}{g},\\
  \al{e}{fg} & = a_{\kv_e}\bar{a}_{\kv_f}\bar{a}_{\kv_g}\del{e}{fg}.
\ea
In this notation, for $\delta$ and $\Omega$, the lower indices indicate a $+$ sign, while the upper indices indicate a $-$ sign. For $\alpha$, the lower indices indicate no conjugate, while the upper indices indicate a conjugated quantity. We will also adopt the convention that an explicit indicated zero index should be simply dropped so that, for example $\del{e0}{g} = \Kron{\kv_e+\kv-\kv_g}$. Further, a red circle will be used in conjunction with the $\alpha$ notation to indicate  an entry of $a$ that is left out, e.g. 
\ba
  \al{ef}{\no{g}} & = a_{\kv_e}a_{\kv_f}\del{ef}{g},\\
  \al{e\no{f}}{g} & = a_{\kv_e}\bar{a}_{\kv_g}\del{ef}{g},\\
  \al{\no{e}}{f\ko} & = \bar{a}_{\kv_f}\akb\del{e}{f0},\\
  \al{e}{f\no{\ko}} & = a_{\kv_e}\bar{a}_{\kv_f}\del{e}{f0}.
\ea
Finally for brevity, repeated sums will be represented using a single summation, e.g.
\ba
  \sum_{e,f} & = \sum_{\jv_e}\sum_{\jv_f},\\
  \sum_{e,f,g,h} & = \sum_{\jv_e}\sum_{\jv_f}\sum_{\jv_g}\sum_{\jv_h}.
\ea
\ese
With this notation in place, Equation~\eqref{eq:akBasic} and its conjugate can be written more compactly as
\bse
\ba
  \p_t \ak  & = -\imag\mu\sum_{1,2}\left[\al{12}{\no{\ko}}e^{-\imag\om{12}{\ko}t}+2\al{1}{2\no{\ko}}e^{-\imag\om{1}{2\ko}t}\right],\\
  \p_t \akb &  = \imag\mu\sum_{1,2}\left[\al{\no{\ko}}{12}e^{-\imag\om{\ko}{12}t}+2\al{2\no\ko}{1}e^{-\imag\om{2\ko}{1}t}\right].
\ea
\ese
For the WKE, we compute 
\bse
\ba
  \p_t|\ak|^2 
    & = \p_t\ak\akb+\ak\p_t\akb,\\
    & = \p_t\ak\akb+c.c.,\\
    & = -\imag\mu\sum_{1,2}\left[\al{12}{\no{\ko}}\akb e^{-\imag\om{12}{\ko}t}+2\al{1}{2\no{\ko}}\akb e^{-\imag\om{1}{2\ko}t}\right] +c.c.,\\
    & = -\imag\mu\sum_{1,2}\left[\al{12}{\ko}e^{-\imag\om{12}{\ko}t}+2\al{1}{2\ko}e^{-\imag\om{1}{2\ko}t}\right] +c.c..\label{eq:ak2}
\ea
\ese
Denoting the Laplace transform as
\[
  \mathcal{L}[f](s) = \int_0^{\infty}e^{-st}f(t)\,dt,
\]
and applying to~\eqref{eq:ak2} then gives
\ba
  s\mathcal{L}[|\ak|^2](s)-|\ak(0)|^2 = -\imag\mu\sum_{1,2}\left[\int_0^{\infty}\al{12}{\ko}e^{-(s+\imag\om{12}{\ko})t}+2\al{1}{2\ko}e^{-(s+\imag\om{1}{2\ko})t}\, dt\right] +c.c..
\ea
Integrating by parts gives
\ba
  & s\mathcal{L}[|\ak|^2](s)-|\ak(0)|^2 \nonumber\\ 
  &  \qquad = -\imag\mu\sum_{1,2}\left[
    \frac{\al{12}{\ko}(0)}{s+\imag\om{12}{\ko}}
    +\frac{2\al{1}{2\ko}(0)}{s+\imag\om{1}{2\ko}}
    +\int_0^{\infty}\frac{\p_t\al{12}{\ko}}{s+\imag\om{12}{\ko}}e^{-(s+\imag\om{12}{\ko})t}
    +2\frac{\p_t\al{1}{2\ko}}{s+\imag\om{1}{2\ko}}e^{-(s+\imag\om{1}{2\ko})t}\, dt\right] +c.c..
    \label{eq:IBP1}
\ea
The definition of the $\alpha$'s then gives 
\bse
\ba
   \p_t\al{12}{\ko} 
    & = \left
      (\p_t\bar{a}_{\kv}a_{\kv_1}a_{\kv_2}
      +\bar{a}_{\kv}\p_ta_{\kv_1}a_{\kv_2}
      +\bar{a}_{\kv}a_{\kv_1}\p_ta_{\kv_2}
    \right)\del{12}{\ko}, \medskip \\
  \p_t\al{1}{2\ko} 
    & = \left
      (\p_t\bar{a}_{\kv}a_{\kv_1}\bar{a}_{\kv_2}
      +\bar{a}_{\kv}\p_ta_{\kv_1}\bar{a}_{\kv_2}
      +\bar{a}_{\kv}a_{\kv_1}\p_t\bar{a}_{\kv_2}
    \right)\del{1}{2\ko}.
\ea
\ese
For convenience recall 
\bse
\ba
  \p_t \akb &  = \imag\mu\sum_{3,4}\left[\al{\no{\ko}}{34}e^{-\imag\om{\ko}{34}t}+2\al{4\no\ko}{3}e^{-\imag\om{4\ko}{3}t}\right],\\
  \p_t a_{\kv_1} & = -\imag\mu\sum_{3,4}\left[\al{34}{\no 1}e^{-\imag\om{34}{1}t}+2\al{3}{4\no1}e^{-\imag\om{3}{41}t}\right],\\
  \p_t a_{\kv_2} & = -\imag\mu\sum_{3,4}\left[\al{34}{\no 2}e^{-\imag\om{34}{2}t}+2\al{3}{4\no2}e^{-\imag\om{3}{42}t}\right],\\
  \p_t \bar{a}_{\kv_2} & = \imag\mu\sum_{3,4}\left[\al{\no2}{34}e^{-\imag\om{2}{34}t}+2\al{4\no2}{3}e^{-\imag\om{42}{3}t}\right],\\
  \p_t \ak  & = -\imag\mu\sum_{3,4}\left[\al{34}{\no{\ko}}e^{-\imag\om{34}{\ko}t}+2\al{3}{4\no{\ko}}e^{-\imag\om{3}{4\ko}t}\right],
\ea
\ese
to obtain
\bse
\ba
  \p_t\al{12}{\ko} =  
    &+\imag\mu\sum_{3,4}
       \al{12}{\no\ko}\left(\al{\no{\ko}}{34}e^{-\imag\om{\ko}{34}t}+2\al{4\no{\ko}}{3}e^{-\imag\om{4\ko}{3}t}\right)\nonumber\\       
    &-\imag\mu\sum_{3,4}
     \al{\no12}{\ko}\left(\al{34}{\no 1}e^{-\imag\om{34}{1}t}+2\al{3}{4\no1}e^{-\imag\om{3}{41}t}\right)\nonumber\\
     &-\imag\mu\sum_{3,4}\al{1\no2}{\ko}\left(\al{34}{\no 2}e^{-\imag\om{34}{2}t}+2\al{3}{4\no2}e^{-\imag\om{3}{42}t}\right),\\
     \p_t\al{1}{2\ko} = 
       & +\imag\mu\sum_{3,4}\al{1}{2\no\ko}\left(\al{\no{\ko}}{34}e^{-\imag\om{\ko}{34}t}+2\al{4\no{\ko}}{3}e^{-\imag\om{4\ko}{3}t}\right)\nonumber\\
       & -\imag\mu\sum_{3,4}\al{\no1}{2\ko}\left(\al{34}{\no 1}e^{-\imag\om{34}{1}t}+2\al{3}{4\no1}e^{-\imag\om{3}{41}t}\right)\nonumber\\
       & -\imag\mu\sum_{3,4}\al{1}{\no2\ko}\left(\al{\no2}{34}e^{-\imag\om{2}{34}t}+2\al{4\no2}{3}e^{-\imag\om{42}{3}t}\right).
\ea
\ese
Putting this all together in~\eqref{eq:IBP1} gives
\ba
  s\mathcal{L}[|\ak|^2](s)-|\ak(0)|^2 &= -\imag\mu\sum_{1,2}\left[
    \frac{\al{12}{\ko}(0)}{s+\imag\om{12}{\ko}}
    +\frac{2\al{1}{2\ko}(0)}{s+\imag\om{1}{2\ko}}\right]\nonumber\\
    &
    -(\imag\mu)^2\sum_{1,2,3,4}\int_0^{\infty}\frac{    
      \al{12}{\no\ko}\al{\no{\ko}}{34}e^{-\imag\om{\ko}{34}t}+2\al{12}{\no\ko}\al{4\no{\ko}}{3}e^{-\imag\om{4\ko}{3}t}
      }{s+\imag\om{12}{\ko}}\, dt\nonumber\\
    &
    +(\imag\mu)^2\sum_{1,2,3,4}\int_0^{\infty}\frac{    
      \al{\no12}{\ko}\al{34}{\no 1}e^{-\imag\om{34}{1}t}+2\al{\no12}{\ko}\al{3}{4\no1}e^{-\imag\om{3}{41}t}
      }{s+\imag\om{12}{\ko}}\, dt\nonumber\\
    &
    +(\imag\mu)^2\sum_{1,2,3,4}\int_0^{\infty}\frac{    
      \al{1\no2}{\ko}\al{34}{\no 2}e^{-\imag\om{34}{2}t}+2\al{1\no2}{\ko}\al{3}{4\no2}e^{-\imag\om{3}{42}t}
      }{s+\imag\om{12}{\ko}}\, dt\nonumber\\
    &
    -2(\imag\mu)^2\sum_{1,2,3,4}\int_0^{\infty}\frac{    
      \al{1}{2\no\ko}\al{\no{\ko}}{34}e^{-\imag\om{\ko}{34}t}+2\al{1}{2\no\ko}\al{4\no{\ko}}{3}e^{-\imag\om{4\ko}{3}t}
      }{s+\imag\om{1}{2\ko}}\, dt\nonumber \\
    &
    +2(\imag\mu)^2\sum_{1,2,3,4}\int_0^{\infty}\frac{    
      \al{\no1}{2\ko}\al{34}{\no 1}e^{-\imag\om{34}{1}t}+2\al{\no1}{2\ko}\al{3}{4\no1}e^{-\imag\om{3}{41}t}
      }{s+\imag\om{1}{2\ko}}\, dt\nonumber\\
    &
    -2(\imag\mu)^2\sum_{1,2,3,4}\int_0^{\infty}\frac{    
      \al{1}{\no2\ko}\al{\no2}{34}e^{-\imag\om{2}{34}t}+2\al{1}{\no2\ko}\al{4\no2}{3}e^{-\imag\om{42}{3}t}
      }{s+\imag\om{1}{2\ko}}\, dt + c.c. .
    \label{eq:IBP1Expand}
\ea
We now apply a second integration by parts which will return a set of jump terms (which are the crux of the derivation and so we will write), and a set of integral terms whose integrand is cubic in $\alpha$'s. The latter point is critical since upon taking the expectation with the assumption of random phases all cubic terms will vanish. As a result the integral portions will be suppressed to avoid unnecessary confusion and complication of presentation. After performing the second IBP, and understanding that all jump terms are evaluated at $t=0$, we have
\ba
  s\mathcal{L}[|\ak|^2](s)-|\ak(0)|^2 &= -\imag\mu\sum_{1,2}
    \frac{\al{12}{\ko}}{s+\imag\om{12}{\ko}}
    +\frac{2\al{1}{2\ko}}{s+\imag\om{1}{2\ko}}\nonumber\\
    &
    -(\imag\mu)^2\sum_{1,2,3,4}
      \frac{\al{12}{\no\ko}\al{\no{\ko}}{34}}{(s+\imag\om{12}{\ko})(s+\imag\om{12}{\ko}+\imag\om{\ko}{34})}+
      \frac{2\al{12}{\no\ko}\al{4\no{\ko}}{3}}{(s+\imag\om{12}{\ko})(s+\imag\om{12}{\ko}+\imag\om{4\ko}{3})}
    \nonumber\\
    &
    +(\imag\mu)^2\sum_{1,2,3,4}
      \frac{\al{\no12}{\ko}\al{34}{\no 1}}{(s+\imag\om{12}{\ko})(s+\imag\om{12}{\ko}+\imag\om{34}{1})}+
      \frac{2\al{\no12}{\ko}\al{3}{4\no1}}{(s+\imag\om{12}{\ko})(s+\imag\om{12}{\ko}+\imag\om{3}{41})} \nonumber\\
    &
    +(\imag\mu)^2\sum_{1,2,3,4}
      \frac{\al{1\no2}{\ko}\al{34}{\no 2}}{(s+\imag\om{12}{\ko})(s+\imag\om{12}{\ko}+\om{34}{2})}+
      \frac{2\al{1\no2}{\ko}\al{3}{4\no2}}{(s+\imag\om{12}{\ko})(s+\imag\om{12}{\ko}+\om{3}{42})}\nonumber\\
    &
    -2(\imag\mu)^2\sum_{1,2,3,4}
      \frac{\al{1}{2\no\ko}\al{\no{\ko}}{34}}{(s+\imag\om{1}{2\ko})(s+\imag\om{1}{2\ko}+\imag\om{\ko}{34})}+
      \frac{2\al{1}{2\no\ko}\al{4\no{\ko}}{3}}{(s+\imag\om{1}{2\ko})(s+\imag\om{1}{2\ko}+\om{4\ko}{3})}\nonumber\\
    &
    +2(\imag\mu)^2\sum_{1,2,3,4}
      \frac{\al{\no1}{2\ko}\al{34}{\no 1}}{(s+\imag\om{1}{2\ko})(s+\imag\om{1}{2\ko}+\imag\om{34}{1})}+
      \frac{2\al{\no1}{2\ko}\al{3}{4\no1}}{(s+\imag\om{1}{2\ko})(s+\imag\om{1}{2\ko}+\imag\om{3}{41})}\nonumber\\
    &
    -2(\imag\mu)^2\sum_{1,2,3,4}
      \frac{\al{1}{\no2\ko}\al{\no2}{34}}{(s+\imag\om{1}{2\ko})(s+\imag\om{1}{2\ko}+\om{2}{34})}+
      \frac{2\al{1}{\no2\ko}\al{4\no2}{3}}{(s+\imag\om{1}{2\ko})(s+\imag\om{1}{2\ko}+\om{42}{3})}\nonumber\\
    &+\int_0^{\infty}\ldots\, dt +c.c. .
    \label{eq:4Suma}
\ea
By symmetry the 2nd and 3rd terms can be combined to yield
\ba
  s\mathcal{L}[|\ak|^2](s)-|\ak(0)|^2 &= -\imag\mu\sum_{1,2}
    \frac{\al{12}{\ko}}{s+\imag\om{12}{\ko}}
    +\frac{2\al{1}{2\ko}}{s+\imag\om{1}{2\ko}}\nonumber\\
    &
    +\mu^2\sum_{1,2,3,4}
      \frac{\al{12}{\no\ko}\al{\no{\ko}}{34}}{(s+\imag\om{12}{\ko})(s+\imag\om{12}{\ko}+\imag\om{\ko}{34})}+
      \frac{2\al{12}{\no\ko}\al{4\no{\ko}}{3}}{(s+\imag\om{12}{\ko})(s+\imag\om{12}{\ko}+\imag\om{4\ko}{3})}
    \nonumber\\
    &
    -2\mu^2\sum_{1,2,3,4}
      \frac{\al{\no12}{\ko}\al{34}{\no 1}}{(s+\imag\om{12}{\ko})(s+\imag\om{12}{\ko}+\imag\om{34}{1})}+
      \frac{2\al{\no12}{\ko}\al{3}{4\no1}}{(s+\imag\om{12}{\ko})(s+\imag\om{12}{\ko}+\imag\om{3}{41})} \nonumber\\
    &
    +2\mu^2\sum_{1,2,3,4}
      \frac{\al{1}{2\no\ko}\al{\no{\ko}}{34}}{(s+\imag\om{1}{2\ko})(s+\imag\om{1}{2\ko}+\imag\om{\ko}{34})}+
      \frac{2\al{1}{2\no\ko}\al{4\no{\ko}}{3}}{(s+\imag\om{1}{2\ko})(s+\imag\om{1}{2\ko}+\om{4\ko}{3})}\nonumber\\
    &
    -2\mu^2\sum_{1,2,3,4}
      \frac{\al{\no1}{2\ko}\al{34}{\no 1}}{(s+\imag\om{1}{2\ko})(s+\imag\om{1}{2\ko}+\imag\om{34}{1})}+
      \frac{2\al{\no1}{2\ko}\al{3}{4\no1}}{(s+\imag\om{1}{2\ko})(s+\imag\om{1}{2\ko}+\imag\om{3}{41})}\nonumber\\
    &
    +2\mu^2\sum_{1,2,3,4}
      \frac{\al{1}{\no2\ko}\al{\no2}{34}}{(s+\imag\om{1}{2\ko})(s+\imag\om{1}{2\ko}+\om{2}{34})}+
      \frac{2\al{1}{\no2\ko}\al{4\no2}{3}}{(s+\imag\om{1}{2\ko})(s+\imag\om{1}{2\ko}+\om{42}{3})}\nonumber\\
    &+\int_0^{\infty}\ldots\, dt +c.c. .
    \label{eq:IBP2}
\ea

As indicated, the idea will now be to take the expectation over assumed random phases. Denoting the expectation as $<\cdot>$ we note that
\bse
\label{eq:expectations}
\ba
  <a_{\kv_i}> & = 0,\\
  <a_{\kv_i}a_{\kv_j}> & = 0,\\
  <a_{\kv_i}\bar{a}_{\kv_j}> & = \begin{cases}
    0 & i\ne j\\
    n_{\kv_i} & i=j,
  \end{cases}\\
  <a_{\kv_i}a_{\kv_j}a_{\kv_l}> & = 0,\\
  <a_{\kv_i}a_{\kv_j}\bar{a}_{\kv_l}> & = 0, \label{eq:cubic1}\\
  <a_{\kv_i}\bar{a}_{\kv_j}a_{\kv_{\kv_l}}> & = 0,\\
  <a_{\kv_i}\bar{a}_{\kv_j}\bar{a}_{\kv_l}> & = 0,\\
  <\bar{a}_{\kv_i}\bar{a}_{\kv_j}\bar{a}_{\kv_l}> & = 0\label{eq:cubic4}.
\ea
The cubic terms indicated in~\eqref{eq:cubic1}--\eqref{eq:cubic4} occur in the unwritten integral terms in~\eqref{eq:IBP2}, as well as imply that
\ba
  <\alpha_{ij}^l> & = 0,\\
  <\alpha_{i}^{jl}> & = 0.
\ea
\ese
The expectations for $a$'s in~\eqref{eq:expectations} can be converted into relevant expectations for quadratic terms in $\alpha$. First consider
\ba
  \sum_{ef}<\al{ij}{\no l}\al{\no l}{ef}> 
    & = \sum_{ef}\left[\del{ij}{l}\del{le}{f}<a_{\kv_i}a_{\kv_j}\bar{a}_{\kv_e}\bar{a}_{\kv_f}>\right],\label{eq:resonant1}\\
    & = \del{ij}{l}\sum_{ef}\left[\del{le}{f}<a_{\kv_i}\bar{a}_{\kv_e}><a_{\kv_j}\bar{a}_{\kv_f}>\right]
        +\del{ij}{l}\sum_{ef}\left[\del{le}{f}<a_{\kv_i}\bar{a}_{\kv_f}><a_{\kv_j}\bar{a}_{\kv_e}>\right],\nonumber\\
    & = \del{ij}{l}n_{\kv_i}n_{\kv_j}+\del{ij}{l}n_{\kv_i}n_{\kv_j},\nonumber\\
    & = 2\del{ij}{l}n_{\kv_i}n_{\kv_j}.\nonumber
\ea
Note that nonzero contributions in~\eqref{eq:resonant1} occur when either $e=i$ and $f=j$, or $f=i$ and $e=j$. For $e=i$ and $f=j$, the corresponding term in the denominator in~\eqref{eq:IBP2} will contain 
\ba
  s+\imag\om{ij}{l}+\imag\om{l}{ef} & = 
    s+\imag(\om{ij}{l}+\om{l}{ij}) = s,
\ea
and similarly for $f=i$ and $e=j$.

The other non-zero terms have a richer matching behavior and we must consider the entire 4-way sum. Note that we will eventually find this richness is essentially meaningless, but for completeness we will cary it through. The details depend superficially on where $\kv$ appears, so first consider 
\ba
  \sum_{ijef}<\al{\no i j}{0}\al{e}{f \no i}> 
    & = \sum_{ijef}\left[\del{ij}{0}\del{e}{fi}<a_{\kv_j}\bar{a}_{\kv}a_{\kv_e}\bar{a}_{\kv_f}>\right],\label{eq:resonant2}\\
    & = \sum_{ijef}\left[\del{ij}{0}\del{e}{fi}<a_{\kv_j}\bar{a}_{\kv}><a_{\kv_e}\bar{a}_{\kv_f}>\right]
        +\sum_{ijef}\left[\del{ij}{0}\del{e}{fi}<a_{\kv_j}\bar{a}_{\kv_f}><a_{\kv_e}\bar{a}_{\kv}>\right],\nonumber\\
    & = \sum_{e}n_{\kv}n_{\kv_e}+\sum_{ij}\del{ij}{0}n_{\kv}n_{\kv_j},\nonumber\\
    & = \sum_{j}n_{\kv}n_{\kv_j}+\sum_{ij}\del{ij}{0}n_{\kv}n_{\kv_j},\nonumber\\
    & = \sum_{ij} n_{\kv}n_{\kv_j}(\del{i}{}+\del{ij}{0}).\nonumber
\ea
Nonzero contributions in~\eqref{eq:resonant2} occur in 2 distinct ways. In the first, $j=0$, and $e=f$, and the corresponding term in the denominator of~\eqref{eq:IBP2} will be
\ba
  (s+\imag\om{ij}{0})(s+\imag\om{ij}{0}+\imag\om{e}{fi}) & = 
    (s+\imag\om{i0}{0})(s+\imag(\om{i0}{0}+\om{e}{ei})) = (s+\imag\om{i}{})s.
\ea
The second type of match occurs when $j=f$, and $e=0$, where the corresponding term in the denominator of~\eqref{eq:IBP2} will be
\ba
  (s+\imag\om{ij}{0})(s+\imag\om{ij}{0}+\imag\om{e}{fi}) & = 
    (s+\imag\om{ij}{0})(s+\imag(\om{ij}{0}+\om{0}{0j})) = (s+\imag\om{ij}{0})s.
\ea
The other way that matching can occur leads to  
\bse
\ba
  \sum_{ijef}<\al{i}{j\no 0}\al{f\no 0}{e}> 
    & = \sum_{ijef}\left[\del{i}{j0}\del{f0}{e}<a_{\kv_i}\bar{a}_{\kv_j}a_{\kv_f}\bar{a}_{\kv_e}>\right],\\
    & = \sum_{ijef}\left[\del{i}{j0}\del{f0}{e}<a_{\kv_i}\bar{a}_{\kv_j}><a_{\kv_f}\bar{a}_{\kv_e}>\right]
        + \sum_{ijef}\left[\del{i}{j0}\del{f0}{e}<a_{\kv_i}\bar{a}_{\kv_e}a_{\kv_f}\bar{a}_{\kv_j}>\right],\nonumber\\
    & = \sum_{ie}n_{\kv_i}n_{\kv_e}\del{}{0}+\sum_{ij}\del{i}{j0}n_{\kv_i}n_{\kv_j},\nonumber\\
    & = \sum_{ij}n_{\kv_i}n_{\kv_j}(\del{}{0}+\del{i}{j0}).\nonumber
\ea
\ese
The end result of this matching is similar to the prior case, and the corresponding terms in the denominator follow similarly. Expectations of all other quadratic terms in $\alpha$ vanish so
\ba
 \sum_{ef}<\al{ij}{\no l}\al{e\no l}{f}> 
  = \sum_{ef}<\al{j\no l}{i}\al{ef}{\no l}> 
  = \sum_{ef}<\al{\no l}{ij}\al{e}{f\no l}>
  = \sum_{ef}<\al{i}{j\no l}\al{\no l}{ef}> = 0.
\ea

Using the expectation identities, explicitly adding in the complex conjugate term, and dropping the low-order terms yields
\ba
  s\mathcal{L}[\nk](s)-\nk(0) \approx &\mu^2\sum_{1,2}\Bigg[
    \frac{2n_{\kv_1}n_{\kv_2}\del{12}{\ko}}{(s+\imag\om{12}{\ko})s}
    -\frac{4\nk n_{\kv_2}\del{12}{\ko}}{(s+\imag\om{12}{\ko})s} 
    +\frac{4n_{\kv_1}n_{\kv_2}\del{1}{2\ko}}{(s+\imag\om{1}{2\ko})s}
    -\frac{4\nk n_{\kv_2}\del{1}{2\ko}}{(s+\imag\om{1}{2\ko})s}
    +\frac{4\nk n_{\kv_1}\del{1}{2\ko}}{(s+\imag\om{1}{2\ko})s}\nonumber\\
    &\qquad\qquad\qquad\qquad 
    -\frac{4\nk n_{\kv_2}\del{1}{}}{(s+\imag\om{1}{})s}
    +\frac{4n_{\kv_1}n_{\kv_2}\del{}{0}}{(s+\imag\om{}{\ko})s} 
    +\frac{4\nk n_{\kv_1}\del{}{2}}{(s+\imag\om{}{2})s}
    \nonumber\\
    &+
    \frac{2n_{\kv_1}n_{\kv_2}\del{12}{\ko}}{(s-\imag\om{12}{\ko})s}
    -\frac{4\nk n_{\kv_2}\del{12}{\ko}}{(s-\imag\om{12}{\ko})s} 
    +\frac{4n_{\kv_1}n_{\kv_2}\del{1}{2\ko}}{(s-\imag\om{1}{2\ko})s}
    -\frac{4\nk n_{\kv_2}\del{1}{2\ko}}{(s-\imag\om{1}{2\ko})s}
    +\frac{4\nk n_{\kv_1}\del{1}{2\ko}}{(s-\imag\om{1}{2\ko})s}\nonumber\\
    &\qquad\qquad\qquad\qquad 
    -\frac{4\nk n_{\kv_2}\del{1}{}}{(s-\imag\om{1}{})s}
    +\frac{4n_{\kv_1}n_{\kv_2}\del{}{0}}{(s-\imag\om{}{\ko})s} 
    +\frac{4\nk n_{\kv_1}\del{}{2}}{(s-\imag\om{}{2})s}
    \Bigg].
    \label{eq:AfterExpecation}
\ea
Upon collecting like terms and simplifying we obtain
\ba
  \mathcal{L}[\nk](s)-\frac{\nk(0)}{s} \approx &\mu^2\sum_{1,2}\Bigg[
    \del{12}{\ko}\frac{4n_{\kv_1}n_{\kv_2}-8\nk n_{\kv_2}}{s(s^2+(\om{12}{\ko})^2)}
    +\del{1}{2\ko}\frac{8n_{\kv_1}n_{\kv_2}-8\nk n_{\kv_2}+8\nk n_{\kv_1}}{s(s^2+(\om{1}{2\ko})^2)}\nonumber\\
    &\qquad\qquad\qquad\qquad 
    -\frac{8\nk n_{\kv_2}\del{1}{}}{s(s^2+(\om{1}{})^2)}
    +\frac{8n_{\kv_1}n_{\kv_2}\del{}{0}}{s(s^2+(\om{}{\ko})^2)} 
    +\frac{8\nk n_{\kv_1}\del{}{2}}{s(s+(\om{}{2})^2)}\Bigg].
\ea
Now taking the inverse Laplace transform yields
\ba
  \nk(t)-\nk(0) \approx &\mu^2\sum_{1,2}\Bigg[
    \del{12}{\ko}\left(4n_{\kv_1}n_{\kv_2}-8\nk n_{\kv_2}\right)\frac{1-\cos(t\om{12}{\ko})}{(\om{12}{\ko})^2}\nonumber\\
    & \qquad +\del{1}{2\ko}\left(8n_{\kv_1}n_{\kv_2}-8\nk n_{\kv_2}+8\nk n_{\kv_1}\right)\frac{1-\cos(t\om{1}{2\ko})}{(\om{1}{2\ko})^2}\nonumber\\
    &
    -4\nk n_{\kv_2}\del{1}{}\frac{1-\cos(t\om{1}{})}{\om{1}{}}
    +8n_{\kv_1}n_{\kv_2}\del{}{0}\frac{1-\cos(t\om{}{0})}{\om{}{0}}
    +8\nk n_{\kv_1}\del{}{2}\frac{1-\cos(t\om{}{2})}{\om{}{2}}\Bigg].
\ea
This can be more compactly written as
\ba
  \nk(t)-\nk(0) \approx &\mu^2\sum_{1,2}\Bigg[
    \del{12}{\ko}\left(2n_{\kv_1}n_{\kv_2}-4\nk n_{\kv_2}\right)
      g(t,\om{12}{\ko}))\nonumber\\
    & \qquad +\del{1}{2\ko}\left(4n_{\kv_1}n_{\kv_2}-4\nk n_{\kv_2}+4\nk n_{\kv_1}\right)
      g(t,\om{1}{20})\nonumber\\
    &
    -2\nk n_{\kv_2}\del{1}{}g(t,\om{1}{})
    +4n_{\kv_1}n_{\kv_2}\del{}{0}g(t,\om{}{0})
    +4\nk n_{\kv_1}\del{}{2}g(t,\om{}{2})\Bigg],
\ea
where
 \[
   g\left(t,\Omega\right) = \frac{\sin^2(\frac{t\Omega}{2})}{\left(\frac{\Omega}{2}\right)^2}.
 \]
 %
 Now note that 
 \bse
 \ba
  \sum_{1,2}\del{12}{\ko}f(\kv_1,\kv_2,\kv) 
    & = \sum_{\jv_1}\sum_{\jv_2}\delta[\kv_1+\kv_2-\kv]f(\kv_1,\kv_2,\kv)\\
    & = \sum_{\jv_1}f(\kv_1,\kv-\kv_1,\kv)\\
    & = \sum_{\jv_1}\int \delta(\kv_1+\kv_2-\kv)f(\kv_1,\kv_2,\kv)\, d\kv_2,
\ea
\ese
where the convolution sum with a Kronecker delta has been replaced with a convolution integral with the Dirac delta. As a result, we can write 
 \ba
  \nk(t)-\nk(0) \approx &\mu^2\sum_{\jv_1}\int
    \delta(\kv_1+\kv_2-\kv)\left(2n_{\kv_1}n_{\kv_2}-4\nk n_{\kv_2}\right)
      g(t,\om{12}{\ko})\nonumber\\
    & \qquad +\delta(\kv_1-\kv_2-\kv)\left(4n_{\kv_1}n_{\kv_2}-4\nk n_{\kv_2}+4\nk n_{\kv_1}\right)
      g(t,\om{1}{2\ko})\nonumber\\
    &
    -2\nk n_{\kv_2}\delta(\kv_1)g(t,\om{1}{})
    +4n_{\kv_1}n_{\kv_2}\delta(\kv)g(t,\om{}{0})
    +4\nk n_{\kv_1}\delta(\kv_2)g(t,\om{}{2})\,d\kv_2.
\ea
The dependence on $L$ can be exposed in the formulation by unrolling the definition of $\mu$ to give
 \ba
  \nk(t)-\nk(0) \approx &\frac{\epsilon^2}{L^2}\sum_{\jv_1}\int
    \delta(\kv_1+\kv_2-\kv)\left(2n_{\kv_1}n_{\kv_2}-4\nk n_{\kv_2}\right)
      g(t,\om{12}{\ko})\nonumber\\
    & \qquad +\delta(\kv_1-\kv_2-\kv)\left(4n_{\kv_1}n_{\kv_2}-4\nk n_{\kv_2}+4\nk n_{\kv_1}\right)
      g(t,\om{1}{2\ko})\,d\kv_2\nonumber\\
    &
    -2\nk n_{\kv_2}\delta(\kv_1)g(t,\om{1}{})
    +4n_{\kv_1}n_{\kv_2}\delta(\kv)g(t,\om{}{0})
    +4\nk n_{\kv_1}\delta(\kv_2)g(t,\om{}{2})\,d\kv_2.
\ea
Because the k-space increments are given by $\Delta k=2\pi/L$, we now obtain 
 \ba
  \nk(t)-\nk(0) \approx &\frac{\epsilon^2}{(2\pi)^2}\frac{(2\pi)^2}{L^2}\sum_{\jv_1}\int
    \delta(\kv_1+\kv_2-\kv)\left(2n_{\kv_1}n_{\kv_2}-4\nk n_{\kv_2}\right)
      g(t,\om{12}{\ko})\nonumber\\
    & \qquad +\delta(\kv_1-\kv_2-\kv)\left(4n_{\kv_1}n_{\kv_2}-4\nk n_{\kv_2}+4\nk n_{\kv_1}\right)
      g(t,\om{1}{2\ko})\,d\kv_2\nonumber\\
    &
    -2\nk n_{\kv_2}\delta(\kv_1)g(t,\om{1}{})
    +4n_{\kv_1}n_{\kv_2}\delta(\kv)g(t,\om{}{0})
    +4\nk n_{\kv_1}\delta(\kv_2)g(t,\om{}{2})\,d\kv_2,
\ea
and taking the limit as $L\to\infty$ gives
 \ba
  \nk(t)-\nk(0) \approx &\frac{\epsilon^2}{(2\pi)^2}\int\int
    \delta(\kv_1+\kv_2-\kv)\left(2n_{\kv_1}n_{\kv_2}-4\nk n_{\kv_2}\right)
      g(t,\om{12}{\ko})\nonumber\\
    & \qquad +\delta(\kv_1-\kv_2-\kv)\left(4n_{\kv_1}n_{\kv_2}-4\nk n_{\kv_2}+4\nk n_{\kv_1}\right)
      g(t,\om{1}{2\ko})\,d\kv_2\nonumber\\
    &
    -2\nk n_{\kv_2}\delta(\kv_1)g(t,\om{1}{})
    +4n_{\kv_1}n_{\kv_2}\delta(\kv)g(t,\om{}{0})
    +4\nk n_{\kv_1}\delta(\kv_2)g(t,\om{}{2})\,d\kv_2\, d\kv_1.
\ea
Taking $t\to\infty$, and relying on the following 
\ba
  \lim_{t\to\infty}
     = g(t,\Omega)  
     = \lim_{t\to\infty}\frac{\sin^2\left(\frac{t\Omega}{2}\right)}{\left(\frac{\Omega}{2}\right)^2}
     = 2\pi t\delta(\Omega)
\ea
leads to 
 \ba
  \nk(t)-\nk(0) \approx &\frac{t\epsilon^2}{2\pi}\int\int
    \delta(\kv_1+\kv_2-\kv)\left(2n_{\kv_1}n_{\kv_2}-4\nk n_{\kv_2}\right)
      \delta(\om{12}{\ko})\nonumber\\
    & \qquad +\delta(\kv_1-\kv_2-\kv)\left(4n_{\kv_1}n_{\kv_2}-4\nk n_{\kv_2}+4\nk n_{\kv_1}\right)
      \delta(\om{1}{2\ko})\nonumber\\
    &
    -2\nk n_{\kv_2}\delta(\kv_1)\delta(\om{1}{})
    +4n_{\kv_1}n_{\kv_2}\delta(\kv)\delta(\om{}{0})
    +4\nk n_{\kv_1}\delta(\kv_2)\delta(\om{}{2})\,d\kv_2\, d\kv_1.
    \label{eq:tToInfinity}
\ea
Here we see that the ``funny'' matching terms eventually drop out since $\om{j}{}>0$ and $\om{}{j}>0$. Finally taking $\tau = t\epsilon^2$, and considering small $\tau$ yields
\ba
  \p_{\tau}\nk \approx &\frac{1}{\pi}\int\int
    \delta(\kv_1+\kv_2-\kv)\left(n_{\kv_1}n_{\kv_2}-2\nk n_{\kv_2}\right)
      \delta(\omega_{\kv_1}+\omega_{\kv_2}-\omega_{\kv})\nonumber\\
    & \qquad +\delta(\kv_1-\kv_2-\kv)\left(2n_{\kv_1}n_{\kv_2}-2\nk n_{\kv_2}+2\nk n_{\kv_1}\right)
      \delta(\omega_{\kv_1}-\omega_{\kv_2}-\omega_{\kv})\,d\kv_2\, d\kv_1.
      \label{eq:tauToZero}
\ea
Also if symmetry hadn't been used in going from~\eqref{eq:4Suma} to~\eqref{eq:IBP2}, the following more symmetric version results
\ba
  \p_{\tau}\nk \approx &\frac{1}{\pi}\int\int
    \delta(\kv_1+\kv_2-\kv)\left(n_{\kv_1}n_{\kv_2}-\nk n_{\kv_1}-\nk n_{\kv_2}\right)
      \delta(\omega_{\kv_1}+\omega_{\kv_2}-\omega_{\kv})\nonumber\\
    & \qquad +\delta(\kv_1-\kv_2-\kv)\left(2n_{\kv_1}n_{\kv_2}-2\nk n_{\kv_2}+2\nk n_{\kv_1}\right)
      \delta(\omega_{\kv_1}-\omega_{\kv_2}-\omega_{\kv})\,d\kv_2\, d\kv_1.
\ea

\end{document}